\documentclass[preprint,number]{elsarticle}
\journal{Journal of Computational Physics}

\usepackage{graphicx}
\graphicspath{{./figures}}

\usepackage[a4paper,margin=1in]{geometry}
\usepackage[base]{babel}
\usepackage{amsmath}
\usepackage{amsfonts}
\usepackage{amsthm}
\usepackage{subcaption} 
\usepackage{xcolor}
\usepackage{mathtools}
\usepackage{hyperref}
\usepackage{cleveref}
\usepackage{nicefrac}
\usepackage{bm}
\usepackage{algorithm}
\usepackage{algpseudocode}
\usepackage{todonotes} 
\usepackage[final]{showlabels} 
\usepackage{lipsum}
\usepackage{booktabs}
\usepackage{siunitx}

\DeclareMathOperator*{\argmin}{arg\,min}

\newcommand{\R}{\mathbb{R}}

\newcommand{\E}{\mathbb{E}}
\newcommand{\mP}{\mathbb{P}}

\newcommand{\bU}{\bm{U}}

\newcommand{\bY}{\bm{Y}}

\newcommand{\bK}{\bm{K}}
\newcommand{\bI}{\bm{I}}

\newcommand{\bone}{\bm{1}}
\newcommand{\bzero}{\bm{0}}
\newcommand{\bn}{\bm{n}}

\newcommand{\bS}{\bm{S}}

\newcommand{\bL}{\bm{L}}

\newcommand{\lbU}{\overline{\bm{U}}}
\newcommand{\lbY}{\overline{\bm{Y}}}

\newcommand{\Qcal}{\mathcal{Q}}
\newcommand{\Ocal}{\mathcal{O}}
\newcommand{\Pcal}{\mathcal{P}}
\newcommand{\Acal}{\mathcal{A}}
\newcommand{\Tcal}{\mathcal{T}}

\newcommand{\Zcal}{\mathcal{Z}}
\newcommand{\Vcal}{\mathcal{V}}
\newcommand{\Wcal}{\mathcal{W}}

\newcommand{\Xcal}{\mathcal{X}}
\newcommand{\Ncal}{\mathcal{N}}

\newcommand{\delt}{\triangle t}

\newcommand{\hVcal}{\widehat{\mathcal{V}}}
\newcommand{\hXcal}{\widehat{\mathcal{X}}}
\newcommand{\hlXcal}{\widehat{\bar{{\mathcal{X}}}}}
\newcommand{\hbY}{\widehat{\bm{Y}}}
\newcommand{\hbZ}{\widehat{\bm{Z}}}
\newcommand{\hbV}{\widehat{\bm{V}}}
\newcommand{\hM}{\widehat{M}}
\newcommand{\htheta}{\hat{\theta}}

\newcommand{\lXcal}{\overline{\mathcal{X}}}

\newcommand{\lH}{\overline{H}}

\newcommand{\norm}[1]{\| #1 \|}
\newcommand{\trnorm}[1]{{\left\vert\kern-0.25ex\left\vert\kern-0.25ex\left\vert #1 
    \right\vert\kern-0.25ex\right\vert\kern-0.25ex\right\vert}}

\newcommand{\rd}{\mathrm{d}}

\newcommand{\mspan}{\mathrm{span}}

\newcommand{\hU}{\widehat{U}}
\newcommand{\hY}{\widehat{Y}}
\newcommand{\hbU}{\widehat{\bm{U}}}

\newcommand{\mdiv}{\mathrm{div}}

\newcommand{\Goh}{\Gamma^{\nicefrac{1}{2}}}
\newcommand{\Gmoh}{\Gamma^{-\nicefrac{1}{2}}}

\newcommand{\tGamma}{\widetilde{\Gamma}}

\newcommand{\D}{\Delta}

\newcommand{\hP}{\widehat{P}}
\newcommand{\hlP}{\widehat{\bar{P}}}
\newcommand{\ham}{\widehat{m}}

\newcommand{\hlam}{\widehat{\bar{m}}}

\newcommand{\Mat}[2][]{\mathrm{Mat}_{#1}(#2)}

\newcommand{\oh}{\nicefrac{1}{2}}

\newcommand{\Law}{\mathrm{Law}}
\newcommand{\true}{\mathrm{true}}

\newtheorem{remark}{Remark}[section]

\newcommand{\md}{\mathrm{d}}

\newcommand{\tr}{\mathrm{tr}}

\newcommand{\DLR}{\mathrm{DLR}}

\newcommand{\KBP}{\mathrm{KBP}}
\newcommand{\ENKF}{\mathrm{ENKF}}

\newcommand{\kbp}{\textsc{Kbp}}

\newcommand{\enkf}{\textsc{Enkf}}
\newcommand{\venkf}{\textsc{Venkf}}
\newcommand{\denkf}{\textsc{Denkf}}
\newcommand{\senkf}{\textsc{Senkf}}
\newcommand{\dlr}{\textsc{Dlr}}

\newcommand{\st}{\,|\,}

\newcommand{\ortho}{\texttt{ortho}}

\newcommand{\moda}[1]{{\color{black} #1}}
\newcommand{\modb}[1]{{\color{black} #1}}
\newcommand{\ra}[1]{{\color{black} #1}}
\newcommand{\rb}[1]{{\color{black} #1}}
\newcommand{\rc}[1]{{\color{black} #1}}

\DeclareSIUnit\mmHg{mmHg}

\floatname{algorithm}{Algorithm}

\algnewcommand{\algorithmicendif}{\textbf{end}}
\algblockdefx[IF]{If}{EndIf}[1]{\algorithmicif\ #1\ \algorithmicthen}{\algorithmicendif}
\algblockdefx[For]{for}{EndFor}[1]{\algorithmicif\ #1\ \algorithmicthen}{\algorithmicendif}

\title{Dynamical Low-Rank Ensemble Kalman Filter \\
	for State/Parameter Estimation}
  \author[1,2]{Thomas Trigo Trindade\corref{cor1}\fnref{fn1}} 
	\author[1]{Sébastien Riffaud}
	\author[1]{Fabio Nobile}

  \affiliation[1]{organization={Institute of Mathematics, École Polytechnique Fédérale de Lausanne (EPFL)},
    postcode={1015},
    postcodesep={},
    city={Lausanne},
  country={Switzerland}}

  \affiliation[2]{organization={Computer, Electrical and Mathematical Sciences and Engineering Division, King Abdullah University of Science and Technology (KAUST)},
    postcode={23955-6900},
    postcodesep={},
    city={Thuwal},
  country={Saudi Arabia}}
  \cortext[cor1]{Corresponding author}
  \fntext[fn1]{T. Trigo Trindade carried out this work while at EPFL.}

\begin{document}

\begin{abstract}
We propose a Dynamical Low-Rank Ensemble Kalman Filter (\dlr-\enkf) for efficient joint state-parameter estimation in high-dimensional dynamical systems. 
The method extends the \dlr-\enkf\ formulation of~\cite{DynamicalLoTrigo2025} to the \modb{augmented state-parameter} framework, \modb{tracking the} filtering density within a dynamically evolving \moda{low-dimensional} subspace. 
Key developments include a time-integration strategy that combines the Basis Update \& Galerkin scheme with forecast/analysis discretisation, and a DEIM-based hyper-reduction technique for efficient evaluation of nonlinear terms. We demonstrate the effectiveness, robustness, and computational advantages of the proposed approach on benchmark problems. 
The results highlight the potential of dynamically evolving reduced bases to achieve accurate filtering and parameter estimation at reduced computational cost.	
\end{abstract}


\maketitle

\section{Introduction}

In Uncertainty Quantification, state and/or parameter identification tasks remain challenging due to their high computational cost. 
This work focuses on the Bayesian data assimilation framework, where a signal $\Xcal_t^{\true}$ evolves according to model dynamics that depend on an unknown parameter $\theta^{\true}$. 
The signal is observed through an observation process $\Zcal_t$, defined as a (linear) functional of the signal and perturbed by additive noise. 
The Bayesian problem \modb{consists in characterising} the posterior density of $(\Xcal_t, \theta_t)$ conditional on the observations. 
In filtering, this posterior is updated sequentially using incoming data, together with the current distributions of the state and parameters. 
Since the exact filtering density is not available in closed form, approximate methods are required. 
Common approaches include particle-based techniques such as Sequential Monte Carlo (particle filters) and ensemble-based methods such as the Ensemble Kalman Filter (\enkf). 
However, when the state dynamics \moda{is} high-dimensional -- for instance, arising from the fine discretisation of a time-dependent PDE -- the computational cost of propagating the dynamical system can become prohibitive for a large number of particles.

\rc{This work focuses on settings in which the filtering distribution concentrates around a low-dimensional (time-dependent) subspace; in this regime, the \enkf\ is expected to track the filtering density reasonably well even with a moderate number of particles. 
Using a larger ensemble nevertheless has advantages, as it reduces the Monte Carlo error and improves accuracy.
Importantly, the \enkf\ does not explicitly exploit this low-rank structure when it is present: the cost of propagating particles still scales with the full state-space dimension, rather than the reduced dimension of the subspace.}
The work~\cite{DynamicalLoTrigo2025} addresses this issue by developing a framework that evolves a \moda{low-dimensional} filtering subspace using the Dynamical Low-Rank (\dlr) approximation method. 
Originally introduced for matrix differential equations~\cite{DynamicalLowRKoch2007}, \dlr\ approximations propagate a low-rank, SVD-like representation of the solution through tailored time-integration schemes~\cite{AnUnconventionCeruti2022,AProjectorSplLubich2014}, thereby avoiding repeated costly SVD truncations at each step. 
The method has also been applied in Uncertainty Quantification for random PDEs~\cite{StabilityPropeKazash2021,ErrorAnalysisMushar2015,DynamicallyOrtSapsis2009} and SDEs~\cite{DynamicalLowZoccolan2023}.
The objective of this work is to design a practical \dlr-\enkf\ algorithm for efficient joint state-parameter estimation \moda{in the context of time-continuous observations.}
To this end, we extend the formulation of~\cite{DynamicalLoTrigo2025} \moda{to nonlinear dynamical systems} and further investigate several questions left open in that work, which are crucial for robust and efficient propagation of the low-rank system. 
In particular, we examine the development of suitable time-integration schemes, the use of hyper-reduction techniques for efficient treatment of nonlinear terms \rb{and the use of rank adaptive algorithms to properly follow the evolution of the system.} 
These aspects will be central to the approach presented in this paper.

The main contributions of this paper are \moda{the following}. 
First, we adapt the \dlr-\enkf\ formulation of~\cite{DynamicalLoTrigo2025} to \moda{nonlinear dynamics and} the augmented-state \enkf\ framework in order to perform joint state-parameter estimation~\cite{AnEnsembleAnderson2001,StateAndParamNusken2019}. 
Second, we propose a time-integration strategy that combines the Basis Update \& Galerkin scheme~\cite{ARankAdaptiveCeruti2022} with the forecast/analysis-based discretisation \rb{and discuss rank adaptive strategy}. 
Third, we introduce a novel DEIM-based approach to accelerate the evaluation of nonlinear terms. 
Finally, we demonstrate the effectiveness, strengths, and limitations of the method through applications to representative model problems.


\rb{This work lies within the broader research area that studies the integration of reduced-order models (ROMs) into state and/or parameter identification problems. Closely related to our problem setting, \cite{EfficientStatePagani2017} combines the POD Reduced Basis technique with an augmented-state \enkf~formulation, supplemented by a Reduced Error Model to mitigate the biases introduced by the model order reduction. Fixed-basis ROM-EnKF approaches can be very efficient once a suitable reduced basis has been constructed. However, fixed reduced basis methods may be less effective when the solution does not exhibit a low-rank profile globally in space and time. This limitation is particularly evident in problems with traveling wave fronts: although the solution may admit a low-rank representation at each fixed time, its global-in-time evolution may not be well captured by a small fixed spatial basis, whereas a dynamically evolving basis can adapt to the relevant solution subspace as it changes over time. Moreover, in filtering problems, the assimilation procedure may progressively concentrate the posterior distribution, and a rank-adaptive formulation can naturally exploit this effect by decreasing the rank as the subspace explored by the filtering ensemble shrinks.

A further difficulty of fixed-basis approaches is the offline training stage. A fixed-basis ROM must be trained before the observations are assimilated, and therefore the training set must anticipate the part of the solution manifold that will be visited by the filter. This is particularly delicate in state/parameter estimation, since the region of the state/parameter space explored by the posterior is not known a priori. Unless prior information is available, this may require sampling a large parameter space with the full order model, which can become computationally prohibitive when the number of parameters increases. In contrast, our approach does not require an offline training phase and employs a reduced basis that evolves dynamically in time, allowing the approximation space to adapt online to the evolving filtering density.}

\modb{
Several other works have explored the integration of \emph{dynamic} reduced basis techniques in data assimilation, specifically in the context of discrete time observations.}
\moda{
	For instance,~\cite{AnAdaptiveSilva2025} proposes an adaptive technique to update the reduced basis in the online phase in the context of multi-fidelity/level Ensemble Kalman filters. 
	More closely related to our setting, a number of works specifically employ dynamical low-rank techniques.
	One example is the blended particle filters framework~\cite{BlendedPaMajda2014}, which builds on the (modified) quasilinear Gaussian-DO method~\cite{BlendedRedSapsis2013,AStatisticaMajda2013} to perform data assimilation of chaotic large-scale systems (with an emphasis on handling non-Gaussian features). 
	Finally, the authors in~\cite{DataAssimilatiPt2Sonder2013,DataAssimilatiSonder2013} use DLR techniques to accelerate the forecast step and fit the forecasted particles to a Gaussian mixture model for the assimilation step.}

\modb{We also mention related work in the broader context} of Bayesian parameter identification for parametrised, \modb{time-dependent PDEs:} \cite{LowRankRiffaud2024} develops a low-rank MCMC method that exploits dynamically evolving low-rank bases to reduce computational costs while maintaining accuracy. Alternatively, \cite{AReducedBasisSilva2023} proposes surrogate models for Ensemble Kalman Inversion (EKI), combined with ROM bias correction. While the former relies on an MCMC method based on low-rank approximations and the latter follows a ROM-based Ensemble Kalman Inversion approach, both methods do not assimilate observations sequentially and are designed for smoothing problems. In contrast, our method assimilates data on-the-fly, which makes it naturally suited for filtering applications.



The paper is organised as follows. 
\Cref{sec:problem-setting} introduces the concepts and notation for the filtering problem. 
In \Cref{sec:dlr-state-param-estimation}, we present the \dlr-\enkf\ equations for parameter identification, together with the proposed time-integration scheme and its DEIM-based acceleration. 
In \Cref{sec:numexp}, we assess the efficiency of the method on two parameter identification problems, a Fisher-KPP model and a blood flow model.

\section{Problem setting} \label{sec:problem-setting}

Let $(\Omega, \Acal , \mP; (\Acal_t)_{t \geq 0})$ be a complete right-continuous filtered probability space. 
Denote by $\Mat{m \times n}$ the set of real-valued $m\times n$ matrices, and $\Mat[+]{m}$, resp. $\Mat[0]{m}$, the set of symmetric positive resp. semi-positive-definite $m \times m$ matrices.
For $A \in \Mat{m \times n}$, denote $\norm{A}_F^2 = \tr(A^{\top} A)$ the Frobenius norm, and by $\Tcal_R(A) \in \mathrm{argmin} \{ \norm{A - B}_F, \mathrm{rank}(B) \leq R \}$ a best rank-$R$ approximation of $A$, which is unique if the $R$-th singular value $\sigma_{R}(A)$ \moda{is strictly larger than the following one} $\sigma_{R+1}(A)$. 

We consider a system whose state $\Xcal^{\true}_t \in \R^d$ (hereafter called the signal) evolves according to some parameter-dependent \textit{true dynamics}
\begin{equation} \label{eqn:true-dynamics}
	\rd \Xcal^{\true}_t = f(\Xcal^{\true}_t, \theta^{\true}) \rd t, \quad \moda{t \geq 0},
\end{equation}
with~$f : \R^d \to \R^d$. 
The dynamics depend on the \textit{true parameter} vector $\theta^{\true} \in \R^{n_\theta}$, which will be assumed to have a small to moderate dimension. 
The signal is \moda{partially} observed via the process $\Zcal_t \in \R^k$ \moda{where $k < d$}, driven by
\begin{equation} \label{eqn:obs-process}
	\rd \Zcal_t = H \Xcal^{\true}_t \rd t + \Goh \rd \widetilde{\Vcal}_t. 
\end{equation}
with $H \in \Mat{k \times d}$ a \textit{linear observation operator} and $\Gamma \in \Mat[+]{k}$, and where $\widetilde{\Vcal}_t$ is a $k$-dimensional \moda{$\Acal_t$-adapted} Brownian motion.
We set $\Zcal_0 = \bm{0} \in \R^k$ and denote $\Acal_{\Zcal_t} = \sigma(\Zcal_s, 0 \leq s \leq t)$ the sub-$\sigma$-algebra generated by the observations up to time $t$, in other words the history of observations, \moda{whereas $f_t^{\star} = f - \E[f \st \Acal_{\Zcal_t}]$ denotes the zero conditional mean part of $f$.}
For $f \in L^1(\Omega)$, the conditional expectation taken with respect to $\Acal_{\Zcal_t}$ shall be denoted by $\E[f \,|\,\Acal_{\Zcal_t}]$. 

\moda{In this work, we consider the situation where $\Xcal_0^{\true}$ and $\theta^{\true}$ are unknown, and use a Bayesian approach to infer the true state and parameter from the data. 
To this end, we model $\Xcal_0^{\true}$ and $\theta^{\true}$ as random variables, with prior distribution $(\Xcal_0, \theta_0) \sim \pi_{\mathrm{pr}}$. 
Then, we aim to efficiently approximate} the state distribution \modb{at any time $t \geq 0$} conditional on the data $\Law(\Xcal_t^{\true} \st \Acal_{\Zcal_t})$ as well as the parameter distribution $\Law(\theta^{\true} \st \Acal_{\Zcal_t})$. 
We propose to do this by following an Ensemble Kalman Filter approach~\cite{OnTheMathematBishop2023}. 
The setup is similar to the one of Ensemble Kalman Inversion~\cite{AnalysisoftheEnsSchillings2017}, with the difference that the observation process varies in time and is assimilated on the fly. 
We refer to~\cite{StateAndParamNusken2019} for a discussion \modb{on} the Ensemble Kalman filter for state/parameter \modb{estimation} in the context of time-varying observations. 

The approach proceeds by (a) constructing an augmented state process $(\Xcal^{\true}_t, \theta^{\true}_t) \equiv \overline{\Xcal}_t \in \R^{d + n_{\theta}}$, and (b) applying the Kalman equations to the augmented state dynamics, \moda{which} verifies
\begin{equation} \label{eqn:aug-true-dynamics}
\rd \overline{\Xcal}^{\true}_t = \bar{f}(\overline{\Xcal}^{\true}_t) \rd t, \quad \moda{t \geq 0, \quad \overline{\Xcal_0}^{\true} = (\Xcal_0^{\true}, \theta^{\true})}
\end{equation}
where $\bar{f}(\overline{\Xcal}^{\true}_t) = [f(\Xcal^{\true}_t, \theta^{\true}_t), \bzero_{n_\theta}]^{\top}$,
and \moda{observation process}
\begin{equation} \label{eqn:aug-obs-process}
	\rd \Zcal_t = \overline{H}\, \overline{\Xcal}^{\true}_t \rd t + \Goh \rd \widetilde{\Vcal}_t, 
\end{equation}
where $\overline{H} = [H, \bzero_{k \times n_{\theta}}]$. 
Note that in the augmented dynamics, $\md \theta^{\true}_t = 0$ with $\theta_0^{\true}  = \theta^{\true}$, and therefore the augmented dynamics and observations is equivalent to the original dynamics and observations. 
\moda{As $\Xcal_0^{\true}$ and $\theta^{\true}_0$ are random (and distributed according to $\pi_{\mathrm{pr}}$), the Kalman filter will induce a non-trivial dynamics on the parameter $\theta_t$, which hopefully will approach $\theta^{\true}$ if the observation process is sufficiently informative.}

\subsection{Augmented State Kalman filtering}

Ensemble Kalman filters are equal-weight interacting particle systems that evolve an empirical measure by combining the model dynamics with a corrective term depending on the observation process.
Even when used with a moderate number of particles, these methods are less susceptible \modb{than standard particle filters to ensemble collapse~\cite{InverseProbSanz2023}}, making these methods especially popular in large-scale data assimilation applications such as geosciences. 
The state/parameter identification equations are derived by applying the \enkf~equations to the augmented state system~\eqref{eqn:aug-true-dynamics} and~\eqref{eqn:aug-obs-process}. 
We refer to e.g.~\cite{StateAndParamNusken2019} for the construction in the filtering context (with model noise), or more broadly~\cite{EnsembleKIglesias2013} for the derivation in \moda{the} related Ensemble Kalman Inversion (EKI) context. 
We shall focus on two variants of the \enkf, which, following the nomenclature in~\cite{OnTheMathematBishop2023}, are referred to as the Vanilla Ensemble Kalman filter (\venkf) and the Deterministic Ensemble Kalman filter (\denkf).
Denote $\hlXcal_t = [\hlXcal_t^{(1)}, \ldots, \hlXcal_t^{({\rc{N}})}]$ a set of ${\rc{N}}$ augmented particles, and let
\begin{equation*}
	\begin{bmatrix}
		\ham_t^x \\
		\ham_t^{\theta}
	\end{bmatrix} 
	= 
	\hlam_t 
	=
	\hlam_t(\hlXcal_t)
	= 
	\E_{\rc{N}}[\hlXcal_t] \moda{\coloneqq} \frac{1}{{\rc{N}}} \sum_{i=1}^{\rc{N}} \hlXcal_t^{(i)}, 
	\quad 
        \hlXcal_t^{\star} = \hlXcal_t - \E_{\rc{N}}[\hlXcal_t],
\end{equation*}
\moda{be} the sample mean and the sample zero-mean \modb{ensemble, respectively}, and 
\begin{equation*}
        \hlP_t 
	= 
	\hlP_t(\hlXcal_t)
	= \frac{{\rc{N}}}{{\rc{N}}-1} \E_{\rc{N}} [\hlXcal^{\star}_t (\hlXcal^{\star}_t)^{\top}]
	= 
	\begin{bmatrix}
		\hP^{xx}_t & \hP^{x \theta}_{\modb{t}} \\
		\hP^{\theta x}_t & \hP^{\theta \theta}_t
	\end{bmatrix}
\end{equation*}
the sample covariance. 
Furthermore, \moda{$ \hVcal_t = \{\hVcal_t^{(p)}\}_{p=1}^{\rc{N}}$ denotes an i.i.d sample of} copies of $\Vcal_t$.
The Ensemble Kalman filters are then written as
\begin{align}
	\rd \hXcal^{(p)}_t &= 
	f(\hXcal^{(p)}_t, \htheta_t^{(p)}) \rd t \nonumber 
	\\
			   &\quad + \hP_t^{xx} H^{\top} \Gamma^{-1} \left( \rd \Zcal_t - \frac{1 + \kappa}{2} H \hXcal^{(p)}_t \rd t  - \kappa \Goh \rd \hVcal_t^{(p)} - \frac{1 - \kappa}{2} H \ham_t^x \rd t \right), \label{eqn:enkf-dynamics-continuous-state} \\
	\rd \htheta_t^{(p)} &= \hP^{\theta x}_t H^{\top} \Gamma^{-1} \left( \rd \Zcal_t - \frac{1 + \kappa}{2} H \hXcal^{(p)}_t \rd t  - \kappa \Goh \rd \hVcal_t^{(p)} - \frac{1 - \kappa}{2} H \ham_t^x \rd t  \right) , \label{eqn:enkf-dynamics-continuous-param}
\end{align}
\moda{where} $\kappa = 1$ \moda{corresponds to} the \venkf~and $\kappa = 0$ \moda{to} the \denkf. 
Observe that the \venkf~\moda{mimicks} the dynamics and observation process by adding a synthetic noise to the observation of the state; on the other hand, the \denkf~has a more consensus-based interpretation, taking the \moda{difference} between a given particle and the sample average. 
While computationally expensive, simulating \moda{the \enkf\ system~\eqref{eqn:enkf-dynamics-continuous-state} and~\eqref{eqn:enkf-dynamics-continuous-param}} with a large number of particles is of interest for several reasons. 
In the asymptotic limit, the filter may converge to a given mean-field limit (if it exists); in particular, when applied to linear-affine dynamics, the mean-field limit of the Ensemble Kalman filter coincides with the true filtering distribution~\cite{OnTheStabilitDelMo2018} . 
This ceases to be true when applying the \enkf~equation to nonlinear dynamics, and the \enkf~may then rather be understood as a dynamic \moda{state/parameter} estimation method.
We refer to~\cite{OnTheMathematBishop2023} and references therein for a more detailed discussion on this topic.  
In any case, a higher number of particles is expected to improve the signal-tracking abilities of the filter.

\begin{remark}
	\moda{Anticipating the \dlr-\enkf\ construction, the mean-field version of~\cref{eqn:enkf-dynamics-continuous-state,eqn:enkf-dynamics-continuous-param} will be of interest.
		Proceeding formally, it is assumed that the empirical measure $\hat{\bar{\mu}}_t^{{\rc{N}}} = {\rc{N}}^{-1}\sum_{i=1}^{\rc{N}} \delta_{\hlXcal_t^{(i)}}$ converges (in some sense) to the mean-field distribution $\bar{\mu}_t$, and denote $\lXcal_t \sim \bar{\mu}_t$ its associated mean-field process. 
		Then, the sample covariance $\hP_t$ converges to $P_t = \E[\lXcal_t \lXcal_t^{\top}]$ (in some sense), and gives rise to the following (McKean-Vlasov) equations
\begin{align}
	\rd \Xcal_t &= 
	f(\Xcal_t, \theta_t) \rd t \nonumber 
	\\
		    &\quad + P_t^{xx} H^{\top} \Gamma^{-1} \left( \rd \Zcal_t - \frac{1 + \kappa}{2} H \Xcal_t \rd t  - \kappa \Goh \rd \Vcal_t - \frac{1 - \kappa}{2} H \E[\Xcal_t \st \Acal_{\Zcal_t}] \rd t \right), \label{eqn:kbp-dynamics-continuous-state} \\
	\rd \theta_t &= P^{\theta x}_t H^{\top} \Gamma^{-1} \left( \rd \Zcal_t - \frac{1 + \kappa}{2} H \Xcal_t \rd t  - \kappa \Goh \rd \Vcal_t - \frac{1 - \kappa}{2} H \E[\Xcal_t \st \Acal_{\Zcal_t}] \rd t  \right). \label{eqn:kbp-dynamics-continuous-param}
\end{align}
In the case of linear-affine dynamics, \cref{eqn:kbp-dynamics-continuous-state} is known as a Kalman-Bucy process~\cite{OnTheMathematBishop2023}. 
	}
\end{remark}

\subsection{Time-discretisation schemes} \label{sec:time-discretisation-schemes}

The numerical discretisation \moda{of~\cref{eqn:enkf-dynamics-continuous-state,eqn:enkf-dynamics-continuous-param}} may be performed in various fashions. 
A straightforward approach is to resort to standard SDE discretisations (Euler-Maruyama scheme or higher order methods~\cite{NumericalSolutKloeden1992}).
Note that a special care must be taken, as \moda{explicit discretisations as} the Euler-Maruyama method may not be stable owing to the nonlinear (non-uniformly Lipschitz) terms~\cite{StrongAnHutzenthaler2011}. 
Furthermore, the above system can become stiff in the case of small observation noise ($\norm{\Gamma}_F \ll 1$), adding a potential source of numerical instability. 

Alternatively, the update may mimic the forecast/analysis procedure typically encountered in the discrete time data assimilation setting. 
In essence, it inverts the steps taken in~\cite{WellPosedKelly2014} \moda{to recover the continuous-time formulation starting from the discrete setting}, and can be viewed as a semi-implicit splitting of the dynamics.
It is in particular known that \moda{for a one-dimensional toy model mimicking the Kalman-Bucy filter, the scheme is strongly} convergent~\cite{AStronglyCBlomker2018}. 
The present work focuses on this latter approach, and we outline the steps of that development here.

On an interval $[t_0, t_1]$, the \modb{\emph{forecast}} step consists in evolving each particle \modb{following} the model dynamics
\begin{equation} \label{eqn:forecast-step}
	\begin{aligned}
		\md \hXcal_t^{f,(p)} &= f(\hXcal_t^{f,(p)}, \htheta_t^{f,(p)}) \md t 
				     && \hXcal_{t_0}^{f,(p)} = \hXcal_{t_0}^{(p)} \\
		\md \htheta^{f,(p)}_t &= 0, && \htheta_{t_0}^{f,(p)} = \htheta_{t_0}^{(p)}
	\end{aligned}
\end{equation}
via e.g. an Euler step, and then \modb{\emph{analyse}} the particles by subjecting them to the corrective dynamics
\begin{equation} \label{eqn:analysis-step}
	\md 
	\begin{bmatrix}
		\hXcal_t^{a,(p)} \\
		\htheta_t^{a,(p)}
	\end{bmatrix}
	= 
	\begin{bmatrix}
		\hP^{xx}_t \\
		\hP^{\theta x}_t
	\end{bmatrix}
	H^{\top} \Gamma^{-1} \left( \md \Zcal_t - \frac{1 + \kappa}{2} H \hXcal_t^{a,(p)} \md t -  \kappa \Goh \md \hVcal^{(p)}_t - \frac{1 - \kappa}{2} H\ham_t^{x} \md t \right)
\end{equation}
with $\hXcal_{t_0}^{a,(p)} = \hXcal_{t_1}^{f,(p)}$ and $\htheta_{t_0}^{a,(p)} = \htheta_{t_1}^{f,(p)} = \htheta_{t_0}^{f,(p)}$. 
This formulation highlights how the update in that substep does not alter the subspace \modb{spanned by} the ensemble of particles, \moda{because of the left-multiplication of the rhs in~\cref{eqn:analysis-step} by $\hP_t^{xx}$, which keeps the dynamics of $\Xcal_t^{a,(p)}$ in $\mspan\{\hXcal_t^{f,(p)}, p = 1, \ldots {\rc{N}}\}$.}
Rewriting the system in augmented state yields
\begin{equation*}
	\md \hlXcal_t^{a,(p)} = \hlP_t \lH^{\top} \Gamma^{-1} \left( \md \Zcal_t - \frac{1 + \kappa}{2} \lH \hlXcal_t^{a,(p)} \md t -  \kappa \Goh \md \hVcal^{(p)}_t - \frac{1 - \kappa}{2} \lH \hlam_t \md t \right).
\end{equation*}
Using a semi-implicit discretisation yields
\begin{multline*}
	(\bI_{d + n_{\theta}} + \delt \moda{\hlP^f_{t_1}} \moda{\lH^{\top}} \Gamma^{-1} \moda{\lH} ) \hlXcal_{t_1}^{a,(p)} 
	\\
	= \moda{\hlXcal_{t_1}^{f,(p)}} + \moda{\hlP^f_{t_1}} \lH^{\top} \Gamma^{-1} \left( \D \Zcal_{t_1} - \kappa \Goh \D \moda{\hVcal^{(p)}_{t_1}} - \frac{1- \kappa}{2} \lH \left(\hlam_{t_1}^f - \hlXcal_{t_1}^{f,(p)} \right) \delt \right),
\end{multline*}
where $\hlam^f_{t_1} = \hlam(\hlXcal_{t_1}^f)$ and $\moda{\hlP^f_{t_1}} = \hlP(\hlXcal_{t_1}^f)$. 
Assuming $\moda{\hlP^f_{t_1}}$ invertible, \moda{we now} left-multiply by $(\moda{\hP^f_{t_1}})^{-1}$ \moda{and use the} Woodbury identity \moda{to rewrite}
\begin{align*}
	\moda{\left[(\moda{\hlP^f_{t_1}} )^{-1} ( \bI_{d + n_{\theta}} +  \delt \moda{\hlP^f_{t_1}}  \lH^{\top} \Gamma \lH ) \right]^{-1}}
	&= \left[(\moda{\hlP^f_{t_1}} )^{-1} + \lH^{\top} (\Gamma / \delt)^{-1} \lH \right]^{-1} 
	\\
										 &= \moda{\hlP^f_{t_1}} - \moda{\hlP^f_{t_1}} \lH^{\top} \left( \Gamma / \delt + \lH \moda{\hlP^f_{t_1}} \lH^{\top} \right)^{-1} \lH \moda{\hlP^f_{t_1}} \\
										 &= \left( \bI_{d+n_{\theta}} - \delt K_{t_1} \lH \right) \moda{\hlP^f_{t_1}},
\end{align*}
where 
\begin{equation*}
	K_{t_1} = \moda{\hlP^f_{t_1}} \lH^{\top} G_{t_1}^{-1} 
	=
	\begin{bmatrix}
		\hP_{t_1}^{xx,f} H^{\top} G_{t_1}^{-1}  \\
		\hP_{t_1}^{\theta x,f} H^{\top} G_{t_1}^{-1}
	\end{bmatrix}
	= 
	\begin{bmatrix}
		K_{t_1}^{xx} \\
		K_{t_1}^{\theta x}
	\end{bmatrix}
	,
	  \quad  
	G_{t_1} = \Gamma + \delt \lH \moda{\hlP^f_{t_1}} \lH^{\top}.
\end{equation*}
It is \moda{also} possible to show \moda{by direct calculation} that
\begin{equation*}
	(\bI_{d+n_{\theta}} - \delt K_{t_1} \lH ) \moda{\hlP^f_{t_1}} \lH^{\top} \Gamma^{-1} = K_{t_1}, 
\end{equation*}
\moda{by} which the analysis step \moda{can be rewritten as} 
\begin{equation} \label{eqn:analysis-interm}
	\hlXcal_{t_1}^{a,(p)} = (\bI_{d+n_{\theta}} - \delt K_{t_1} \lH ) \hlXcal_{t_1}^{f,(p)} + K_{t_1} \left( \D \Zcal_{t_1} - \kappa \Goh \D \moda{\hVcal^{(p)}_{t_1}} - \frac{1- \kappa}{2} \lH\left(\hlam_{t_1}^f - \hlXcal_{t_1}^{f,(p)} \right) \delt \right), 
\end{equation}
or equivalently, introducing $\eta = 1$, 
\begin{equation} \label{eqn:analysis}
	\begin{bmatrix}
		\hXcal_{t_1}^{a,(p)} \\
		\htheta_{t_1}^{a,(p)}
	\end{bmatrix}
	=  
	\begin{bmatrix}
		\hXcal_{t_1}^{f,(p)} \\
		\htheta_{t_1}^{f,(p)}
	\end{bmatrix}
	+  
	\begin{bmatrix}
		K_{t_1}^{xx} \\
		K_{t_1}^{\theta x}
	\end{bmatrix}
	\left( \D \Zcal_{t_1} -  \frac{1 + \kappa}{2} H \hXcal_{t_1}^{f,(p)} \delt - \eta \kappa \Goh \D \moda{\hVcal^{(p)}_{t_1}} - \eta \frac{1- \kappa}{2} H \ham_{t_1}^{x,f} \moda{\delt} \right). 
\end{equation}
Following~\cite[Lemma 8.5]{DataAssimilLaw2015}, under the \moda{verifiable} assumption of \moda{invertibility of} $\Gamma$, it is possible to extend the analysis update formula to the case of singular $\hlP_t$. 

We will consider the following parameter combinations in our numerical experiments.
\begin{enumerate}
	\item $\eta = 1, \kappa = 1$ (\venkf) 
	\item $\eta = 1, \kappa = 0$  (\denkf)
	\item $\eta = 0, \kappa = 1$ (\senkf) 
\end{enumerate}
The \venkf~discretisation corresponds to the usual forecast/analysis dynamics, and the same discretisation scheme is used to obtain the \denkf~update.
While the \senkf~(``Sequential optimiser'' \enkf) has a straightforward interpretation of computing the update using the observation-forecast mismatch of each particle, it is known~\cite{EnsembleDataAWhitak2002} that it may underestimate analysis-error covariances, possibly causing filter divergence. 
On the other hand, this has an interpretation as a ``sequential optimisation'' which, for each particle, seeks the minimum of a weighted functional combining the deviation-from-observation mismatch and the deviation-from-forecast mismatch~\cite{InverseProbSanz2023}. 
This procedure may be relevant in the case where the filtering distribution is not well-approximated by Gaussians.

Several Ensemble Square Root Kalman filters (including the Ensemble Adjustment Kalman Filter, the Ensemble Transform Kalman Filter, the unperturbed \enkf) are known to converge to the continuous time \denkf~\cite{OnTheContLange2019} in the limit of \moda{continuous} observations (and suitably noise scalings). 
Adapting these filters to the low-rank context and studying their efficiency is of interest, but is beyond the scope of this work, and we shall not discuss this point further.

\subsubsection{Connection to Discrete Observations}

The setting devised in the above section assumes the observations to be given by a process, and is well-defined in the limit of $\delt \rightarrow 0$.
In practice however, it can be of interest to recast the problem in a setting of discrete observations at each time step; this is in particular natural when the data stream consists in direct observations of the state rather than of stochastic increments. 
To that end, we detail how, for a fixed time-discretisation, the discretised analysis step can be understood as the discrete time Kalman-gain assimilation step, with a rescaled observation noise.
For a fixed $\delt > 0$, define
\begin{equation*}
	y_{\moda{t_1}} \coloneq \frac{\D \Zcal_{\moda{t_1}}}{\delt} = H \Xcal_{\moda{t_1}}^{\true} + \left( \frac{\Gamma}{\delt} \right)^{\oh} \xi_{\moda{t_1}}, \quad \xi_{\moda{t_1}} \sim \Ncal(0, I_k).
\end{equation*}
Therefore, defining the rescaled observation noise level $\tGamma = (\Gamma / \delt)$,~\cref{eqn:analysis} is equivalent to
\begin{equation} \label{eqn:analysis-disc}
	\begin{bmatrix}
		\hXcal_{t_1}^{a,(p)} \\
		\htheta_{t_1}^{a,(p)}
	\end{bmatrix}
	=  
	\begin{bmatrix}
		\hXcal_{t_1}^{f,(p)} \\
		\htheta_{t_1}^{f,(p)}
	\end{bmatrix}
	+  
	\begin{bmatrix}
		\widetilde{K}_{t_1}^{xx} \\
		\widetilde{K}_{t_1}^{\theta x}
	\end{bmatrix}
	\left( y_{\modb{t_1}} - \frac{1 + \kappa}{2} H \hXcal_{t_1}^{f,(p)} - \eta \kappa \moda{\tGamma^{\nicefrac{1}{2}} \hat{\xi}_{t_1}^{(p)}} - \eta \frac{1- \kappa}{2} H \ham_{t_1}^{x,f}  \right),
\end{equation}
with $\eta = 1$ and 
\begin{equation*}
  \widetilde{K}_{t_1} = \hlP^f_{t_1} \lH^{\top} \widetilde{G}_{t_1}^{-1}
	=
	\begin{bmatrix}
    \hP_{t_1}^{xx,f} H^{\top} \widetilde{G}_{t_1}^{-1}  \\
    \hP_{t_1}^{\theta x,f} H^{\top} \widetilde{G}_{t_1}^{-1}
	\end{bmatrix}
	= 
	\begin{bmatrix}
		K_{t_1}^{xx} \\
		K_{t_1}^{\theta x}
	\end{bmatrix}
	, 
	\quad
  \widetilde{G}_{t_1} = \tGamma + \lH \hlP^f_{t_1} \lH^{\top}.
\end{equation*}
Those correspond to the Ensemble Kalman update formulae in the discrete time setting.



We conclude this section by noting that the cost of simulating ${\rc{N}}$ augmented state particles scales like $\Ocal(d{\rc{N}})$ (assuming the parameter dimension $n_{\theta} \ll d$). 
The computational cost may become prohibitively expensive when the state space dimension $d$ is large, and offset only by choosing a small to moderate number of particles ${\rc{N}}$.  
Unlike particle filters which often require a prohibitively large number of particles to prevent particle collapse, the Ensemble Kalman as introduced in~\cite{SequentialDataEvense1994} and its many variants are popular data assimilation tools in high dimensional problems, as empirically they perform satisfactorily well under a small to moderate ensemble size.
Note that both methods are still subject to catastrophic filter divergence and may require some stabilisation mechanism -- that issue will not be addressed further in this work. 
Nevertheless, the possibility of increasing the number of particles remains of interest, as doing so can cause the ensemble to converge to a mean-field limit (if it exists) and hence subject to smaller sample variance. In the following section, we propose to perform a model order reduction by way of Dynamical Low-Rank approximations adapted to the Data Assimilation context, \modb{which will allow us to afford evolving larger ensembles thanks to the reduced computational cost of the model.}

\section{Dynamical Low-Rank State/Parameter estimation} \label{sec:dlr-state-param-estimation}

This section firstly establishes the \dlr-\enkf\ system used to efficiently propagate an ensemble of low-rank state particles alongside the ensemble of parameters.
The small to moderate dimension of the parameter implies that no relevant gains should be expected when applying model order reduction techniques to the parameter ensemble; therefore, the \dlr\ technique is only applied to the state particles. 
\moda{On the other hand we expect that the uncertainty on the parameters is reduced over time thanks to the continuous incoming stream of information.
Hence, the state dynamics will also be less and less uncertain, justifying the use of dynamical low-rank methods.} 
Additionally, we introduce a time-integration scheme that suitably combines a forecast/analysis-type discretisation with a robust time-integration scheme of the DLR system for the forecast step.
A hyper-reduction technique is furthermore proposed to efficiently evaluate nonlinear terms in the \moda{(augmented state) dynamics}. 

In~\cite{DynamicalLoTrigo2025}, the DLR methodology was used to perform a Dynamical Low-Rank approximation of \modb{the Kalman-Bucy Process (\kbp), which corresponds to the dynamics \eqref{eqn:kbp-dynamics-continuous-state} with an affine drift term $f$ and possibly an extra additive noise term $\Sigma^{\nicefrac{1}{2}} \md \Wcal_t$.}
The \modb{\dlr} method \modb{proposed in \cite{DynamicalLoTrigo2025}} approximates 
\begin{equation*}
	\Xcal_t^{\KBP} \approx \Xcal_t^{\DLR} = U_t^0 + \sum_{i=1}^R U_t^i Y_t^i \equiv U_t^0 + \bU_t \bY_t^{\top},
\end{equation*}
by modifying the dynamics in such a way that it is possible to characterise the evolution of the factors in the right-hand-side.
The above decomposition separates the process in its \moda{conditional mean} $U_t^0 \in \R^d$ and \moda{a low-rank zero conditional mean part. 
Therefore,}
the stochastic modes $\bY_t \in [L^2(\Omega)]^R$ verify $\E[Y_t^i \st \Acal_{\Zcal_t}] = 0$ in order to avoid redundancy with the mean.
In this setting, the physical modes $\bU_t \in \R^{d \times R}$ are orthogonal at all times, and consequently carry the information of the range of the covariance. 
The \modb{resulting} process is \modb{named} \dlr-\kbp. 

The \dlr-\enkf\ method \modb{also proposed in \cite{DynamicalLoTrigo2025} is a particle approximation of \dlr-\kbp. 
It} approximates the \moda{the law of $\Xcal_t^{\DLR}$ by the ensemble $\hXcal_t^{\DLR} = [\hXcal_t^{\DLR,(1)}, \ldots, \hXcal_t^{\DLR,({\rc{N}})}] \in \R^{d \times {\rc{N}}}$ where} 
\begin{equation*}
	\moda{\hXcal_t^{\DLR, (p)} 
	= 
	\widehat{U}_t^0 + \sum_{i=1}^R \hU_t^i \hY_t^{i,(p)}
	\equiv
\widehat{U}_t^0 + \hbU_t (\hbY_t^{(p)})^{\top}}.
\end{equation*}
\moda{Notice that the deterministic modes $\{\hU^i_t\}_{i=0}^R$ are common to all the particles, that is, all particles belong to the same affine low-dimensional space $\hU_t^0 + \mspan(\hU_t^1, \ldots, \hU_t^{R})$.
	We denote $\hbY_t^{(p)} =  \left[\widehat{Y}^{1,(p)}_t, \ldots, \widehat{Y}^{R,(p)}_t \right]$
	the $p$-th $R$-dimensional particle in the reduced space $\mspan(\hU_t^1, \ldots, \hU_t^{R})$ and $\hbY_t = [\hY_t^{(1)} ;  \ldots ; \hY_t^{({\rc{N}})}] \in \R^{{\rc{N}} \times R}$ the ensemble of reduced particles.
}
In analogy to the mean-field case, $\widehat{U}^0_t = \E_{\rc{N}}[\hXcal_t^{\DLR}]$ is the sample \modb{conditional} mean, $\hbY_t$ verify $\E_{\rc{N}}[\hbY_t] = 0$, and the sample \modb{conditional} covariance is given \moda{by $\widehat{P}_t = \hbU_t \hP_{\hbY_t} \hbU_t^{\top}$, where $\hP_{\hbY_t} = \frac{1}{{\rc{N}}-1} \hbY_t^{\top} \hbY_t$ is the reduced sample \modb{conditional} covariance}. 

Extending the \moda{\dlr-\kbp\ and}  \dlr-\enkf\ \moda{dynamics proposed in~\cite{DynamicalLoTrigo2025}} to nonlinear dynamics involves a few caveats. 
\moda{
	The paper~\cite{DynamicalLoTrigo2025} \modb{actually} proposes a DLR formulation for a general \modb{McKean-Vlasov} SDE of the form
\begin{equation} \label{eqn:gen-sde}
	\md \Xcal_t  = a(\Xcal_t, t, \modb{\mu_t}) \md t + b(\Xcal_t, t, \modb{\mu_t}) \md \Vcal_t + c(\Xcal_t, t, \modb{\mu_t}) \md \Zcal_t
\end{equation}
\modb{with $\mu_t = \Law(X_t)$,} involving a Brownian motion $\Vcal_t$ and a general observation process $\Zcal_t$. 
Starting from this, it derives the \dlr-\kbp\ and \dlr-\enkf\ processes for linear-affine dynamics.
This latter assumption allows one to perform important simplifications. 
In this work we focus on nonlinear dynamics, therefore we start again from the DLR formulation for~\cref{eqn:gen-sde} proposed in~\cite{DynamicalLoTrigo2025} and ``formally'' derive suitable \dlr-\kbp\ and \dlr-\enkf\ processes for nonlinear dynamics.
We focus first on the mean-field \venkf\ process~\eqref{eqn:kbp-dynamics-continuous-state} with $\kappa=1$, which we rewrite here for convenience in the form of~\cref{eqn:gen-sde}, i.e.
}
%
\begin{align*}
	\rd \Xcal_t &= (f(\Xcal_t, \theta_t) - P_t S \Xcal_t) \rd t 
	- P_t H^{\top} \Gmoh \rd \Vcal_t
	+ P_{t} H^{\top} \Gamma^{-1} \rd \Zcal_t 
		\\
		    &=  a \rd t + b \rd \Vcal_t + c \rd \Zcal_t,
\end{align*}
where $S = H^{\top} \Gamma^{-1} H$ and $P_t$ is the (conditional) covariance of $\Xcal_t$.
Anticipating the developments below, the following quantities will be of interest.
\begin{align*}
	\E[a \st \Acal_{\Zcal_t}] &= \E [ f(\Xcal_t, \theta_t) \st \Acal_{\Zcal_t}] - P_t S \E[\Xcal_t \st \Acal_{\Zcal_t}], &
	\moda{a_t^{\star}} &= \moda{f^{\star}_t}(\Xcal_t, \theta_t) - P_t S \Xcal_t^{\star}, 
		     \\
	\E[c \st \Acal_{\Zcal_t}] &= P_t H^{\top} \Gamma^{-1}, 
				  & \moda{c^{\star}_t}& = 0,
\end{align*}
Using the above relations in the mean separated DO equations derived in~\cite[Section 3]{DynamicalLoTrigo2025}, we obtain the following system
\begin{align*}
	\md U_t^0 &= \E [ f(\Xcal_t, \theta_t) \st \Acal_{\Zcal_t}] \rd t + P_t H^{\top} \Gamma^{-1} (\rd \Zcal_t - H U^0_t \rd t),  \\
	\md \bU_t &= P^{\perp}_{\bU_t} \E[ \moda{f^{\star}_t}(\Xcal_t, \theta_t) (\moda{\bY_t M_{\bY_t}^{-1}}) \st \Acal_{\Zcal_t}]  \md t ,
        \\
	\md \moda{\bY_t^{\top}} &= \bU_t^{\top} \moda{f^{\star}_t}(\Xcal_t, \theta_t) \md t 
	- \bU_t^{\top} P_{t} H^{\top} \Gamma^{-1}( H \Xcal_t^{\star} + \Goh \md \Vcal_t ), 
\end{align*}
where $M_{\bY_t} = \E[Y_t^i Y_t^j \st \Acal_{\Zcal_t}]$ is the (continuous) Gram matrix.  
Based on \moda{these} equations, \modb{and introducing a particle approximation, we obtain the} \dlr-\venkf~system proposed in~\cite{DynamicalLoTrigo2025} \moda{\modb{for} nonlinear \modb{dynamics}}
\begin{equation} \label{eqn:dlr-venkf-nonlin}
	\begin{aligned}
		\rd \widehat{U}_t^0 &=  \E_{\rc{N}}[f(\hXcal_t, \htheta_t)] \rd t 
	+ \widehat{P}_t H^{\top} \Gamma^{-1} (\rd \Zcal_t - H \moda{\hU^0_t} \rd t - \Goh \md \E_{\rc{N}} [\widehat{\Vcal}_t]),
		\\
		\rd \moda{\hbU_t} &= P^{\perp}_{\moda{\hbU_t}} \E_{\rc{N}}[ \moda{f_t^{\star}}(\hXcal_t, \htheta_t) \hbY_t] \moda{\hM_{\hbY_t}^{-1}} \md t ,
		\\
		\rd \hbY_t^{\top} &= 
		\moda{\hbU_t^{\top}} \moda{f_t^{\star}}(\hXcal_t, \htheta_t) \md t
		- \moda{\hbU_t^{\top}} \widehat{P}_{t} H^{\top} \Gamma^{-1} ( H \hXcal_t^{\star} + \Goh \md \hVcal_t^{\star}).
	\end{aligned}
\end{equation}
where $\moda{\hM_{\hbY_t} = \frac{1}{{\rc{N}}} \widehat{\bY}_t^{\top} \widehat{\bY}_t}$  \moda{denotes the empirical} Gram matrix. 
Observe that the sample mean of the Brownian motions \moda{is} subtracted from the equation for the stochastic modes to ensure \moda{that} they remain \moda{(empirical)} zero-mean, and \modb{is} added back in the sample mean equation for consistency.

Performing the same steps for the \moda{\denkf~process} yields the following \dlr-\denkf~system
\begin{equation} \label{eqn:dlr-denkf-nonlin}
	\begin{aligned}
		\rd \widehat{U}_t^0 &=  \E_{\rc{N}}[f(\hXcal_t, \htheta_t)] \rd t 
		+ \widehat{P}_t H^{\top} \Gamma^{-1} (\rd \Zcal_t - H U^0_t \rd t), 
		\\
		\rd \moda{\hbU_t} &= P^{\perp}_{\moda{\hbU_t}} \E_{\rc{N}}[ \moda{f_t^{\star}}(\hXcal_t, \htheta_t) \hbY_t] \moda{\hM_{\hbY_t}^{-1}} \md t ,
		\\
	\rd \hbY_t^{\top} &= 
	\moda{\hbU_t^{\top}} \moda{f_t^{\star}}(\hXcal_t, \htheta_t) \md t
	- \frac{1}{2} \moda{\hbU_t^{\top}} \widehat{P}_t H^{\top} \Gamma^{-1} H \hXcal_t^{\star} \md t.
\end{aligned}
\end{equation}


\subsection{Forecast-Analysis time-discretisation} \label{sec:forecast-analysis-time-discretisation}

The inverse Gram matrix $\moda{\hM_{\hbY_t}^{-1}}$ appearing in the above equations constitutes an additional source of stiffness when considering suitable time-integration of the above DLR systems. 
\modb{This consideration has} led to the development of various robust time-integrators in the DLR literature, of which we mention the Basis Update and Galerkin (BUG) integrator~\cite{ARankAdaptiveCeruti2022,RobustHiRiffaud2025}. 
\moda{The} time-integrator proposed in~\cite{DynamicalLoTrigo2025} has a BUG-like structure combined with Euler-Maruyama steps. 
In this work, we consider a different algorithm based on the forecast/analysis structure. 
The algorithm is again described over the time-interval $[t_0, t_1]$.

\textbf{Forecast step.}
For both \dlr-\venkf~and \dlr-\denkf, the dynamics \moda{is given} by the first terms \moda{of~\cref{eqn:dlr-venkf-nonlin}, dropping all the terms involving the Kalman gain}
\begin{align*}
	\rd \moda{\hU_t^0} &=  \E_{\rc{N}}[f(\hXcal_t, \htheta_t)] \rd t, 
		  &
	\rd \moda{\hbU_t} &= P^{\perp}_{\moda{\hbU_t}} \E_{\rc{N}}[ \moda{f_t^{\star}}(\hXcal_t, \htheta_t) \hbY_t] \moda{\hM_{\hbY_t}^{-1}} \md t ,
		  &
	\rd \hbY_t^{\top} &= \moda{\hbU_t^{\top}} f^{\star}(\hXcal_t, \htheta_t) \md t.
\end{align*}
\moda{Rather} than \moda{updating the mean separately}, \moda{we treat the} solution \moda{$\hU_t^0 + \hbU_t \hbY_t^{\top}$} as a rank-$R+1$ \moda{function}, which allows to straightforwardly use the BUG algorithm.
The first step consists in computing augmented subspaces.
To this end, denote $\hbZ = [\bone, \hbY_{t_0}]$ and $\hbV = [U^0_{t_0}, \bU_{t_0}]$. 
Following~\cite{ARankAdaptiveCeruti2022} and adapting it to the present two-field decomposition, on $[t_0, t_1]$ \moda{we solve
\begin{align*}
	\bK(t_1) &= \bK(t_0) + \delt \E_{\rc{N}}[f(\hbV \hbZ^{\top}, \htheta_{t_0}) \hbZ], &
	\bL(t_1)^{\top} &= \bL(t_0)^{\top} + \delt \hbV^{\top} f(\hbV \hbZ^{\top}, \htheta_{t_0}),
\end{align*}
and} construct the augmented bases
\begin{align*}
	\lbU &= \ortho([\hbV, \bK(t_1)]),
	     &
	\lbY &= \ortho([\hbZ, \bL(t_1)]).
\end{align*}
The solution is then updated \modb{by Galerkin projection} as follows
\begin{equation*}
	\moda{\bS(t_1) =  \bS(t_0) + \delt \lbU^{\top} \E_{\rc{N}} [ f(\hXcal_{t_0} , \htheta_t) \lbY], }
\end{equation*}
The augmented solution \moda{is then} $\hXcal_{t_1}^{f,\mathrm{aug}} = \lbU \bS(t_1) \lbY^{\top}$, \moda{which can also be written as $\hXcal_{t_1}^{f,\mathrm{aug}}  = (\hU^0_{t_1})^f + \lbU \tilde{S}(t_1) \lbY^{\top}$ with $\E_{\rc{N}}[\lbU \tilde{S}(t_1) \lbY^{\top}] = 0$ since} the constant vector $\bone$ belongs to $\lbY$. 
Finally, truncate the zero-mean part to rank-$R$ using a generalised SVD, yielding 
\begin{equation}
\label{eq:truncation}
\Tcal_R \left( \moda{\lbU \tilde{S}(t_1) \lbY^{\top}} \right) = \moda{\hbU_{t_1}^f} (\hbY_{t_1})^{\top}.
\end{equation}
The forecasted low-rank solution is given as $\hXcal_{t_1}^f = (U_{t_1}^0)^f + \moda{\hbU_{t_1}} (\hbY_{t_1}^{f})^{\top}$. 

\textbf{Analysis step.} The analysis dynamics \moda{is}
\begin{align*}
	\rd \moda{\hU_t^0} &= \widehat{P}_t H^{\top} \Gamma^{-1} (\rd \Zcal_t - H \moda{\hU^0_t} \rd t - \kappa \Goh \md \E_{\rc{N}}[ \moda{\hVcal_t} ]), \\
	\rd \hbY_t^{\top} &= \bU_t^{\top} \widehat{P}_t H^{\top} \Gamma^{-1} \left( - \frac{1 + \kappa}{2}  H \hXcal_t^{\star} \md t \moda{-} \kappa \Goh \md \moda{\hVcal_t^{\star}} \right),
\end{align*}
\moda{where we recall that $\hP_t$ is the sample covariance. 
Note that we recover} the \dlr-\venkf~scheme with $\kappa = 1$, and \dlr-\denkf~with $\kappa = 0$. 
The dynamics of the mean \moda{is} structurally similar to the \enkf~analysis dynamics \moda{as described in~\cref{sec:time-discretisation-schemes}}, and consequently following the same steps as for the \venkf~ resp. \denkf~\moda{to derive~\cref{eqn:analysis-interm}} yields the update 
\begin{equation}
	(\moda{\hU_{t_1}^0)^{a}} = (\moda{\hU_{t_1}^0)^f}  + K_{t_1} \left( \D \Zcal_{t_1} - H (\moda{\hU_{t_1}^0)^f} \delt  - \kappa \Goh \E_{\rc{N}} [\D \moda{\hVcal_{t_1}}]  \right), \label{eqn:mean-analysis}
\end{equation}
The Kalman gain $K_{t_1}$ has a low-rank structure that will be detailed below.
Conversely, the dynamics for $\hbY_t$ can be interpreted as an update in the subspace $\mathrm{span}(\moda{\hbU_{t_1}})$.
Indeed, recalling that $\hP_{t_1}^f = \moda{\hbU_{t_1}} \moda{\hP_{\hbY_{t_1}^f}} \moda{\hbU_{t_1}^{\top}}$, \moda{where $\hP_{\hbY_{t_1}^f} = \frac{1}{{\rc{N}}-1}(\hbY_{t_1}^f)^{\top} \hbY_{t_1}^f$}, consider the semi-implicit discretisation
\begin{align*}
	(\hbY_{t_1}^{a})^{\top} &= (\hbY_{t_0}^{a})^{\top} + \moda{\hbU_{t_1}^{\top}} \moda{\widehat{P}^f_{t_1}} H^{\top} \Gamma^{-1} \left( - H (\hXcal_{t_1}^{\star})^a \delt + \frac{1 - \kappa}{2} H (\hXcal_{t_1}^{\star})^f \delt  - \kappa \Goh \D \moda{\hVcal_{t_1}^{\star}} \right) \\
				&= (\hbY_{t_0}^{a})^{\top} + \moda{\hP_{\hbY_{t_1}^f}} H_{\bU}^{\top} \Gamma^{-1} \left( - H_{\bU} (\moda{\hbY_{t_1}^a})^{\top} \delt + \frac{1 - \kappa}{2} H_{\bU} (\hbY_{t_1}^f)^{\top} - \kappa \Goh \D \moda{\hVcal_{t_1}^{\star}} \right) 
\end{align*}
where $H_{\bU} = H \bU_{t_1}$. 
Therefore, \moda{following the \modb{same steps used to} derive~\cref{eqn:analysis-interm}} it holds
\begin{equation}
	(\hbY_{t_1}^{a})^{\top} = (\hbY_{t_1}^{f})^{\top} + K_{t_1}^{\bU} \left( -  \frac{1 + \kappa}{2} H_{\bU} (\hbY_{t_1}^f)^{\top} \delt \moda{-} \kappa \Goh \D \moda{\hVcal_{t_1}^{\star}} \right), \label{eqn:zmean-analysis-subspace}
\end{equation}
where
\begin{equation*}
	K_{t_1}^{\bU} = \moda{\hP_{\hbY_{t_1}^f}} H_{\bU}^{\top} (\Gamma + \delt H_{\bU} \moda{\hP_{\hbY_{t_1}^f}} H_{\bU}^{\top})^{-1}.
\end{equation*}
Note that therefore, 
\begin{align}
	\moda{\hbU_{t_1}} (\hbY_{t_1})^{a} &= \moda{\hbU_{t_1}} (\hbY_{t_1}^{f})^{\top} + \moda{\hbU_{t_1}} K_{t_1}^{\bU} \left( -  \frac{1 + \kappa}{2} H_{\bU} (\hbY_{t_1}^f)^{\top} \delt \moda{-} \kappa \Goh \D \moda{\hVcal_{t_1}^{\star}} \right) \nonumber \\
					   &= \moda{\hbU_{t_1}} (\hbY_{t_1}^{f})^{\top} +  K_{t_1} \left( -  \frac{1 + \kappa}{2} H \moda{\hbU_{t_1}} (\hbY_{t_1}^f)^{\top} \delt \moda{-} \kappa \Goh \D \moda{\hVcal_{t_1}^{\star}} \right) \label{eqn:zmean-analysis-ambient}
\end{align}
Combining~\cref{eqn:mean-analysis,eqn:zmean-analysis-ambient}, the \moda{analysis} state $\hXcal_{t_1}^{a} = (\moda{\hU_{t_1}^0)^{a}}  + \moda{\hbU_{t_1}} (\hbY_{t_1}^{a})^{\top}$ verifies the same state update relation as~\cref{eqn:analysis}, which we rewrite here for clarity as
\begin{equation} \label{eqn:analysis-low-rank}
	\hXcal_{t_1}^a = \hXcal_{t_1}^f + \moda{\hbU_{t_1}} K_{t_1}^{\bU}\left( \D \Zcal_{t_1} -  \frac{1 + \kappa}{2} H \hXcal_{t_1}^{f} \delt - \kappa \Goh \D \moda{\hVcal_{t_1}}  - \frac{1- \kappa}{2} \delt H \moda{\hU_{t_1}^0} \right).
\end{equation}

\begin{remark}
	In the forecast step, it is also possible to forego the explicit inclusion of the mean in the basis, keeping the remainder of the algorithm unchanged. 
	While the mean is not \moda{explicitly tracked}, it results in a more straightforward implementation, closer to the original BUG scheme~\cite{ARankAdaptiveCeruti2022}. 
\end{remark}

\begin{remark}
\moda{In the truncation step \eqref{eq:truncation}, when the state has several variables, we truncate each variable separately rather than applying a single ``global'' truncation, since they may have different orders of magnitude. For example, in Section \ref{sec:numexp-3}, the variables $A$ and $u$ appear in the augmented state as $(\hXcal_{t_1}^{f,\mathrm{aug}})^{\star}=\begin{bmatrix}(\hXcal_{t_1}^{f,\mathrm{aug}})^{\star}_A\\ (\hXcal_{t_1}^{f,\mathrm{aug}})^{\star}_u\end{bmatrix}$, and the truncated state is obtained from the separate truncations $\Tcal_R \left((\hXcal_{t_1}^{f,\mathrm{aug}})^{\star}_A\right)$ and $\Tcal_R \left((\hXcal_{t_1}^{f,\mathrm{aug}})^{\star}_u\right)$. More precisely, assuming $2R\le\min\{d,{\rc{N}}\}$, we set
\begin{equation*}
\hbU_{t_1}^f =
\begin{bmatrix}
Q_A V_A & \\
& Q_u V_u
\end{bmatrix}
\quad \textrm{and} \quad
\hbY_{t_1}^{f} = \lbY
\begin{bmatrix}
W_A \Sigma_A & W_u \Sigma_u
\end{bmatrix},
\end{equation*}
where $\lbU = \begin{bmatrix} \lbU_A  \\ \lbU_u \end{bmatrix}$, $\lbU_A = Q_A R_A$ (resp. $\lbU_u = Q_u R_u$) is the QR factorisation of $\lbU_A$ (resp. $\lbU_u$), and $R_A \tilde{S}(t_1) = V_A \Sigma_A W_A^\top$ (resp. $R_u \tilde{S}(t_1) = V_u \Sigma_u W_u^\top$) is the rank-R truncated SVD of $R_A \tilde{S}(t_1)$ (resp. $R_u \tilde{S}(t_1)$)}.
\end{remark}

\rb{  
\subsubsection{Initial rank \& rank-adaptivity mechanisms}
\label{sec:initial-rank-rank-adaptivity}

The rank plays a central role in the proposed \dlr-\enkf~method. In a filtering context, the assimilation procedure is expected to progressively concentrate the posterior distribution, and therefore the ensemble of particles may eventually be well represented in a relatively low-dimensional subspace. This is precisely the regime in which a dynamically evolving low-rank representation is advantageous. However, if the rank is chosen too small, the range of the covariance is artificially restricted, which may lead to covariance underestimation, overconfidence, and biased parameter estimates. Conversely, choosing a rank that is unnecessarily large increases the computational cost without significantly improving the accuracy.

In a fixed-rank setting, choosing the rank requires some practical calibration. Possible strategies include running short full-order simulations to estimate the rank, or propagating several \dlr~approximations with increasing ranks and retaining the smallest rank beyond which the estimates no longer vary significantly. 

Rank-adaptive methods provide a more flexible alternative by allowing the rank to vary during the simulation. In this work, we employ the rank-adaptive BUG strategy of~\cite{ARankAdaptiveCeruti2022}. This algorithm requires two user-defined parameters: an initial rank $R_{\mathrm{init}}$ and a truncation threshold $\vartheta>0$. The rank is then updated according to
\begin{equation} \label{eqn:new-rank}
R_{\mathrm{new}} = \argmin_{R \geq 1} \left\{ R \;\bigg|\, \sum_{j = R+1}^{2 R_{\mathrm{old}}} \sigma_j^2(\tilde{S}(t_1)) \leq \vartheta^2 \right\},
\end{equation}
where $R_{\mathrm{new}}$ and $R_{\mathrm{old}}$ denote the new and current ranks, respectively. The initial rank should be chosen sufficiently large to capture the dominant directions of the ensemble during the first iterations; if it is too small, the recovery of missing directions may be delayed, whereas an excessively large initial rank mainly increases the computational cost. Additionally, $\vartheta$ controls the amount of singular-value energy discarded during the truncation step. Choosing the truncation threshold too large may discard important dynamics, while choosing it too small may result in an unnecessary increase of the rank.

This adaptive strategy works well in practice. It has nevertheless two well-known limitations. First, since the augmented basis generated by the BUG integrator has rank at most $2R_{\mathrm{old}}$, the rank can increase by at most a factor two at each iteration. Consequently, if the intrinsic rank of the solution grows much faster than this, the algorithm may temporarily fail to capture all relevant directions. 
This is unlikely to happen during the simulation if $\delt$ is small enough, but it could happen at the first few iterations if the initial rank was chosen too small.
Second, the discarded singular values only quantify the truncation error within the augmented BUG subspace, and do not provide a complete error estimate for components orthogonal to this subspace. 

These limitations also suggest practical diagnostics: a persistent increase of the rank, and, in particular, repeated occurrences of $R_{\mathrm{new}}=2R_{\mathrm{old}}$, indicate that the current rank is too small or that the initial rank was chosen too conservatively. More sophisticated strategies, such as predictor-corrector approaches~\cite{APredictorCorrectorHauck2023}, periodically estimate the neglected dynamics outside the current subspace and enrich the basis accordingly, thereby providing additional information on whether the current rank is sufficient, at the expense of an increased computational cost.
In our numerical experiments we did not see the need for this more sophisticated and costly rank-adaptive technique. 
}


\rb{
\subsubsection{DLR-ENKF signal-tracking abilities}
Our} proposed method has the potential to perform efficient data assimilation by dynamically evolving a low-rank analog of the \enkf. 
The accuracy of our method is contingent on the true filtering measure effectively exhibiting a \modb{nearly} low-rank profile at all times, and in this context we discuss two linked but distinct potential effects that may degrade the signal-tracking ability of the \dlr-\enkf. 
The first is the neglect of observation information, and the second is covariance under-approximation.
Both of these effects stem from the low-rank approximation of the (range of the) covariance, which is the mechanism making our method efficient in the first place; after detailing the issues, we discuss possible mitigations.

The first effect expresses the fact that our method only assimilates the \modb{component of the} (rescaled and state-space-embedded) observations $\widetilde{\D \Zcal_{t_1}} = H^{\top} \Gamma^{-1} \D \Zcal_{t_1}$ lying in the range of $\moda{\hbU_{t_1}}$. 
This is readily seen from~\cref{eqn:analysis-low-rank}; focussing on the Kalman gain part with the data increments,
\begin{align*}
	K_{t_1}^{\bU} \D \Zcal_{t_1} &= \moda{\hP_{\hbY_{t_1}^f}} H_{\bU}^{\top} (\Gamma + \delt H_{\bU} \moda{\hP_{\hbY_{t_1}^f}} H_{\bU}^{\top})^{-1} \D \Zcal_{t_1}  \\
				     &= \underbrace{\moda{\hP_{\hbY_{t_1}^f}} \left[ \bI_{R} - H_{\bU}^{\top} \Gamma^{-1} H_{\bU} \left( \moda{\hP_{\hbY_{t_1}^f}^{-1}} / \delt + H_{\bU}^{\top} \Gamma^{-1} H_{\bU}  \right)^{-1} \right]}_{\modb{=T}} \moda{\hbU_{t_1}^{\top}} (H^{\top} \Gamma^{-1} \D \Zcal_{t_1}) \\
				     &= T \moda{\hbU_{t_1}^{\top}} \Pcal_{\moda{\hbU_{t_1}}}\widetilde{\D \Zcal_{t_1}},
\end{align*}
having used the Woodbury identity in the second equality, and where $\moda{\Pcal_{\hbU_{t_1}} = \hbU_{t_1} \hbU_{t_1}^{\top}}$. 
In particular, when $H = \bI_{d} = \Gamma$, it holds
\begin{equation*}
	K_{t_1}^{\bU} \D \Zcal_{t_1} = T \moda{\hbU_{t_1}^{\top}} \Pcal_{\moda{\hbU_{t_1}}} \D \Zcal_{t_1} = K_{t_1}^{\bU} \Pcal_{\moda{\hbU_{t_1}}} \D \Zcal_{t_1},
\end{equation*}
and hence the \dlr-\enkf\ only assimilates the \modb{component of the} observations in \moda{the} subspace of $\moda{\hbU_{t_1}}$. 
\moda{Note that the \enkf\ operates under the same principle, assimilating the observations only in the subspace spanned by the particles.
The distinction is that, in the common regime ${\rc{N}} \ll d$, this subspace is (at most) ${\rc{N}}$-dimensional rather than being limited at $R$.
However, when the ${\rc{N}}$ particles concentrate on an $R$-dimensional subspace that the \dlr-\enkf\ identifies, the assimilations performed by the \enkf\ and the \dlr-\enkf\ are expected to be comparable.}

The second effect is related to the potential under-approximation of the sample covariance due to rank underestimation. 
Consider the setting where the full order model \enkf\ has a covariance $\hP_t^{\ENKF}$ with a suitable low-rank approximation of rank $R'$. 
If the system is simulated using the \dlr-\enkf\ method using a rank $R < R'$, and for simplicity assuming that $\hP_t^{\DLR} \approx \Tcal_R(\hP_t^{\ENKF})$, it holds
\begin{equation*}
	\norm{\hP^{\DLR}_t}_F^2 \approx \norm{\Tcal_R(\hP^{\ENKF}_t)}_F^2 = \sum_{i=1}^R \lambda_i^2(\hP^{\ENKF}_t) < \sum_{i=1}^{R'} \lambda_i^2(\hP^{\ENKF}_t) \approx \norm{\hP^{\ENKF}_t}^2_F.
\end{equation*}
In effect, the \dlr-\enkf\ may underestimate the true covariance of the system, causing it to express higher confidence in a given state than warranted; furthermore, increasing the rank may at first cause the uncertainty to \textit{increase} before stabilising when the rank is sufficiently large.
In our numerical \modb{results}, the issue of uncertainty under-approximation seems to be more prevalent than the issue of neglecting the signal. 
This can however be \rb{avoided or mitigated} by choosing a sufficiently large rank \textit{a priori}, \rb{or by resorting to rank-adaptive strategies as discussed in~\cref{sec:initial-rank-rank-adaptivity}.}

\subsection{Hyper-reduction}

When $\moda{f_t^{\star}}(\hXcal_t, \htheta_t)$ is nonlinear with respect to $\hXcal_t$, the evaluation of $\E_{\rc{N}}[ \moda{f_t^{\star}}(\hXcal_t, \htheta_t) \hbY_t ]$ and $\bU_t^{\top} \moda{f_t^{\star}}(\hXcal_t, \htheta_t)$ scales with $d \times {\rc{N}}$ (instead of $d \times R$ and $R \times {\rc{N}}$). 
To reduce this cost, we employ a CUR approximation based on the discrete empirical
interpolation method (DEIM) \cite{AnEmpiricalIBarrau2004,NonlinearModelChatur2010}. 
The idea is to approximate the nonlinear term $\moda{f_t^{\star}}(\hXcal_t, \htheta_t) \in \Mat{d \times {\rc{N}}}$ by evaluating it only at a few row $\{{\sigma^1_t},\ldots,{\sigma^n_t}\}$ and column $\{{\varsigma^1_t}, \ldots, {\varsigma^m_t}\}$, with $n \ll d$ and $m \ll {\rc{N}}$. 
Let $\moda{{\mathbf P}_\sigma} \in \Mat{n \times d}$ (resp. $\moda{{\mathbf P}_\varsigma} \in \Mat{m \times \rc{N}}$) denote the matrix extracting the rows indexed by $\moda{\sigma = \{{\sigma^1_t},\ldots,{\sigma^n_t}\}}$ (resp. columns indexed by $\moda{\varsigma = \{{\varsigma^1_t}, \ldots, {\varsigma^m_t}\}}$). The nonlinear term is approximated by the CUR factorisation
\begin{equation}
\label{eq:hyper-reduction}
\moda{f_t^{\star}}(\hXcal_t, \htheta_t) \approx \left( \moda{f_t^{\star}}(\hXcal_t, \htheta_t) \moda{{\mathbf P}^\top_\varsigma} \right) \left( \moda{{\mathbf P}_\sigma} \moda{f_t^{\star}}(\hXcal_t, \htheta_t) \moda{{\mathbf P}^\top_\varsigma} \right)^{+} \left(\moda{{\mathbf P}_\sigma} \moda{f_t^{\star}}(\hXcal_t, \htheta_t) \right)^\top,
\end{equation}
where $(\cdot)^+$ stands for the Moore--Penrose pseudoinverse. The indices $\{{\sigma^1_t},\ldots,{\sigma^n_t}\}$ and $\{{\varsigma^1_t}, \ldots, {\varsigma^m_t}\}$ are selected by the DEIM, with $n=R$ and $m=2R$. 
Specifically, $R$ columns $\moda{\tilde\varsigma} = \{\varsigma^1_t,\ldots,\varsigma^R_t\}$ are first chosen by Algorithm \ref{al:deim} with input $\ortho(\hbY_t)$. \moda{This corresponds to choosing $R$ particles.} Then, $R$ rows $\moda{\sigma = \{{\sigma^1_t},\ldots,{\sigma^R_t}\}}$ are selected using Algorithm \ref{al:deim} with input $\ortho(\moda{f_t^{\star}}(\hXcal_t, \htheta_t) \moda{{\mathbf P}^\top_{\tilde\varsigma}})$, where $\moda{{\mathbf P}^\top_{\tilde\varsigma}} \in \Mat{{\rc{N}} \times R}$ extracts the \moda{columns} indexed by $\moda{\tilde\varsigma}$. Finally, $R$ additional columns $\{\varsigma^{R+1}_t,\ldots,\varsigma^{2R}_t\}$ are chosen using Algorithm \ref{al:deim} with input $\ortho(\moda{{\mathbf P}_\sigma} \moda{f_t^{\star}}(\hXcal_t, \htheta_t))$. \moda{Thus, the final ${\mathbf P}^\top_\varsigma$ in \eqref{eq:hyper-reduction} extracts the $2R$ columns indexed by $\varsigma = \{\varsigma^{1}_t,\ldots,\varsigma^{2R}_t\}$.}
Compared to the approach described in \cite{ObliqueProjBabaee2023}, this CUR approximation presents two main differences: (1) the sampling indices $\{{\sigma^1_t},\ldots,{\sigma^R_t}\}$ are determined from the left basis of the nonlinear term $\moda{f^{\star}_t}(\hXcal_t, \htheta_t) \moda{{\mathbf P}^\top_\varsigma}$ (which is more appropriate than using the left basis of the solution at the previous time-step), and (2) an oversampling of $R$ additional columns is used to improve accuracy.

\begin{algorithm}[H]
	\caption{DEIM point selection algorithm \cite{NonlinearModelChatur2010}}
	\label{al:deim}
	\begin{algorithmic}[1]
		\Require{An orthonormal matrix ${\mathbf Q} \in \R^{n \times m}$}
		\Ensure {A vector ${\mathbf p} \in \R^m$ containing the indices $\moda{\{\sigma^1,\ldots,\sigma^m\}}$}
		\State ${\mathbf q} = {\mathbf Q}(:,1)$
		\State $[\sim,\moda{\sigma^1}] = \max(|{\mathbf q}|)$
		\State ${\mathbf p} = [\moda{\sigma^1}]$
		\For {$j=2,\ldots,m$}
		\State ${\mathbf q} = {\mathbf Q}(:,j)$
		\State ${\mathbf r} = {\mathbf q} - {\mathbf Q}(:,1:j-1) \big({\mathbf Q}({\mathbf p},1:j-1)\big)^{-1} {\mathbf q}({\mathbf p})$
		\State $[\sim,\moda{\sigma^j}] = \max(|{\mathbf r}|)$
		\State ${\mathbf p} = [{\mathbf p};\moda{\sigma^j}]$
		\EndFor
	\end{algorithmic}
\end{algorithm}

\section{Numerical experiments} \label{sec:numexp}

This section assesses the qualities of the proposed \dlr-\enkf~method for parameter identification by benchmarking \moda{it} on two examples, a Fisher-KPP problem and a reduced one-dimensional blood flow model of the 55 arteries in the human body. 
These two models already showcase several relevant properties of the method, which are discussed in detail in the sections. 
To summarise our findings, for a moderate rank $R$ the \dlr-\enkf~consistently performs the same level of data assimilation and parameter identification as its \enkf~counterpart, and consistently does so at a reduced computational cost. 
The DEIM/CUR hyper-reduction proves effective in reducing the cost for hyper-quadratic nonlinearities (but offers no substantial gains for quadratic nonlinearities). 
One key finding of our numerical experiments is that the accuracy of \dlr-\enkf~is highly dependent on the approximation rank $R$, \moda{which} must be chosen \moda{large enough} to capture the rank of the state ensemble \textit{a priori}. 
If the rank is too small, the \dlr-\enkf~may under-estimate the covariance, eventually leading to inaccurate assimilation and parameter identification.  
\moda{This naturally calls for rank-adaptive strategies, which will be explored in this work and investigated further in future works.}



\subsection{Fisher-Kolmogorov-Petrovski-Piskunov system}

The following experimental setting is inspired by the one in~\cite{EfficientStatePagani2017}. 
The dynamics \modb{is} given by (the discretisation of) a Fisher-Kolmogorov-Petrovski-Piskunov (Fisher-KPP) reaction-diffusion equation, which \modb{models} the evolution of biological population invasions, genes, or chemical reactions over space and time. 
On the domain
\begin{equation*}
	D = \left\{ (x,y) \in \R^2 \st 1 \leq \sqrt{x_1^2 + x_2^2} \leq 1.5,\, x_1 \geq 0, x_2 \geq 0 \right\},
\end{equation*}
denote $u(t, x, \theta)$ the solution on $(0,T)$ with $T = 0.154$ to 
\begin{align*}
&	\partial_t u(t, x, \theta) - \mdiv( \nu(x, \theta) \nabla u) - 75 u (1 - u) = 0, && x \in D, t \in (0,T) \\
&	\nu(x, \theta) \nabla u(t, x, \theta) \cdot \bn = 0, && x \in D, t \in (0,T) \\
&	u(0, x, \theta) = \exp( -(x_1 - 1.5)^2 - 50 x_2^2 ) && x \in D.
\end{align*}
In this problem setting, the concentration $u$ diffuses with a preferential direction along the \moda{tangential direction} of the annulus arc with a wave-front structure; by $T = 0.154$, the wave-front has traversed almost the whole domain. 
The physical discretisation is performed via the Finite Element Method, using $\mP_1^C(\Tcal_h)$ the space of $\mP_1$-continuous finite elements on the quasi-uniform mesh $\Tcal_h \subset D$. 
Furthermore, denote $d = \mathrm{dim}(\mP_1^C(\Tcal_h)) = 540$.
The time-step discretisation is set to $\delt = 4.4 \cdot 10^{-5}$, and we use the explicit Euler method \moda{for the forecast step in the \enkf}, \modb{and} \moda{the explicit BUG} \modb{method for the forecast in} \moda{the \dlr-\enkf\ detailed in~\cref{sec:forecast-analysis-time-discretisation}.}

\begin{figure}[htbp]
\centering
\subfloat[$\xi_1$]{\includegraphics[width=0.32\linewidth]{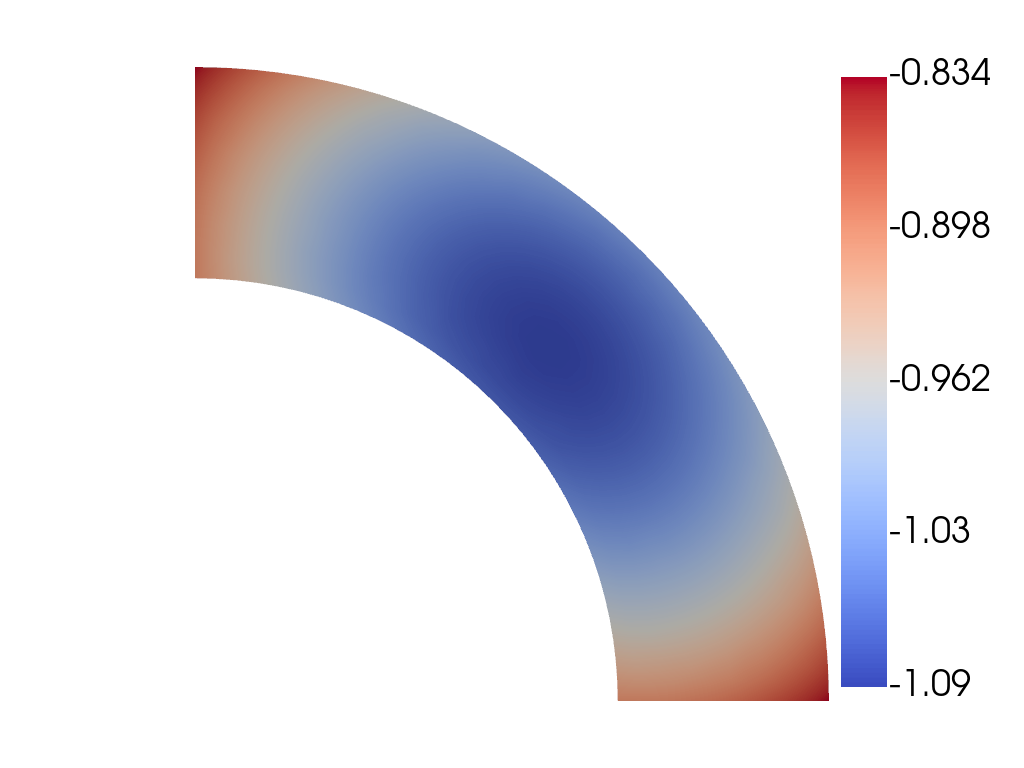}} 
\subfloat[$\xi_2$]{\includegraphics[width=0.32\linewidth]{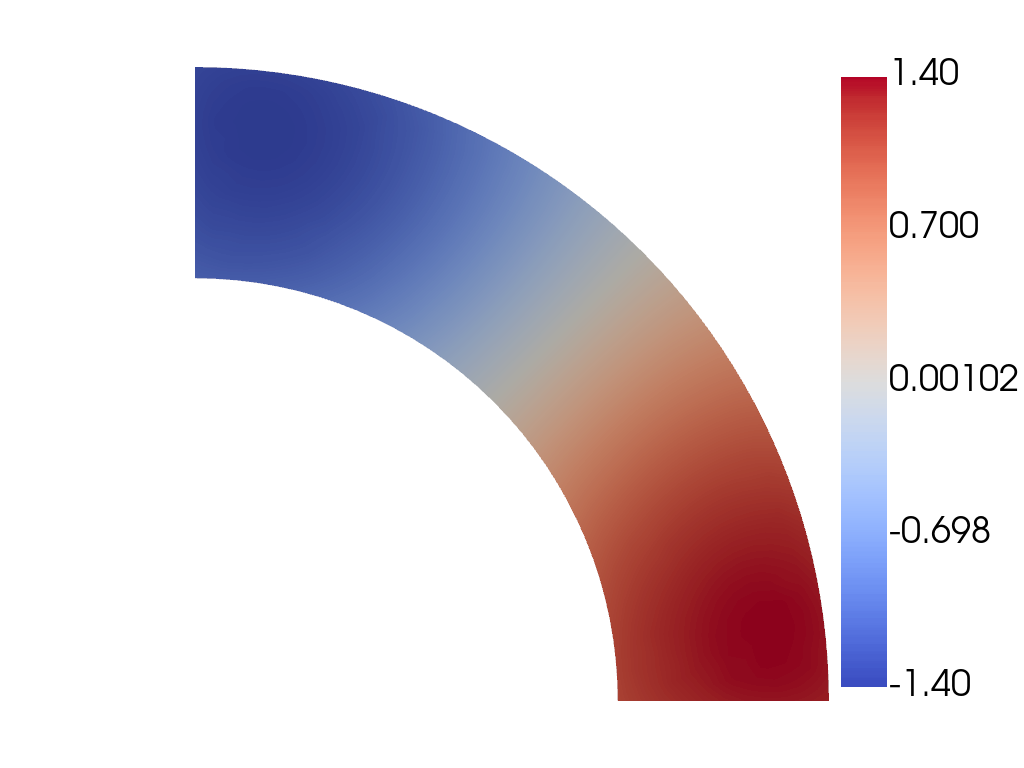}} 
\subfloat[$\xi_3$]{\includegraphics[width=0.32\linewidth]{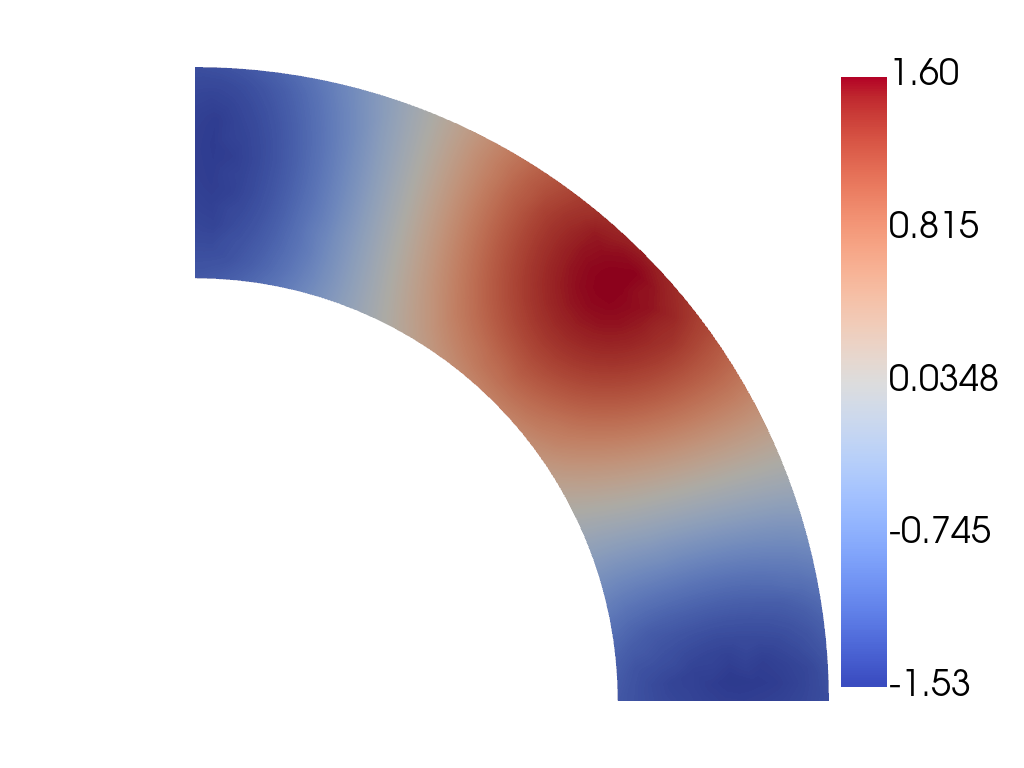}} \\
\subfloat[$\xi_4$]{\includegraphics[width=0.32\linewidth]{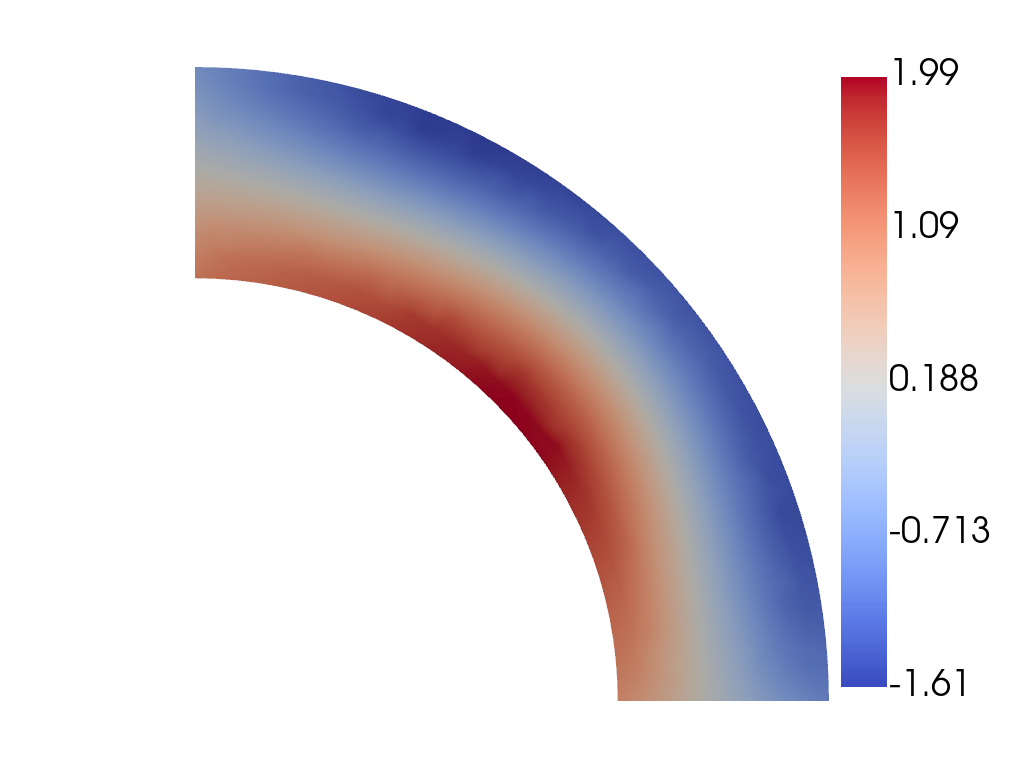}} 
\subfloat[$\xi_5$]{\includegraphics[width=0.32\linewidth]{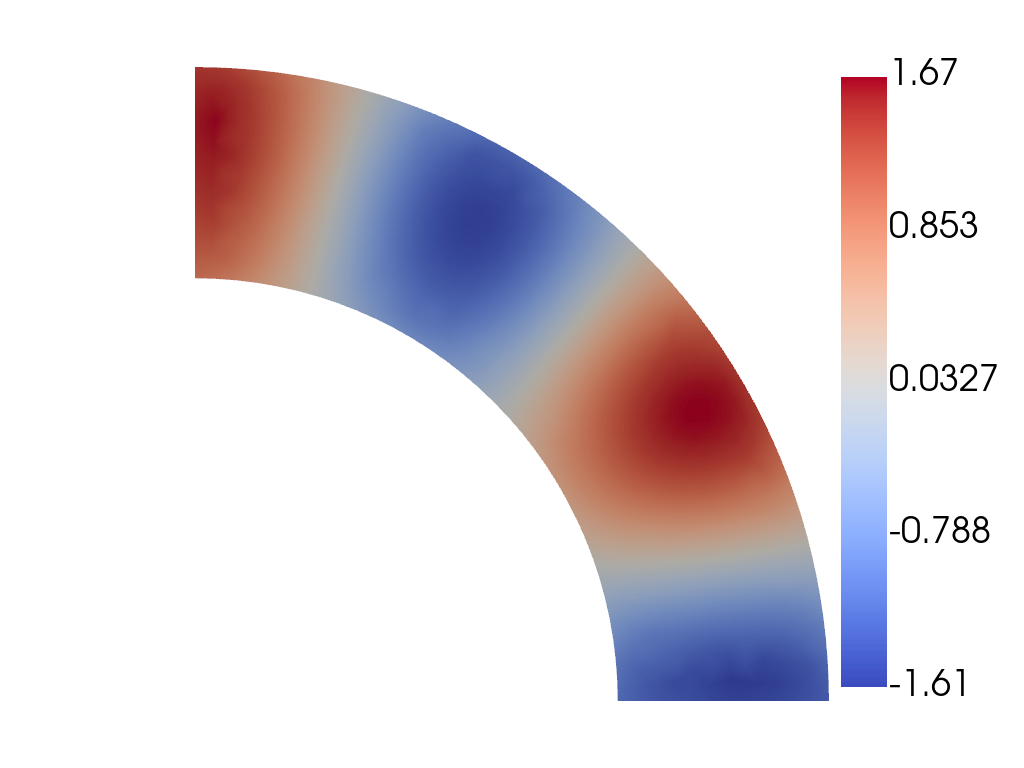}} 
\subfloat[$\xi_6$]{\includegraphics[width=0.32\linewidth]{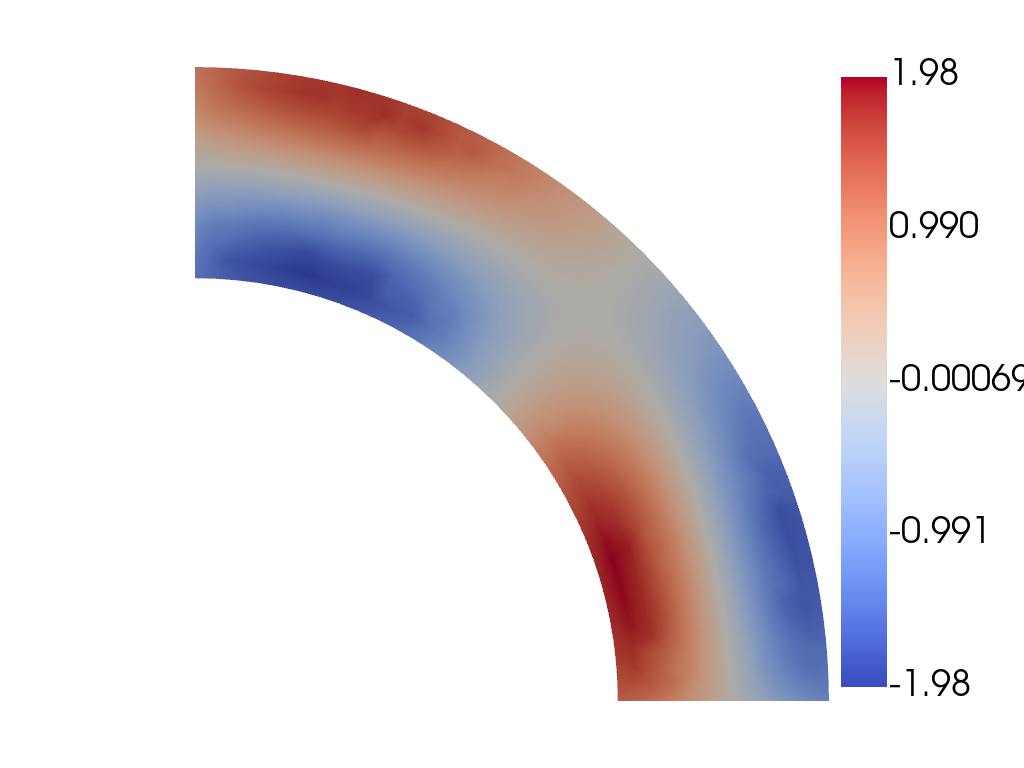}} 
\caption{\ra{Fisher-KPP: modes of the diffusion coefficient $\nu(x)$}}
\label{fig:eigenfunctions}
\end{figure}

\moda{For the diffusion coefficient, we follow~\cite{EfficientStatePagani2017} and define the parametric diffusion coefficient at the discrete level by the Karhunen-Loève-like expansion}
\begin{equation*}
	\nu(x, \theta) = \sqrt{2} + \sum_{i=1}^{n_{\theta}} \theta_i \sqrt{\lambda_i} \xi_i(x),
\end{equation*}
where $\xi_i$, $i= 1, \ldots, n_{\theta}$ are the first $n_{\theta}$ ordered eigenfunctions of 
\begin{equation*}
	C_{ij} = a \exp\left(  - \frac{\norm{x^i - x^j}}{2 b^2} \right)	+ c \delta_{ij}, \quad 1 \leq i,j \leq d,
\end{equation*}
where $\delta_{ij}$ is a Kronecker symbol, $a = 1$, $b = 1$, $c = 0.1$ and $\{x^i\}_{i=1}^{d} \subset \Tcal_h$ are the mesh nodes.
The eigenfunctions $\xi_i$ are depicted in~\Cref{fig:eigenfunctions}. 
In order to ensure $\nu(x, \theta) > 0$, the parameters are constrained to the hypercube
\begin{equation*}
	\Qcal_{n_{\theta}} =	\left\{ \theta \in \R^{n_{\theta}} \st \max_{i=1, \ldots, n_{\theta}} |\theta_i| \leq \sqrt{2} \left( \sum_{i=1}^{n_{\theta}} \sup_{x \in D} \sqrt{\lambda_i} |\xi_i(x)| \right)^{-1} \right\}. 
\end{equation*}
In this example, we consider $n_\theta = 6$, and \moda{set} 
\begin{equation}
	\moda{\theta^{\true} = 
		\left(0.271,\, 0.266,\, 0.504, \, -0.111, \, -0.014, \, -0.086 \right)^{\top}
\in \Qcal_{n_{\theta}}.}
\end{equation}
Then, we sample a perturbed $\theta^{\mathrm{perturbed}} \sim \Ncal(\theta^{\true}, 0.05^2 \bI_{n_{\theta}})$ and set the initial conditions of the ensemble as $\htheta_0^{(i)} \sim \Ncal(\theta^{\mathrm{perturbed}}, 0.05^2 \bI_{n_{\theta}})$ for $i = 1, \ldots, {\rc{N}} = 200$. 
\moda{The ensemble of state particles} \modb{is} \moda{initialised with the same initial condition $\Xcal_0^{(i)} = I_h(\exp( -(x_1 - 1.5)^2 - 50 x_2^2 ))$} \modb{for all the particles} \moda{$i =1, \ldots, {\rc{N}}$, where $I_h$ is the FE interpolation operator.}

We consider two observation scenarios, full ($h^{\mathrm{full}} = id$) and partial observations.  
The partial observation operator $h^{\mathrm{part}}$ is defined as 
\begin{equation*}
	(h^{\mathrm{part}} u(x))_i = \frac{30}{0.05 \pi} \int_{D}  \exp \left( - \frac{ (x - x_i)^2 + (y - y_i)^2}{ 2 (0.05)^2} \right) u(x) \md x \md y,
\end{equation*}
for $(x_i, y_i) = (\rho \cos \alpha, \rho \sin \alpha)$ and where $\rho \in \{1, 1.5\}$ and  $\alpha \in \{\pi/2, \pi/3, \pi/4, \pi/6\}$ for all possible combinations.
Consequently, $h^{\mathrm{part}}u \in \R^{8}$.
\moda{The true signal is} \modb{generated} \moda{using an explicit Euler scheme with the time-discretisation $\delt = 4.4 \cdot 10^{-5}$, the observations are computed for that time-discretisation.
	Additionally, in either cases of full ($k = d$) or partial} \modb{($k = 8$)} \moda{observations, we set $\Gamma = \gamma \bI_{k}$ with $\gamma = 10^{-8}$.
}


\rb{We do not employ hyper-reduction in this example. Since the nonlinearities in the augmented state are quadratic and the ranks considered below are small, the reduced nonlinear terms can be evaluated efficiently without hyper-reduction, and applying DEIM/CUR brings only marginal gains. The purpose of this example is to validate, as a first step, the \dlr-\enkf~formulation itself, rather than the DEIM/CUR hyper-reduction method.}

\begin{figure}[htbp]
\centering
\subfloat[True]{\includegraphics[width=0.32\linewidth]{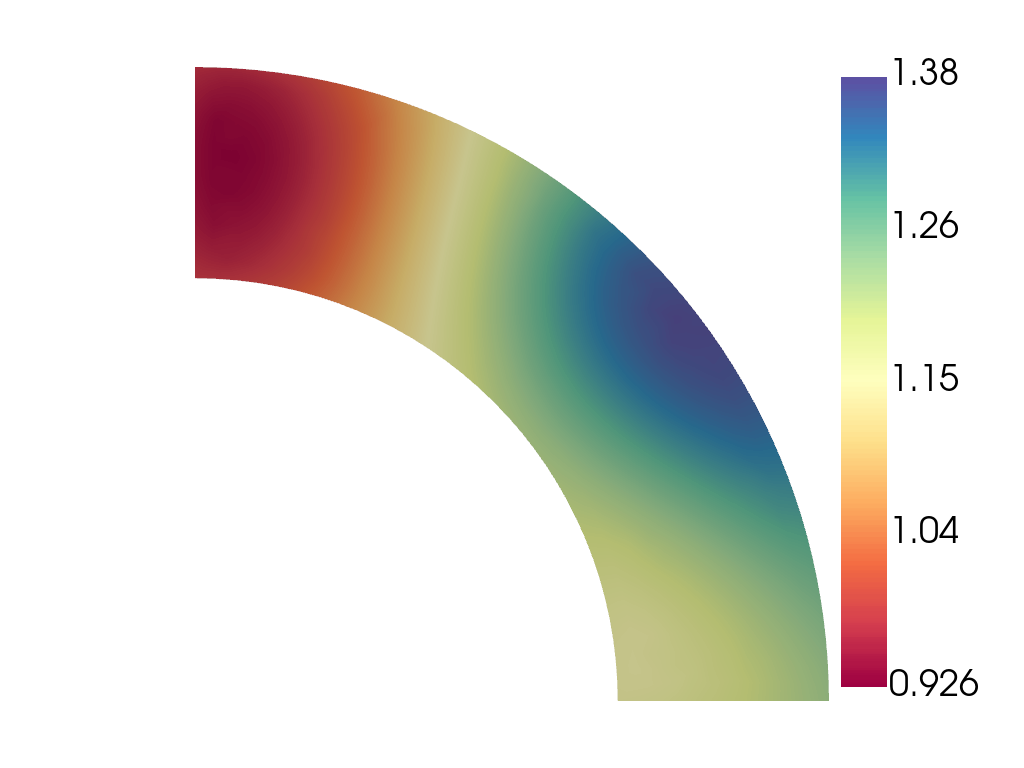}} 
\subfloat[\senkf]{\includegraphics[width=0.32\linewidth]{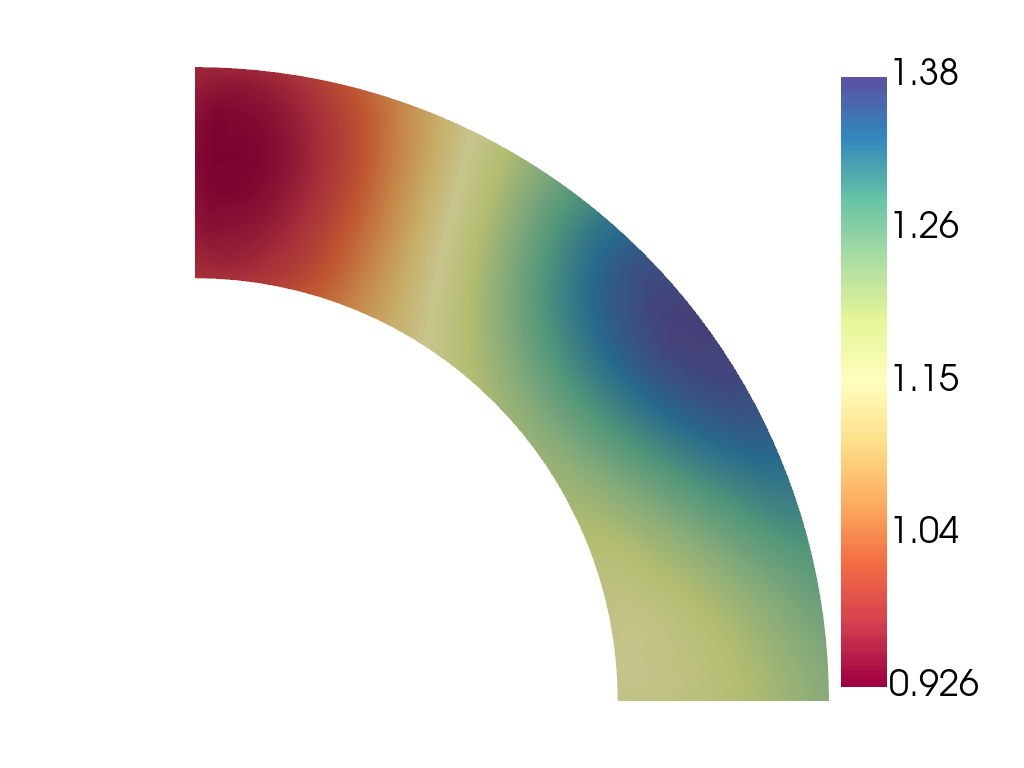}} 
\subfloat[\dlr-\senkf\ ($R = 7$)]{\includegraphics[width=0.32\linewidth]{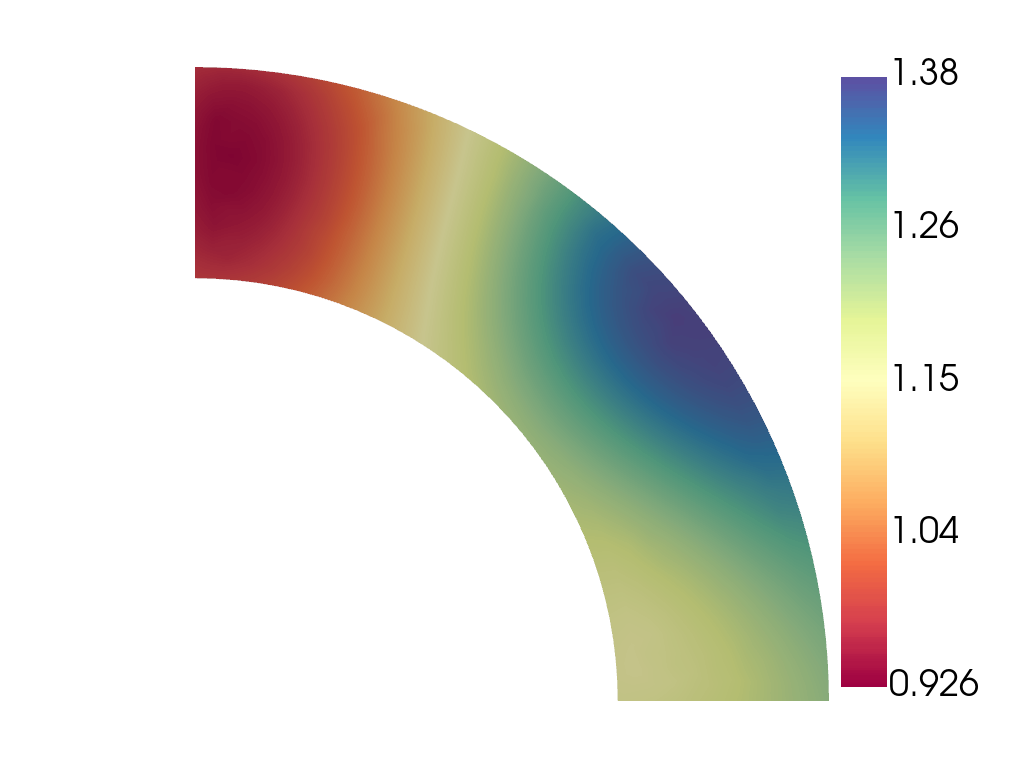}} 
\caption{\ra{Fisher-KPP: true and identified diffusion coefficient $\nu(x, \theta)$ at final time (full observations)}}
\label{fig:diffusion-coefficient}
\end{figure}

\begin{figure}[htbp]
\centering
\subfloat[$\theta_1$]{\includegraphics[width=0.32\linewidth]{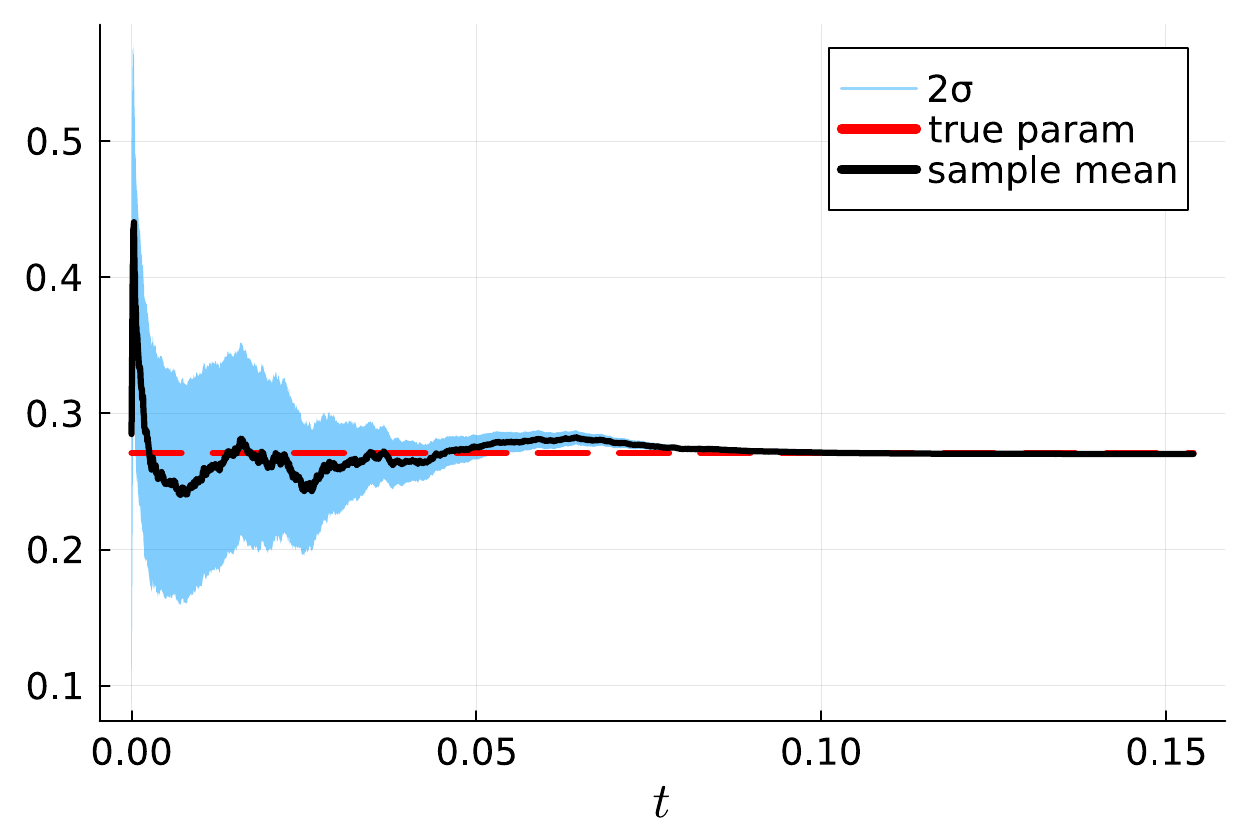}} 
\subfloat[$\theta_2$]{\includegraphics[width=0.32\linewidth]{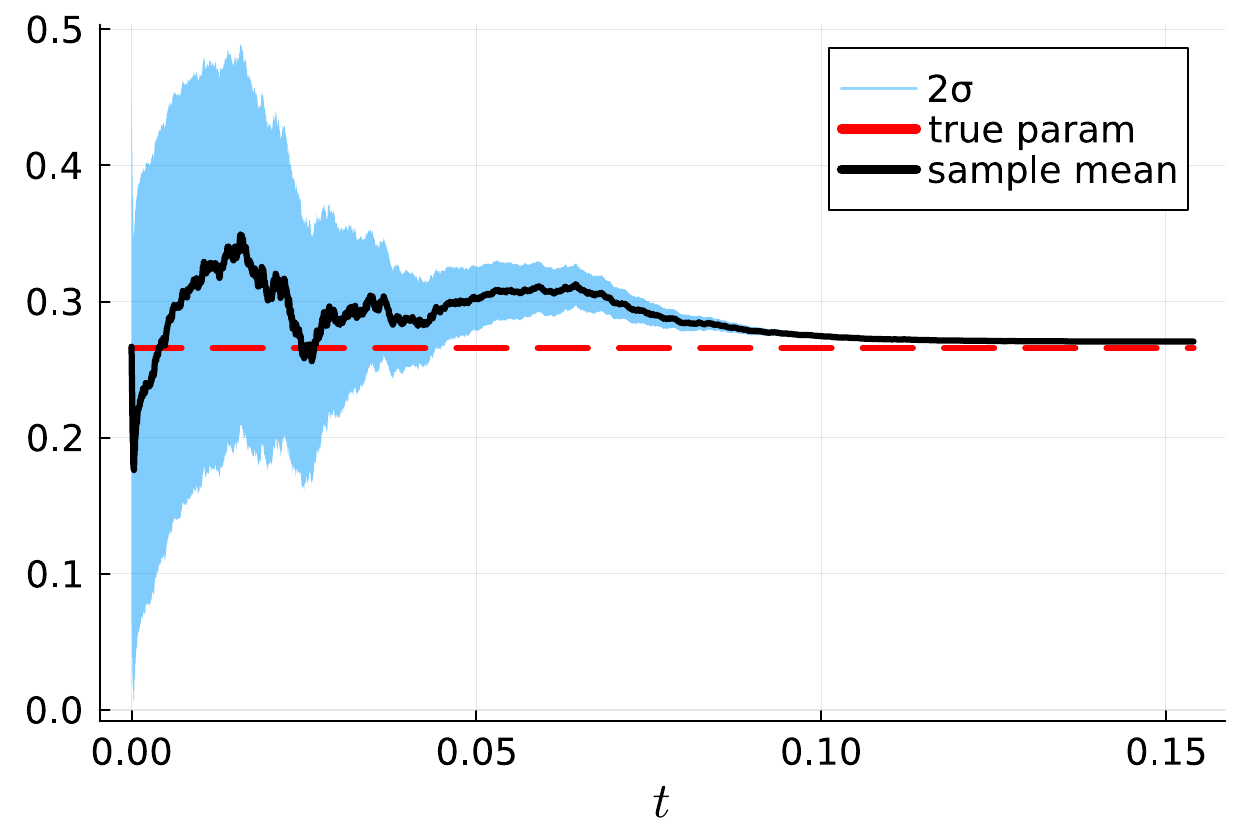}} 
\subfloat[$\theta_3$]{\includegraphics[width=0.32\linewidth]{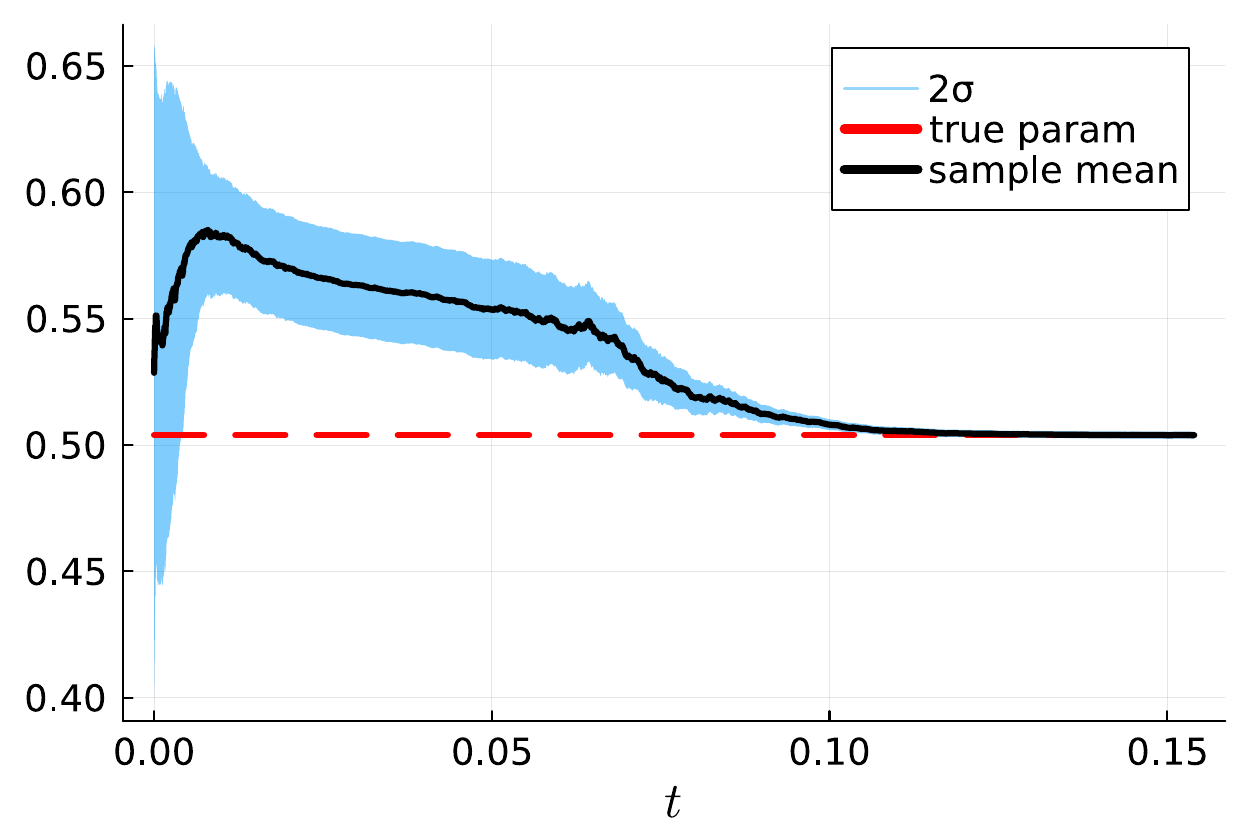}} \\
\subfloat[$\theta_4$]{\includegraphics[width=0.32\linewidth]{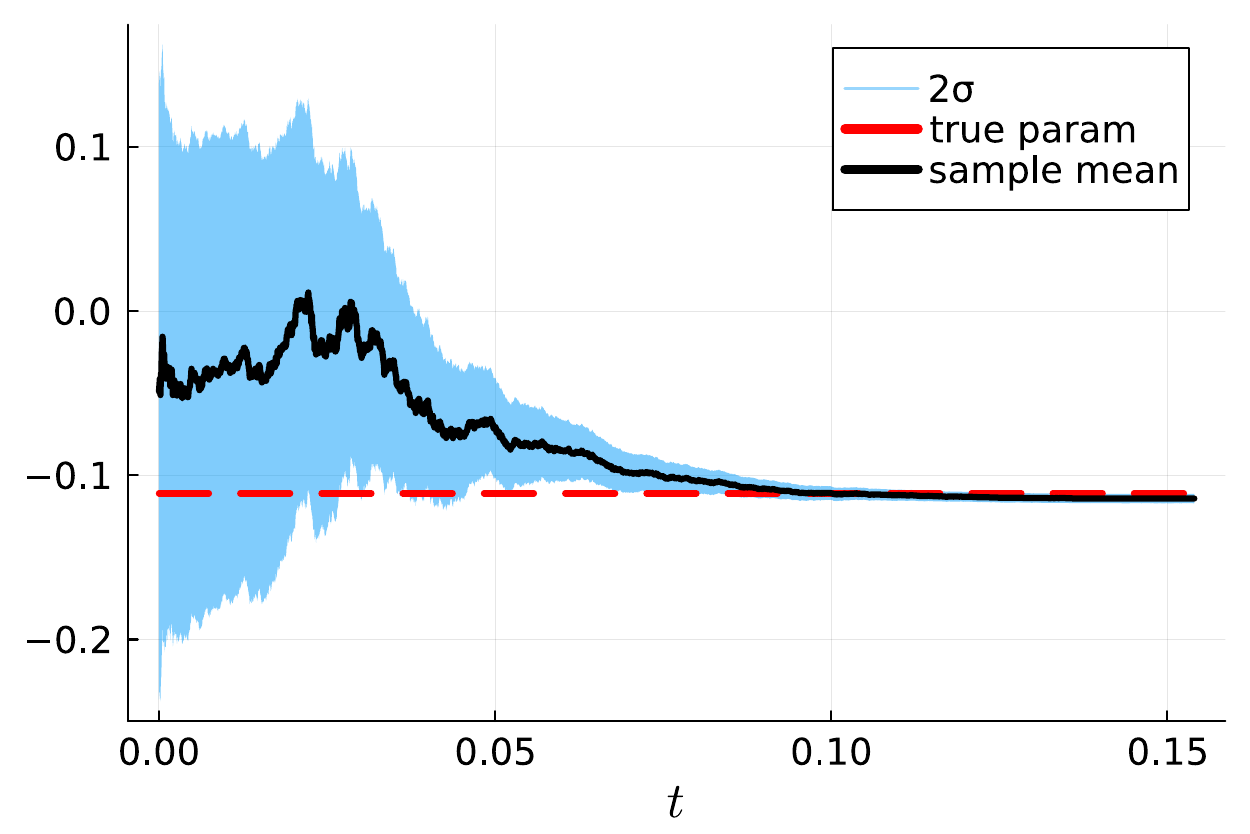}} 
\subfloat[$\theta_5$]{\includegraphics[width=0.32\linewidth]{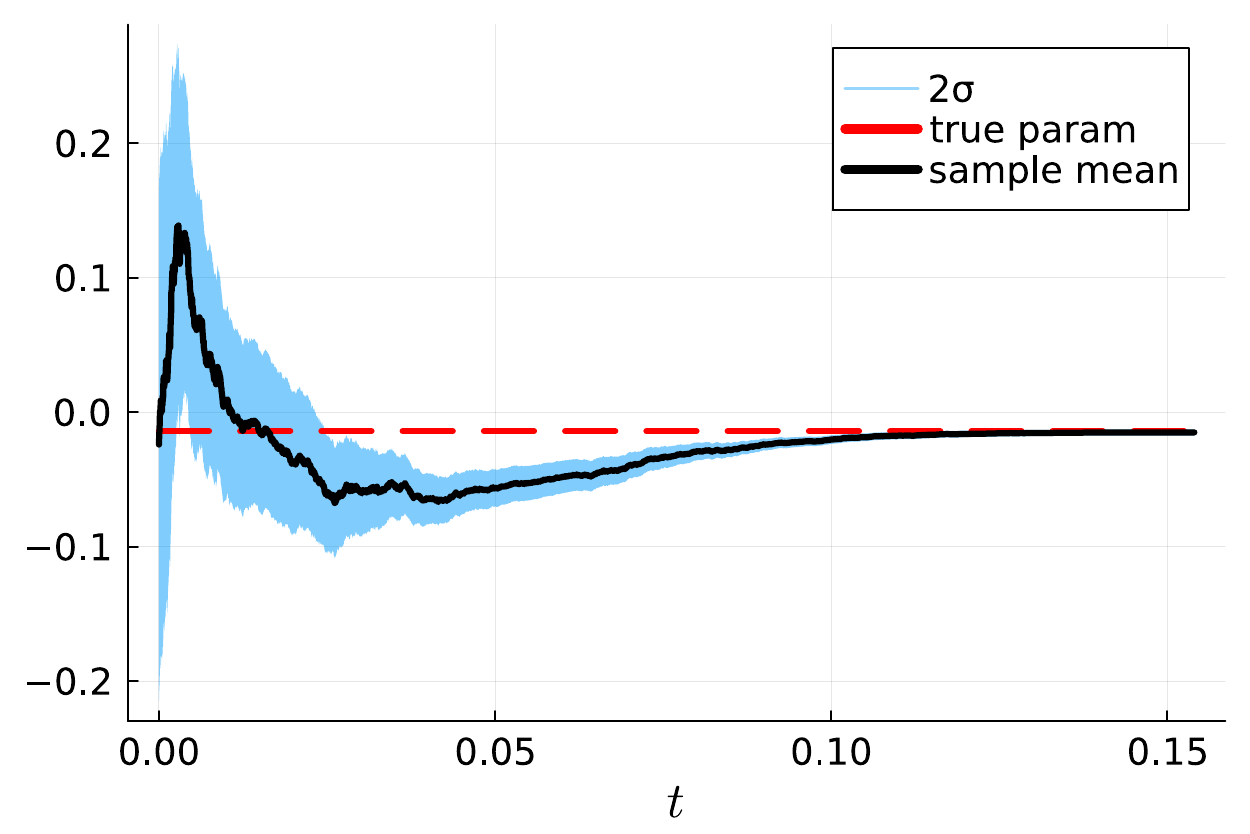}} 
\subfloat[$\theta_6$]{\includegraphics[width=0.32\linewidth]{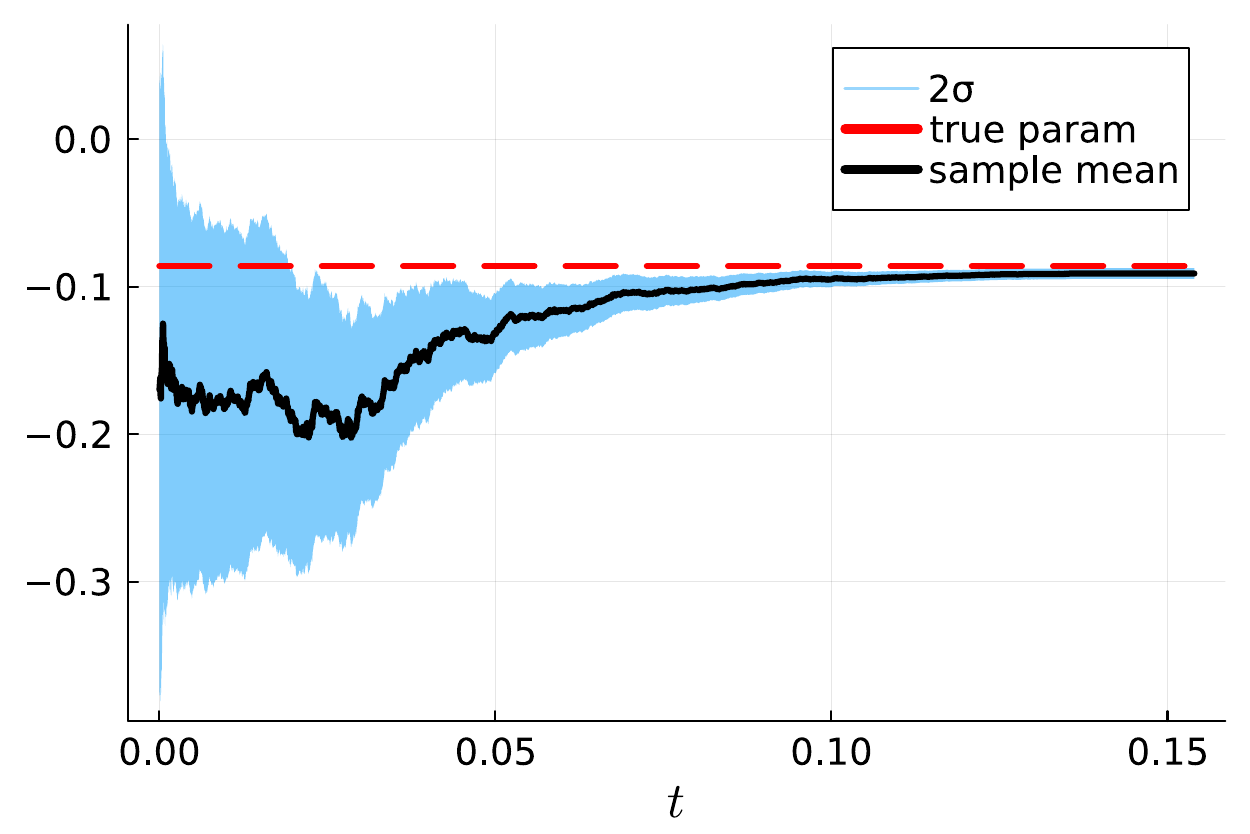}} 
\caption{\ra{Fisher-KPP: Full Order Model \senkf~parameter identification with ${\rc{N}} = 200$ particles (full observations)}}
\label{fig:fom-parameter-identification}
\end{figure}

\begin{figure}[htbp]
\centering
\subfloat[$\theta_1$]{\includegraphics[width=0.32\linewidth]{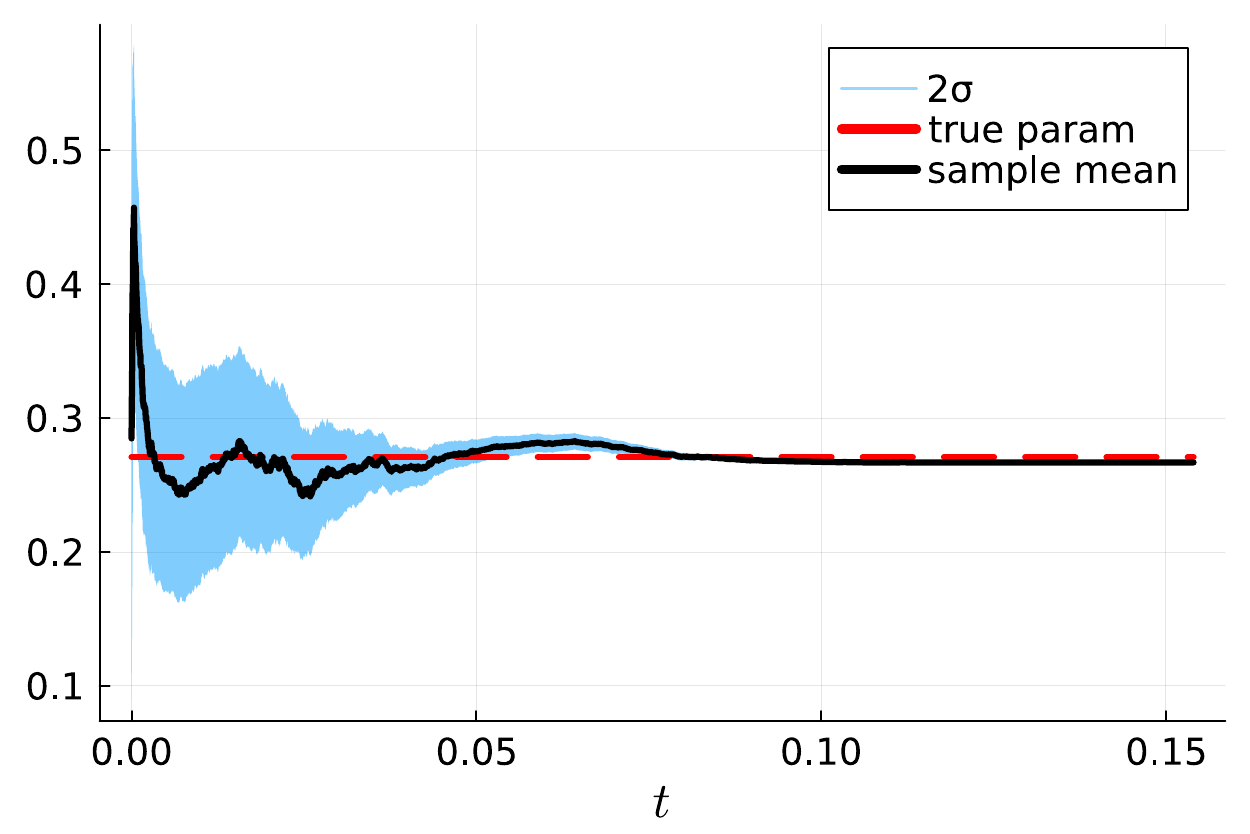}} 
\subfloat[$\theta_2$]{\includegraphics[width=0.32\linewidth]{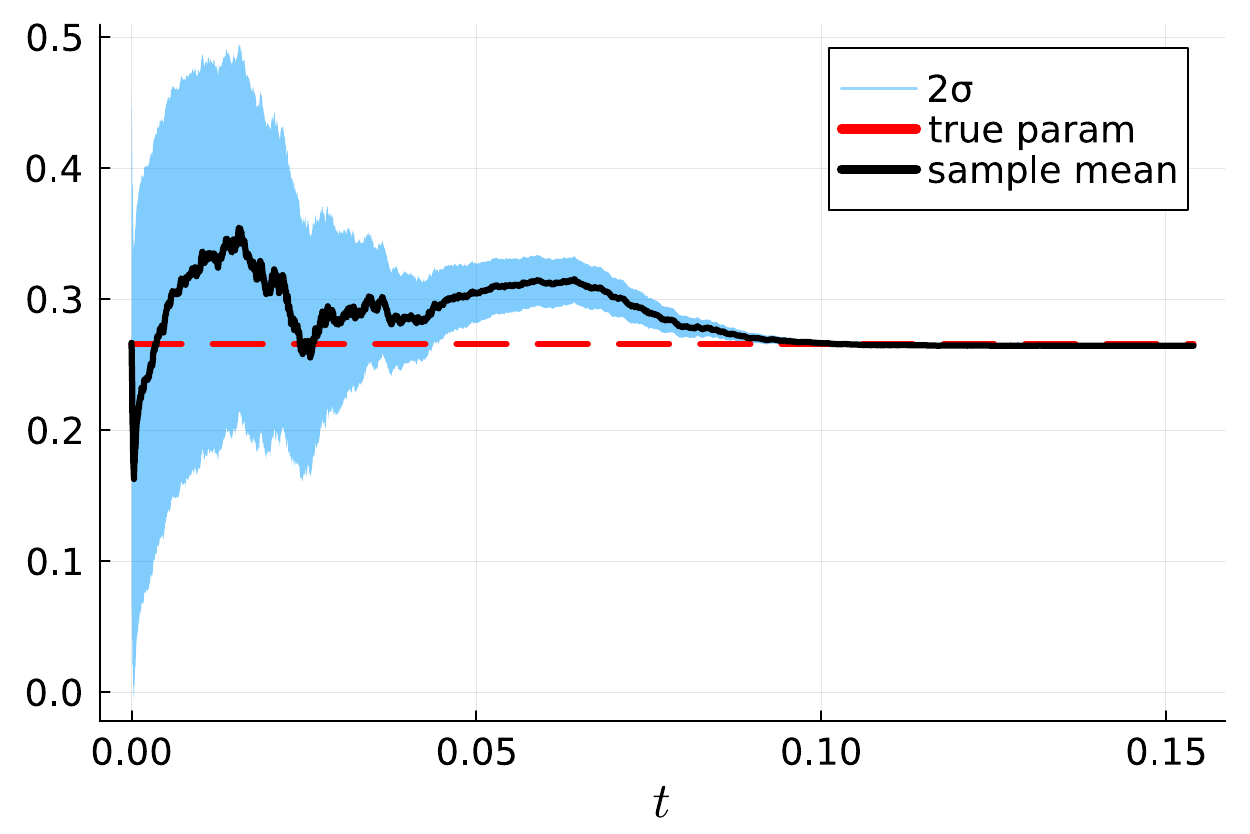}} 
\subfloat[$\theta_3$]{\includegraphics[width=0.32\linewidth]{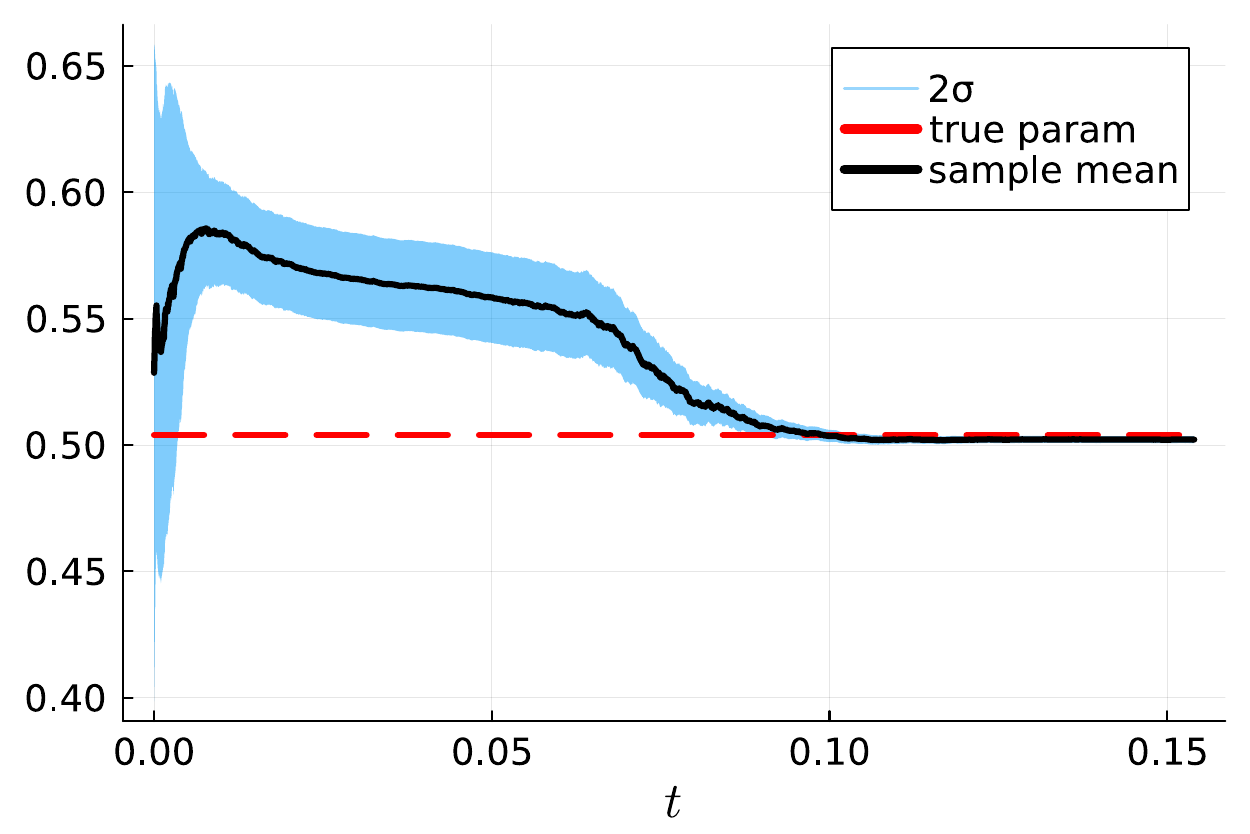}} \\
\subfloat[$\theta_4$]{\includegraphics[width=0.32\linewidth]{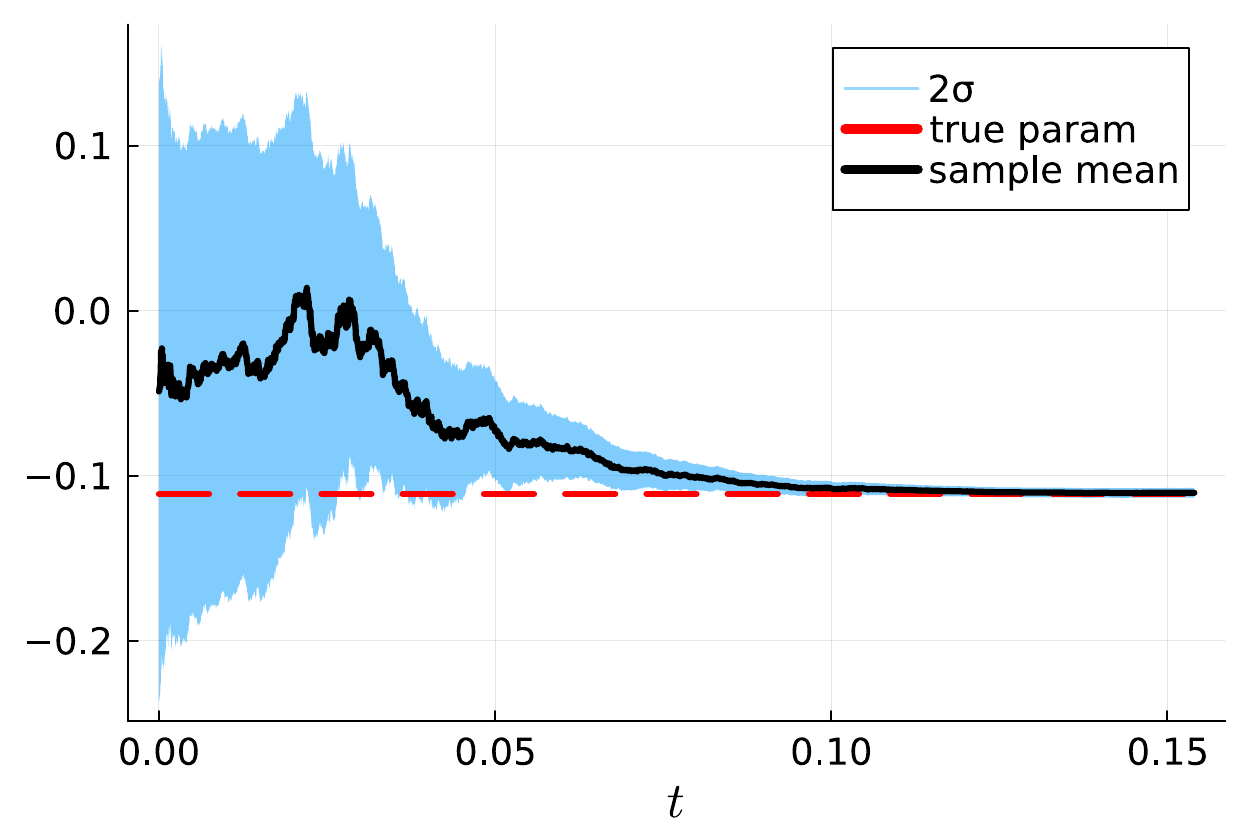}} 
\subfloat[$\theta_5$]{\includegraphics[width=0.32\linewidth]{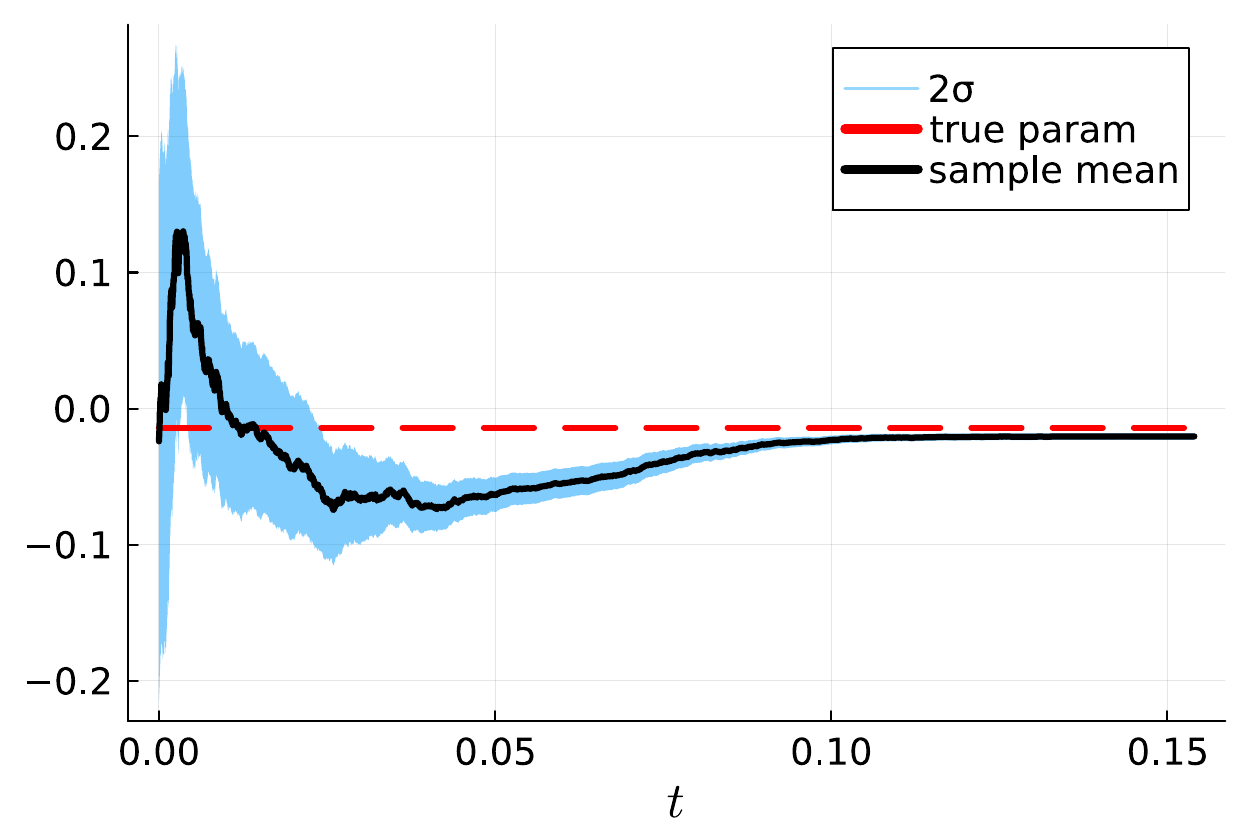}} 
\subfloat[$\theta_6$]{\includegraphics[width=0.32\linewidth]{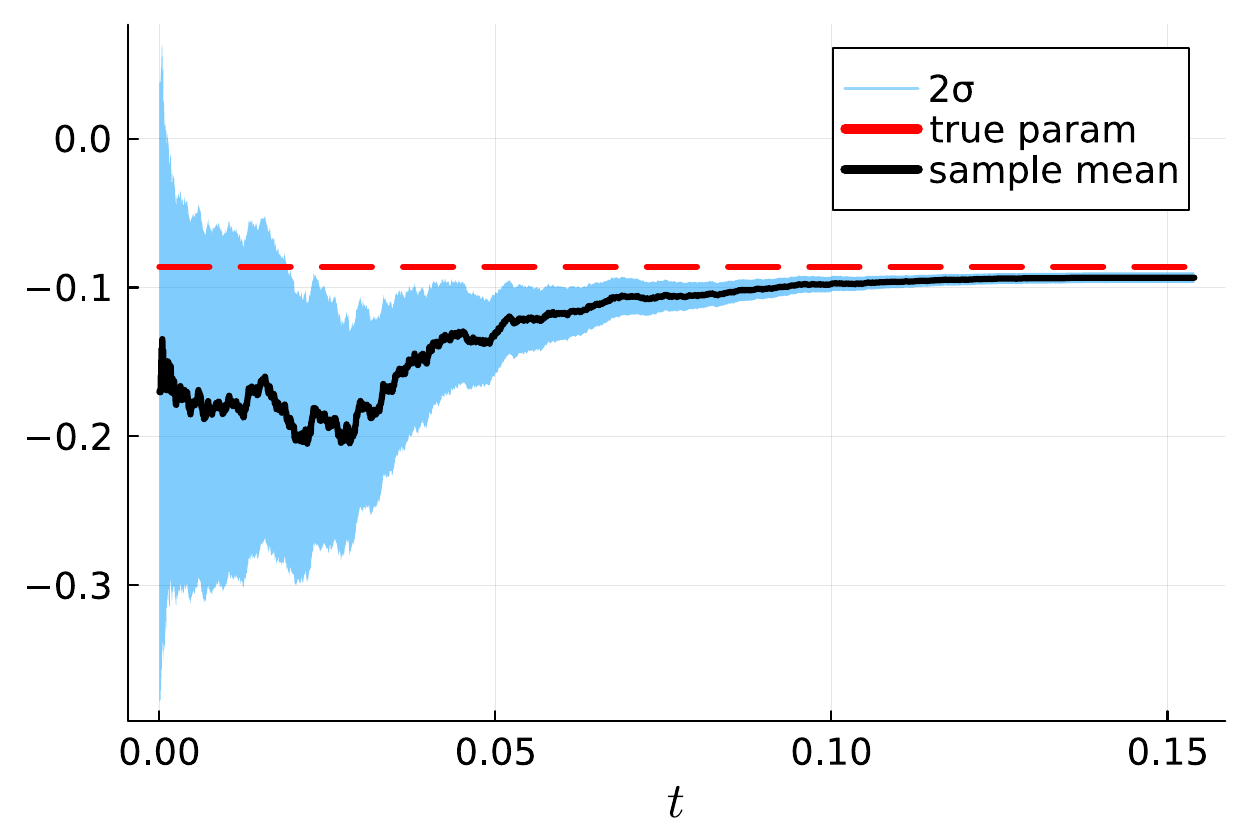}} 
\caption{\ra{Fisher-KPP: \dlr-\senkf~($R = 7$) parameter identification with ${\rc{N}} = 200$ particles (full observations)}}
\label{fig:dlr-parameter-identification}
\end{figure}

For both observation scenarios, we compare the \venkf, \denkf~and \senkf~to their \dlr~counterparts for varying ranks $R = 2, 5, 7$. 
\moda{For the \dlr-\enkf\ schemes, we used the explicit BUG scheme combined with the analysis step as outlined in~\cref{sec:forecast-analysis-time-discretisation}. 
	To keep the comparisons fair, all of the other parameters are identical to the \enkf\ case, in particular using ${\rc{N}} = 200$ particles and $\delt= 4.4 \cdot 10^{-5}$.}
In~\Cref{fig:diffusion-coefficient} we compare the true field $\nu(x)$ to the one obtained by the \senkf~algorithm, as well as the~\dlr-\senkf. 
Their similarity is further verified when examining the evolution of the ensemble of parameters. 
\Cref{fig:fom-parameter-identification} displays the parameter identification achieved by the \senkf. 
The filter correctly identifies all the parameters (with minor biases for $\theta_2$, $\theta_4$, $\theta_6$). 
Similarly, \Cref{fig:dlr-parameter-identification} displays the same quantities when using the \dlr-\senkf~filter with rank $7$.
Observe that the plots are very similar, suggesting that this rank is already enough to suitably capture the filtering density as it evolves in time.


\begin{figure}
	\subfloat[][$R = 2$]{\includegraphics[width=.24\textwidth]{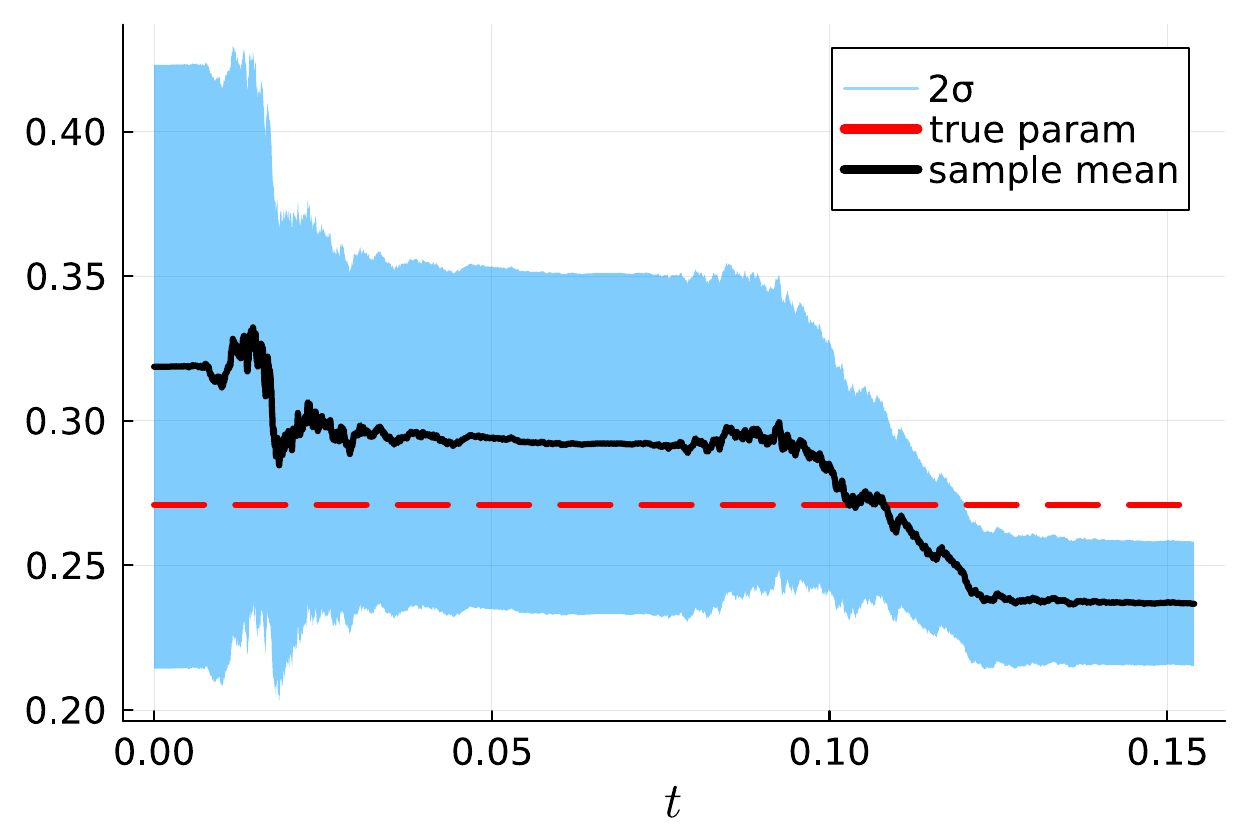}}
	\subfloat[][$R = 5$]{\includegraphics[width=.24\textwidth]{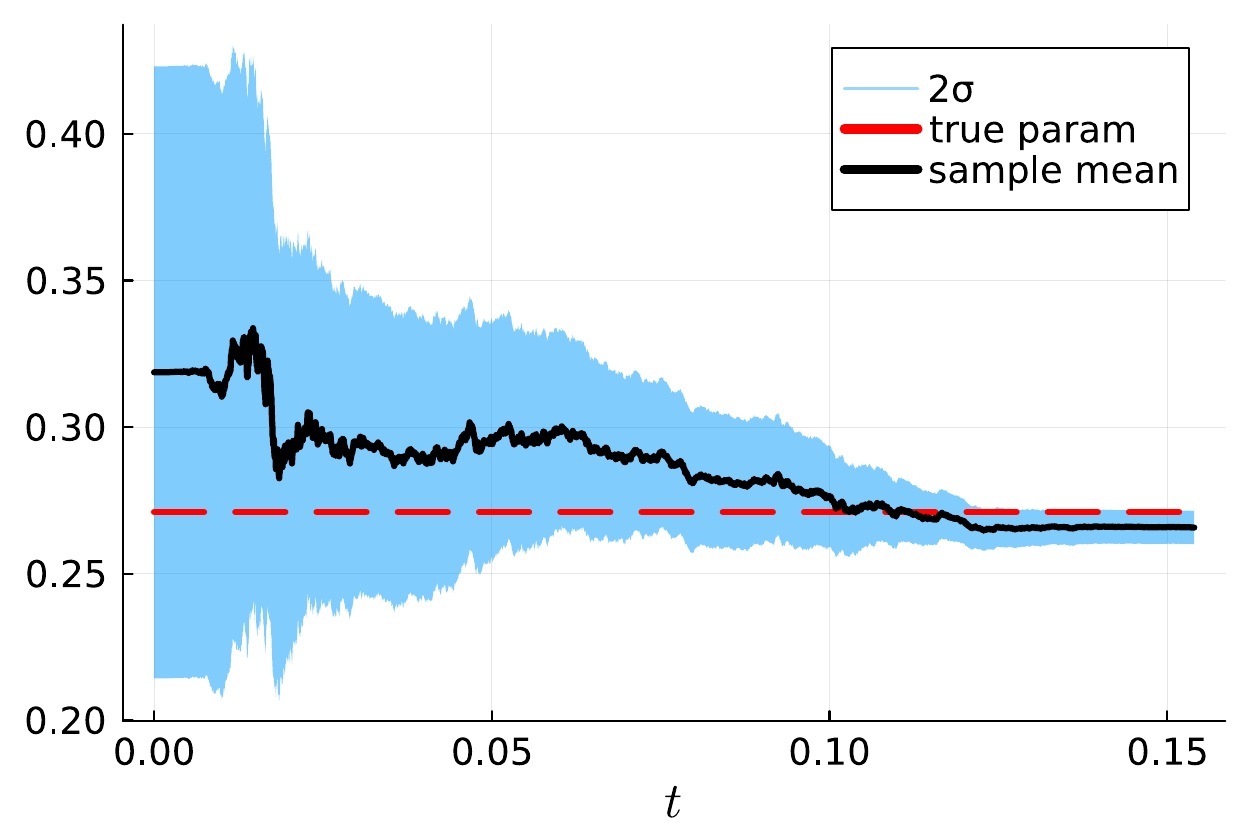}} 
	\subfloat[][$R = 7$]{\includegraphics[width=.24\textwidth]{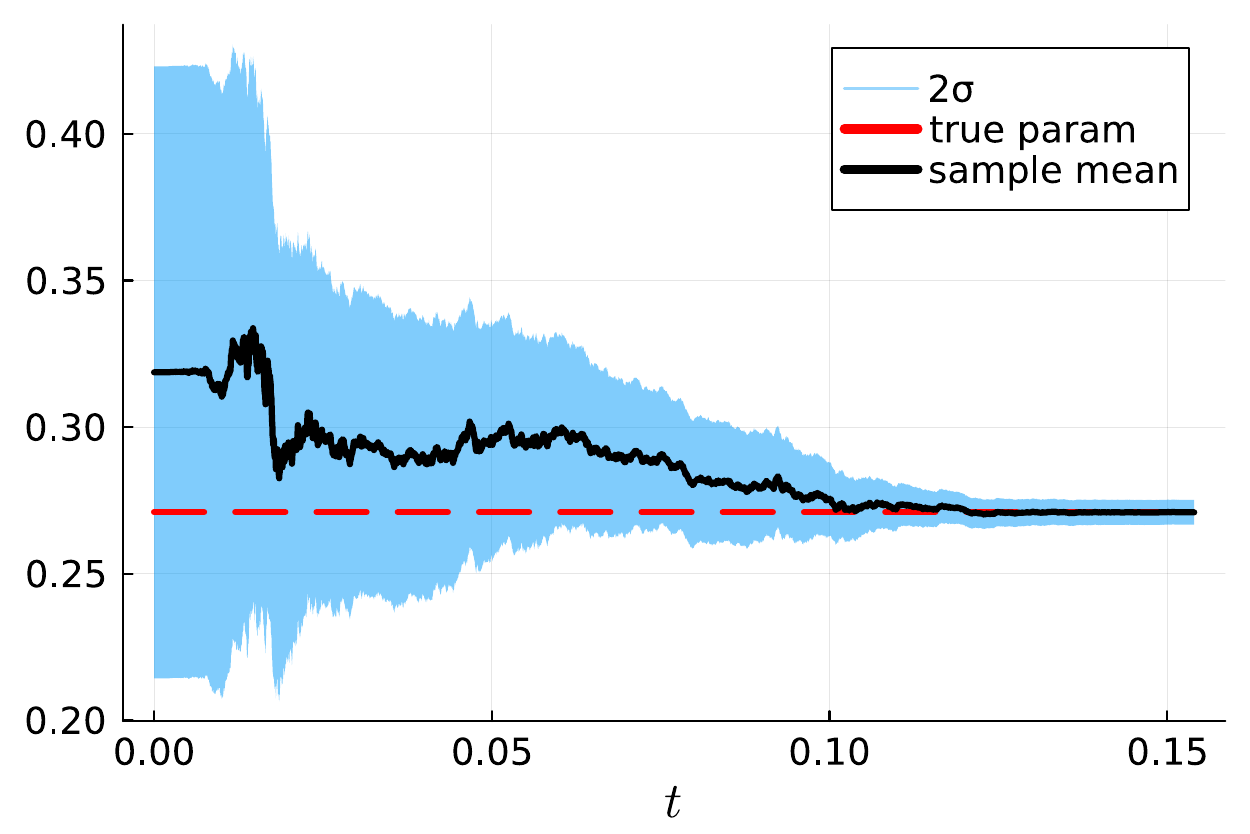}}
	\subfloat[][FOM]{\includegraphics[width=.24\textwidth]{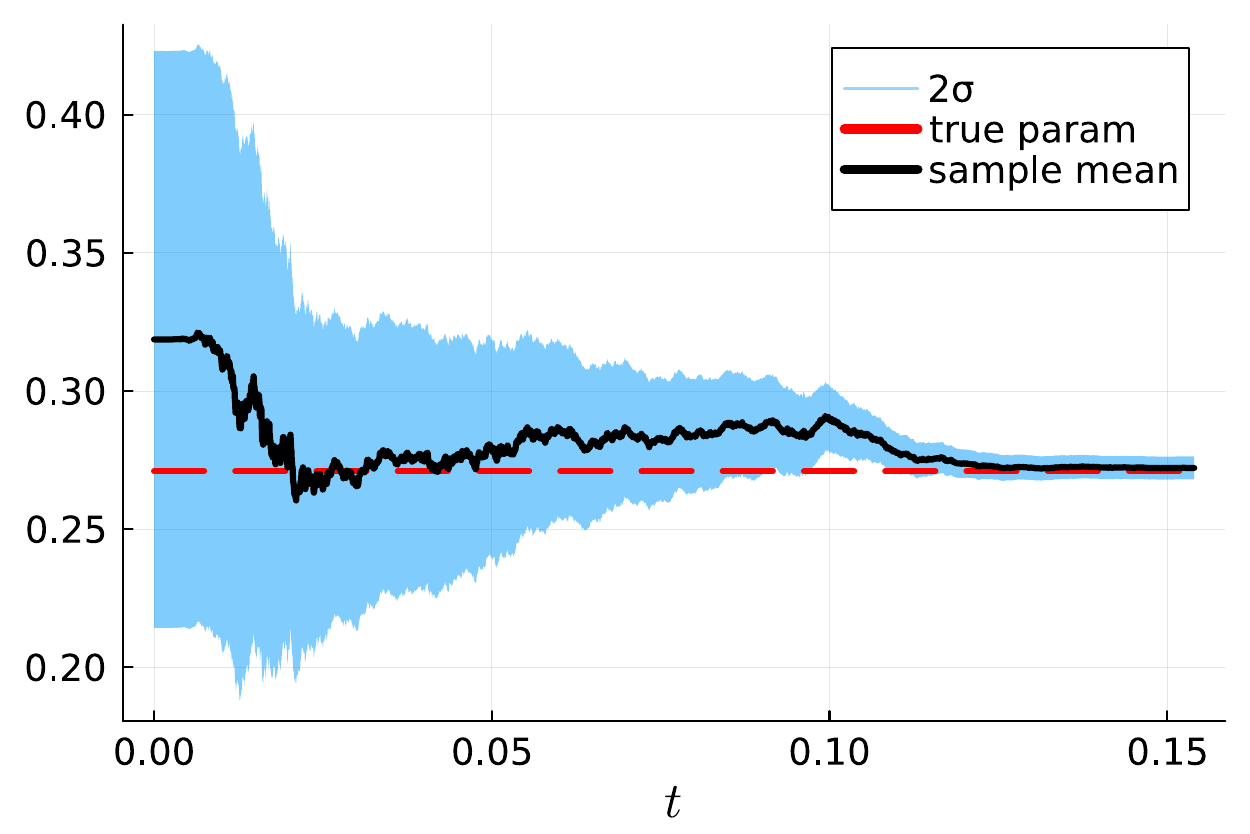}}
	\caption{\ra{Fisher-KPP: $\theta_1$ for partial observations (${\rc{N}} = 200$); (a)--(c) \dlr-\senkf\ (d) \senkf}}
	\label{fig:partial-observations-theta1}
\end{figure}

\begin{figure}
	\subfloat[][$R = 2$]{\includegraphics[width=.24\textwidth]{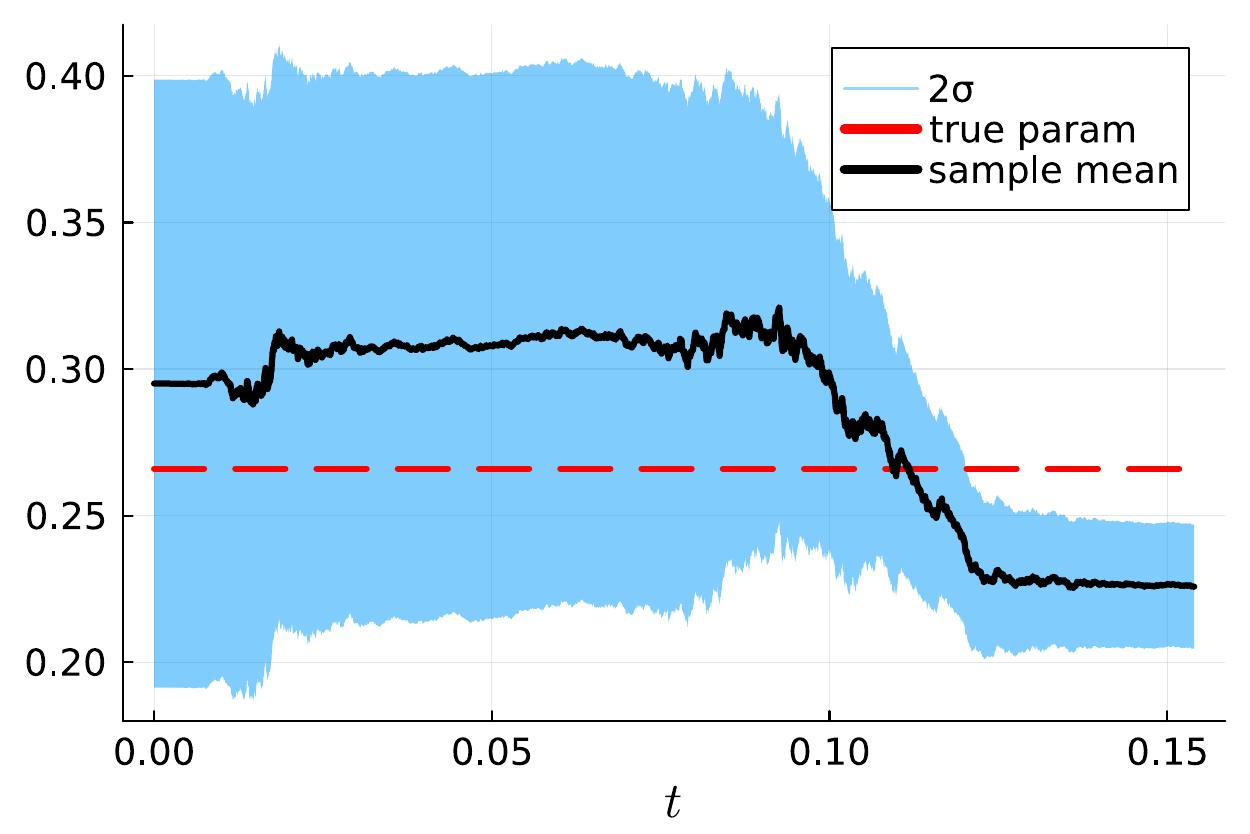}}
	\subfloat[][$R = 5$]{\includegraphics[width=.24\textwidth]{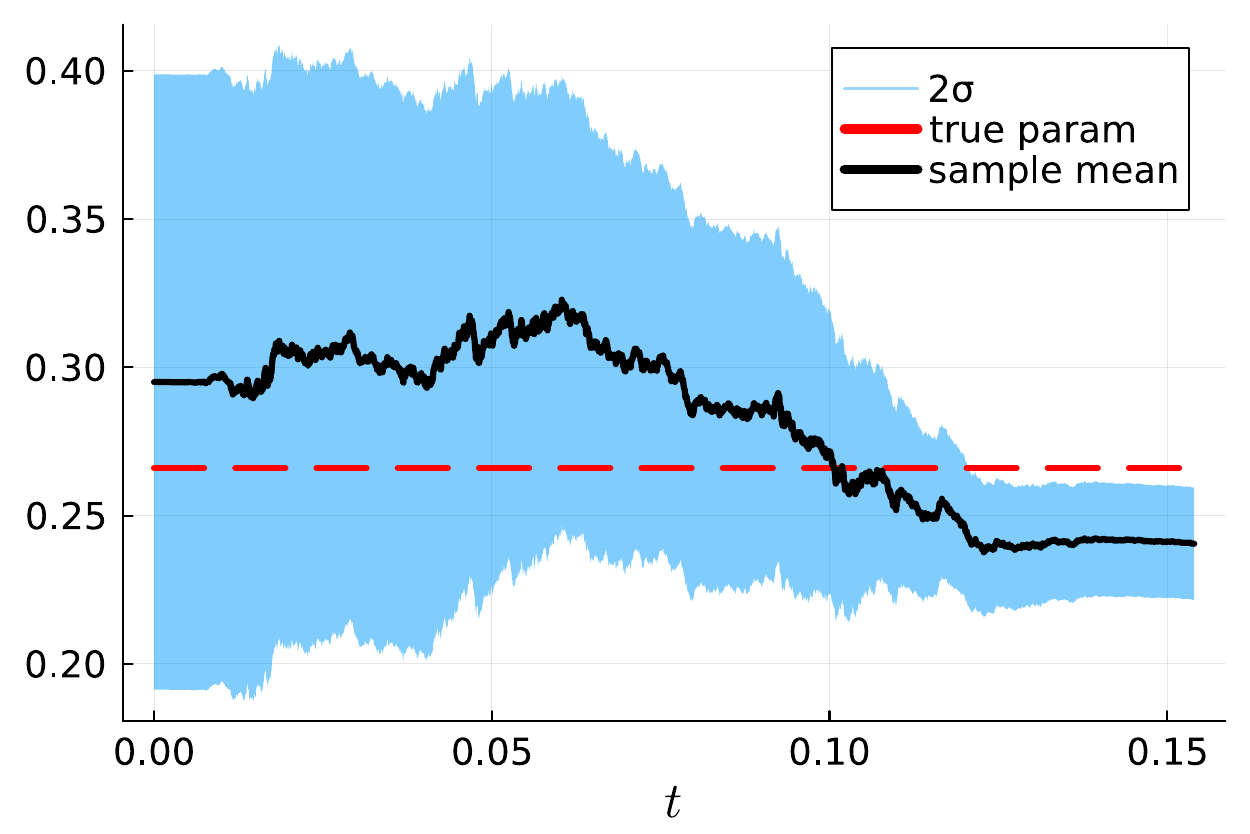}} 
	\subfloat[][$R = 7$]{\includegraphics[width=.24\textwidth]{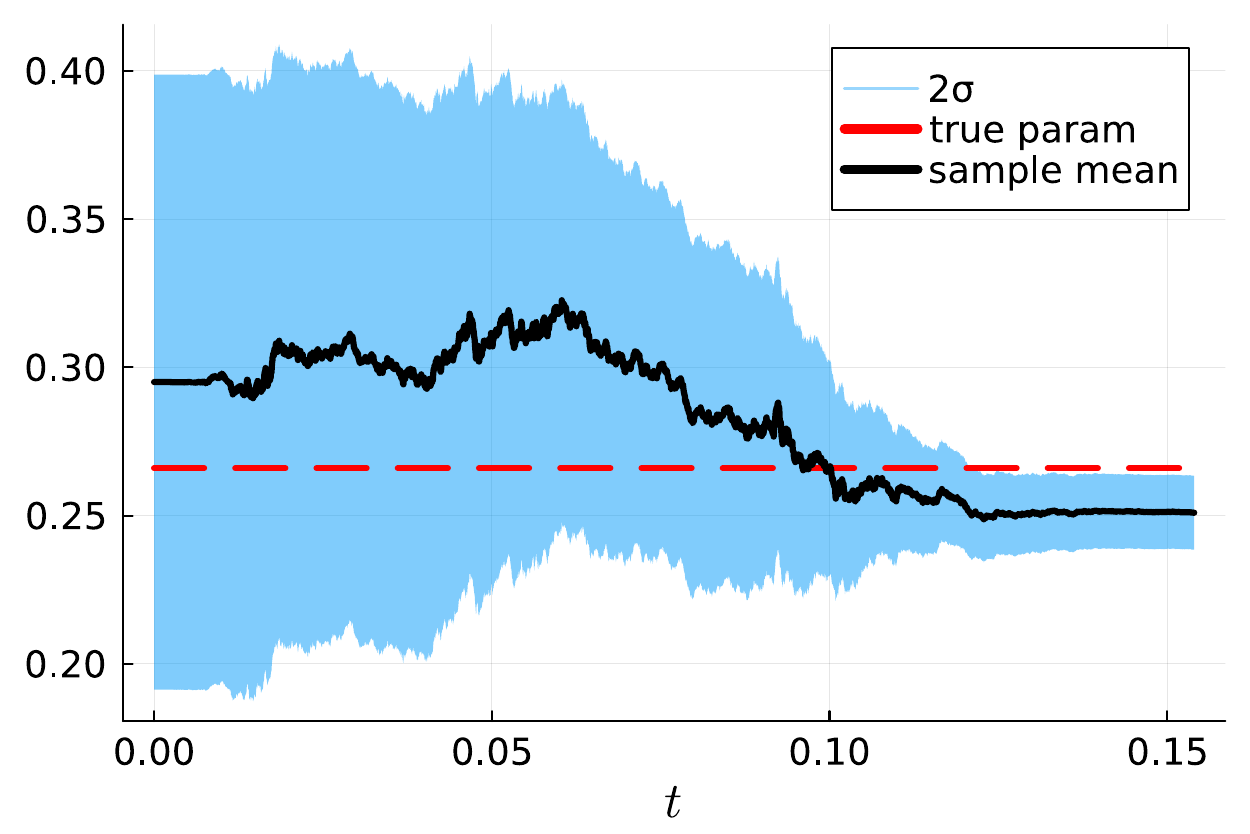}}
	\subfloat[][FOM]{\includegraphics[width=.24\textwidth]{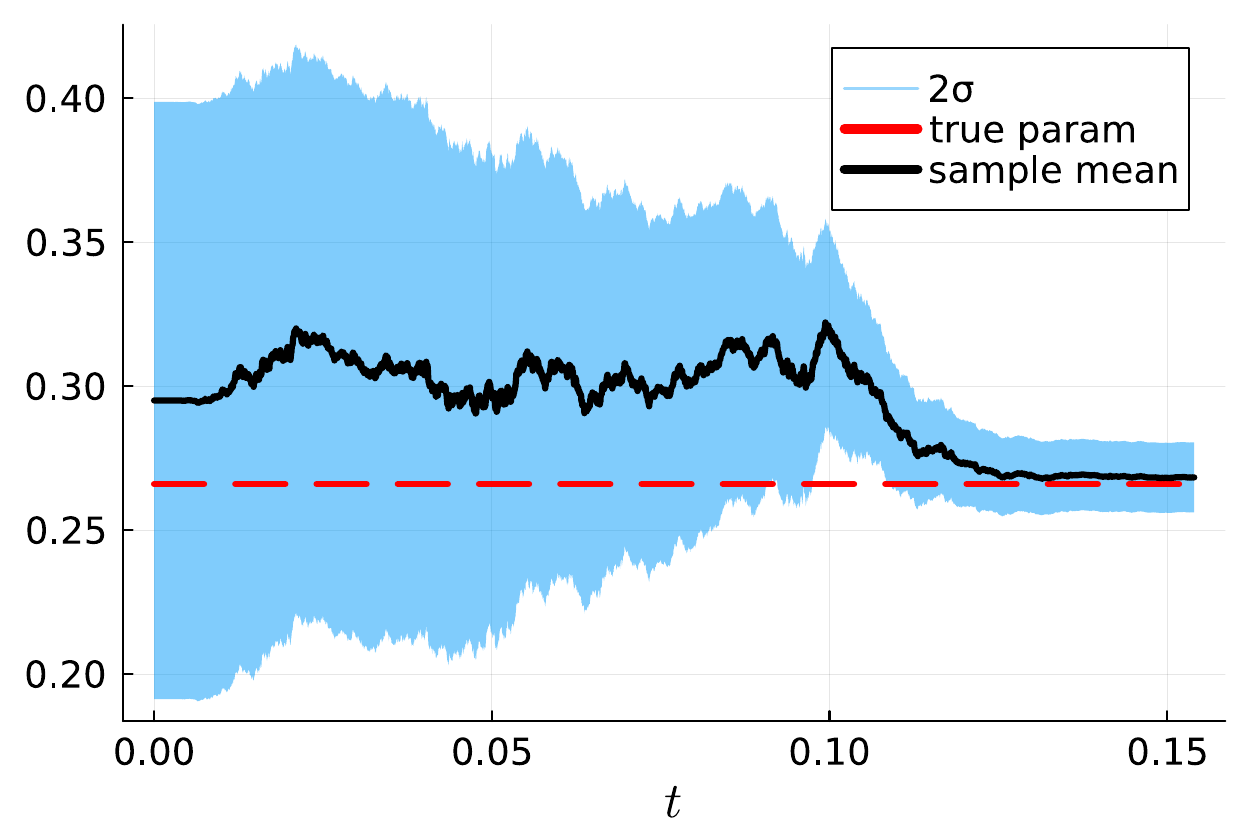}}
	\caption{\ra{Fisher-KPP: $\theta_2$ for partial observations (${\rc{N}} = 200$); (a)--(c) \dlr-\senkf\ (d) \senkf}}
	\label{fig:partial-observations-theta5}
\end{figure}

The relative errors $\norm{\E_{\rc{N}}[{\htheta}_{T}] - \theta^{\true}} / \norm{\theta^{\true}}$ are displayed in~\Cref{tab:rel-errs} for all combinations of methods \moda{(\senkf, \venkf, \denkf)} and both observation scenarios \moda{(averaged over 10 runs)}.
Throughout the methods, it is seen that increasing the rank correlates with a decreasing relative error.
In the case of full observation, the relative error approaches that obtained by the full order model. 
\moda{
	It is noteworthy that in the case of partial observations, the \enkf\ is (at least initially) not well-approximated by a subspace of rank $7$ or less. 
	This is indirectly seen in~\cref{fig:partial-observations-theta1,fig:partial-observations-theta5}, as the evolution of the ensemble of parameter of $\theta_1$ \modb{and $\theta_2$} for the \dlr-\senkf\ are distinct from the \senkf. 
\modb{Despite} developing a distinct solution, the \dlr-\enkf\ correctly} identifies \modb{$\theta_1$, but fails for $\theta_2$}.
	\modb{This showcases two effects: (a) the \dlr-\enkf\ distinct from the true solution \emph{can} still perform an adequate signal-tracking and parameter identification (notably, it was observed several times in simulations that the relative errors of the \dlr-\enkf\ filters were eventually \emph{lower} than their \enkf\ counterparts -- not shown in figures); (b) however, augmenting the rank is sometimes necessary to ensure consistency with the underlying \enkf\ (and in this case, ensure successful parameter identification).}

Additionally, the timings are displayed in~\Cref{tab:timings}. 
Thanks to the low-rank structure, the \moda{\dlr-\enkf\ versions} are always faster than their full order \enkf\ counterparts, \moda{and, predictably,} the computational cost increases with the rank. 
It is expected that speed-ups can be even more important in the case of stronger nonlinearities, \moda{much finer FE discretisations,} and combinations with the DEIM technique.

\begin{table}
	\centering
	\subfloat[Full observations]{
\begin{tabular}{@{}c|ccc|l@{}}
  \multicolumn{1}{l|}{} & \multicolumn{3}{c|}{DLR}                                     & FOM                 \\ \midrule
  \multicolumn{1}{l|}{} & $R = 2$              & $R =5$               & $R = 7$              &                     \\ \midrule
  \senkf & 0.585 & 0.013 & 0.014 & \multicolumn{1}{c}{0.012} \\ \midrule
  \venkf                & 0.660 & 0.018 & 0.011 & 0.008 \\ \midrule
  \denkf                & 0.275 & 0.014 & 0.008 & 0.006 \\ \bottomrule
\end{tabular}
} \quad\qquad
	\subfloat[Partial observations]{
\begin{tabular}{@{}c|ccc|l@{}}
            \multicolumn{1}{l|}{} & \multicolumn{3}{c|}{DLR}                                     & FOM                 \\ \midrule
            \multicolumn{1}{l|}{} & $R = 2$              & $R =5$               & $R = 7$              &                     \\ \midrule
            \senkf & 0.191 & 0.096 & 0.092 & \multicolumn{1}{c}{0.077} \\ \midrule
            \venkf                & 0.240 & 0.112 & 0.109 & 0.088 \\ \midrule
            \denkf                & 0.232 & 0.113 & 0.111 & 0.086 \\ \bottomrule
            \end{tabular}
}
  \caption{\ra{Fisher-KPP: relative errors $\norm{\E_{\rc{N}}[{\htheta}_{T}] - \theta^{\true}} / \norm{\theta^{\true}}$ (averaged over 10 runs)}}
	\label{tab:rel-errs}
\end{table}

\begin{table}
	\centering
	\subfloat[Full observations]{
            \begin{tabular}{@{}c|ccc|l@{}}
            \multicolumn{1}{l|}{} & \multicolumn{3}{c|}{DLR}                                     & FOM                 \\ \midrule
            \multicolumn{1}{l|}{} & $R = 2$              & $R =5$               & $R = 7$              &                     \\ \midrule
            \senkf & 4.7 & 9.2 & 14.4 & \multicolumn{1}{c}{18.2} \\ \midrule
            \venkf                & 4.7 & 9.3 & 13.8 & 18.6 \\ \midrule
            \denkf                & 4.4 & 8.6 & 13.9 & 17.6 \\ \bottomrule
            \end{tabular}
  } \quad\qquad
	\subfloat[Partial observations]{
            \begin{tabular}{@{}c|ccc|l@{}}
            \multicolumn{1}{l|}{} & \multicolumn{3}{c|}{DLR}                                     & FOM                 \\ \midrule
            \multicolumn{1}{l|}{} & $R = 2$              & $R =5$               & $R = 7$              &                     \\ \midrule
            \senkf & 5.7 & 11.2 & 15.0 & \multicolumn{1}{c}{19.3} \\ \midrule
            \venkf                & 4.3 & 9.1 & 14.4 & 18.4 \\ \midrule
            \denkf                & 4.2 & 9.0 & 13.9 & 18.0 \\ \bottomrule
            \end{tabular}
  }
  \caption{\ra{Fisher-KPP: timings [mins]}}
	\label{tab:timings}
\end{table}

%
%

\begin{figure}[htbp]
\centering
\subfloat[\ra{$t \in [0,0.154]$ (global)}]{\includegraphics[width=0.45\linewidth]{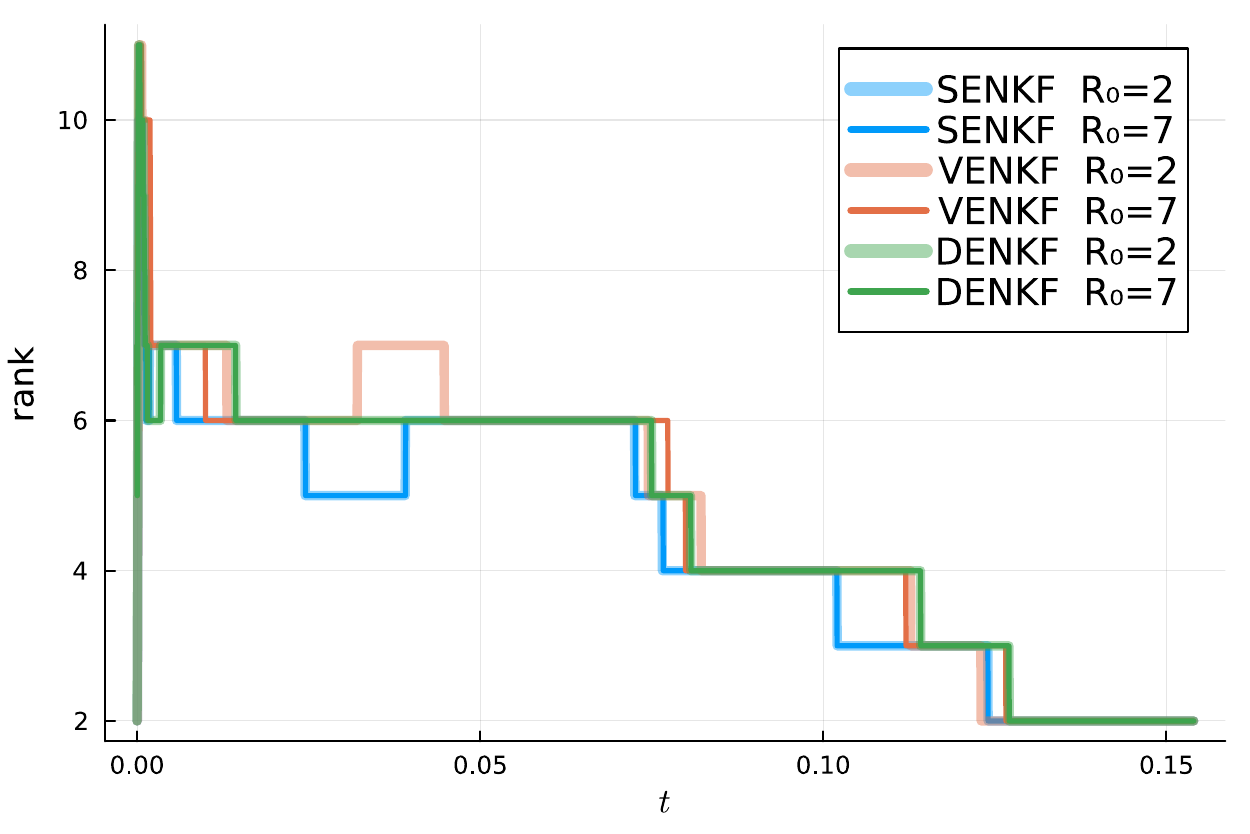}} 
\subfloat[\ra{$t \in [0, 2.64 \cdot 10^{-3}]$ (first 60 steps)}]{\includegraphics[width=0.45\linewidth]{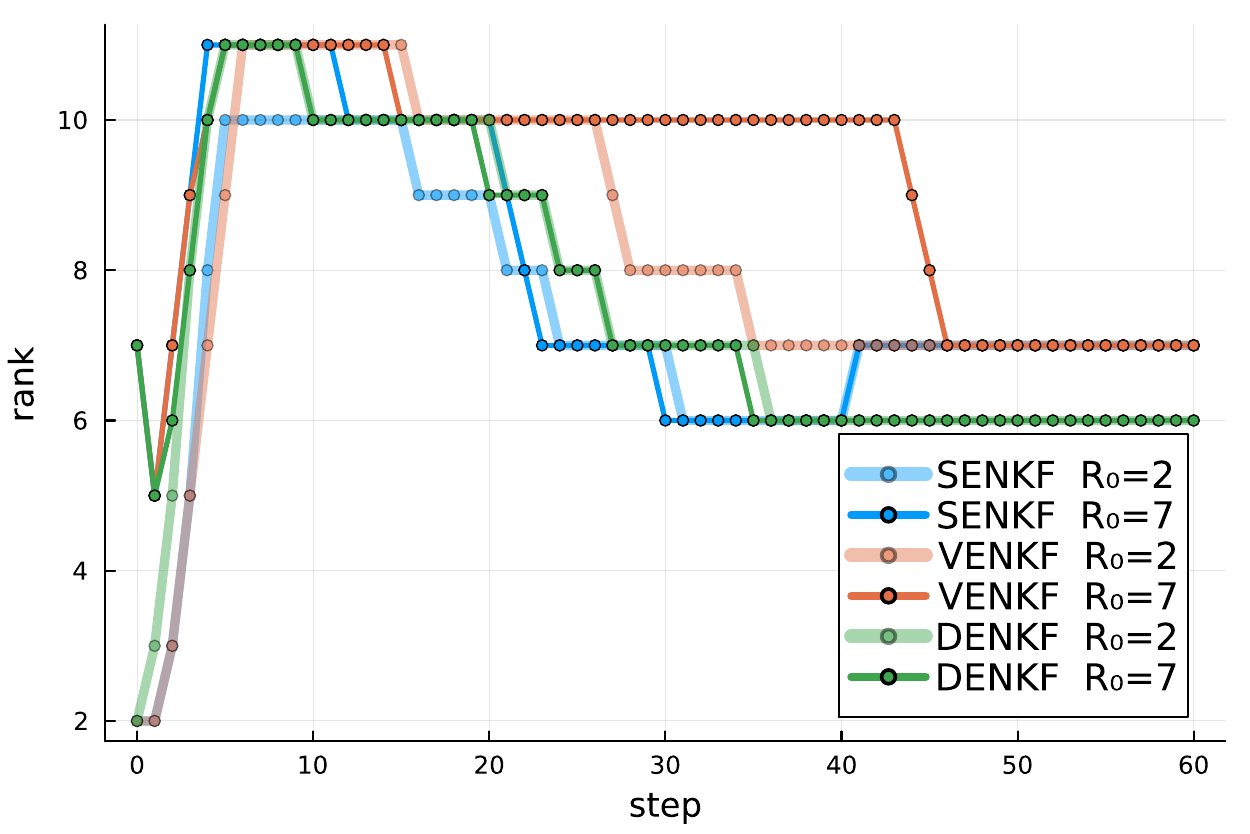}} 
\caption{\ra{Fisher-KPP: ranks for rank-adaptive schemes (full observations, $\vartheta = 2 \cdot 10^{-8}$).}}
\label{fig:ranks}
\end{figure}

Finally, we explore the integration of rank-adaptivity in the assimilation scheme.
	Rank-adaptive methods are of particular interest in the filtering context, as it is common for the filtering distribution to eventually concentrate on a low-dimensional subspace (or even collapse onto a single point); using an increasingly smaller rank may prove to be sufficient to accurately evolve the filtering distribution. 
  \rb{On the other hand, upward rank-adaptivity may mitigate the effects of covariance under-approximation.}
	This experiment considers the case of full observations, and the truncation threshold is fixed to $\vartheta = 2 \cdot 10^{-8}$; 
  \rb{
   the system is initialised with $R_{\mathrm{initial}} = 2$ or $7$.}
	The rank evolutions are displayed in~\Cref{fig:ranks}.
  \rb{Across all methods and initial ranks, the ranks first increase and then decrease, suggesting it adequately captures the initial ensemble spread-out and subsequent filtering density collapse. 
  This also indicates the scheme is less susceptible to covariance under-approximation, as evidenced by the relative accuracy of the adaptive \dlr-\senkf~$R=2$ case, which roughly matches the corresponding fixed-rank $R=7$ case. 
  Furthermore, the scheme showcases a relative robustness to the choice of the initial rank, as the spike-then-decay behaviour is independent of that parameter.}
  The timings are displayed in~\Cref{tab:rank-adaptive-timings} and show a modest speed-up \rb{with respect to the fixed rank $R=7$}, which is not surprising given the relatively small initial rank. 
	Significantly higher speed-ups are expected for problems with higher initial ranks and finer discretisations. 

%

\begin{table}
	\centering
  \subfloat[\rb{$\norm{\E_{\rc{N}}[{\htheta}_{T}] - \theta^{\true}} / \norm{\theta^{\true}}$} \label{tab:rank-adaptive-relative-errors}]{
\begin{tabular}{@{}l|cc|cc@{}}
\toprule
 & \multicolumn{2}{c|}{$R_0 = 2$} & \multicolumn{2}{c}{$R_0 = 7$} \\
method & fixed & adaptive & fixed & adaptive \\ \midrule
\senkf & 0.264 & 0.030 & 0.026 & 0.030 \\
\venkf & 1.292 & 0.014 & 0.015 & 0.018 \\
\denkf & 0.180 & 0.009 & 0.006 & 0.009 \\  \bottomrule
\end{tabular}
  } \quad
	\subfloat[\rb{Timings [min]} \label{tab:rank-adaptive-timings}]{
\begin{tabular}{@{}l|cc|cc@{}}
\toprule
 & \multicolumn{2}{c|}{$R_0 = 2$} & \multicolumn{2}{c}{$R_0 = 7$} \\
method & fixed & adaptive & fixed & adaptive \\ \midrule
\senkf & 7.5 & 12.1 & 20.8 & 12.0 \\
\venkf & 7.2 & 13.3 & 20.9 & 12.8 \\
\denkf & 7.9 & 14.0 & 26.3 & 13.5 \\  \bottomrule
\end{tabular}
  }
  \caption{\ra{Fisher-KPP: fixed/adaptive ranks (single run, full observations, \\ truncation threshold $\vartheta = 2 \cdot 10^{-8}$).}
	}
	\label{tab:rank-adaptive}
\end{table}

\subsection{One-dimensional reduced blood flow model}
\label{sec:numexp-3}

\begin{figure}[htbp]
\centering
\includegraphics[width=0.8\linewidth]{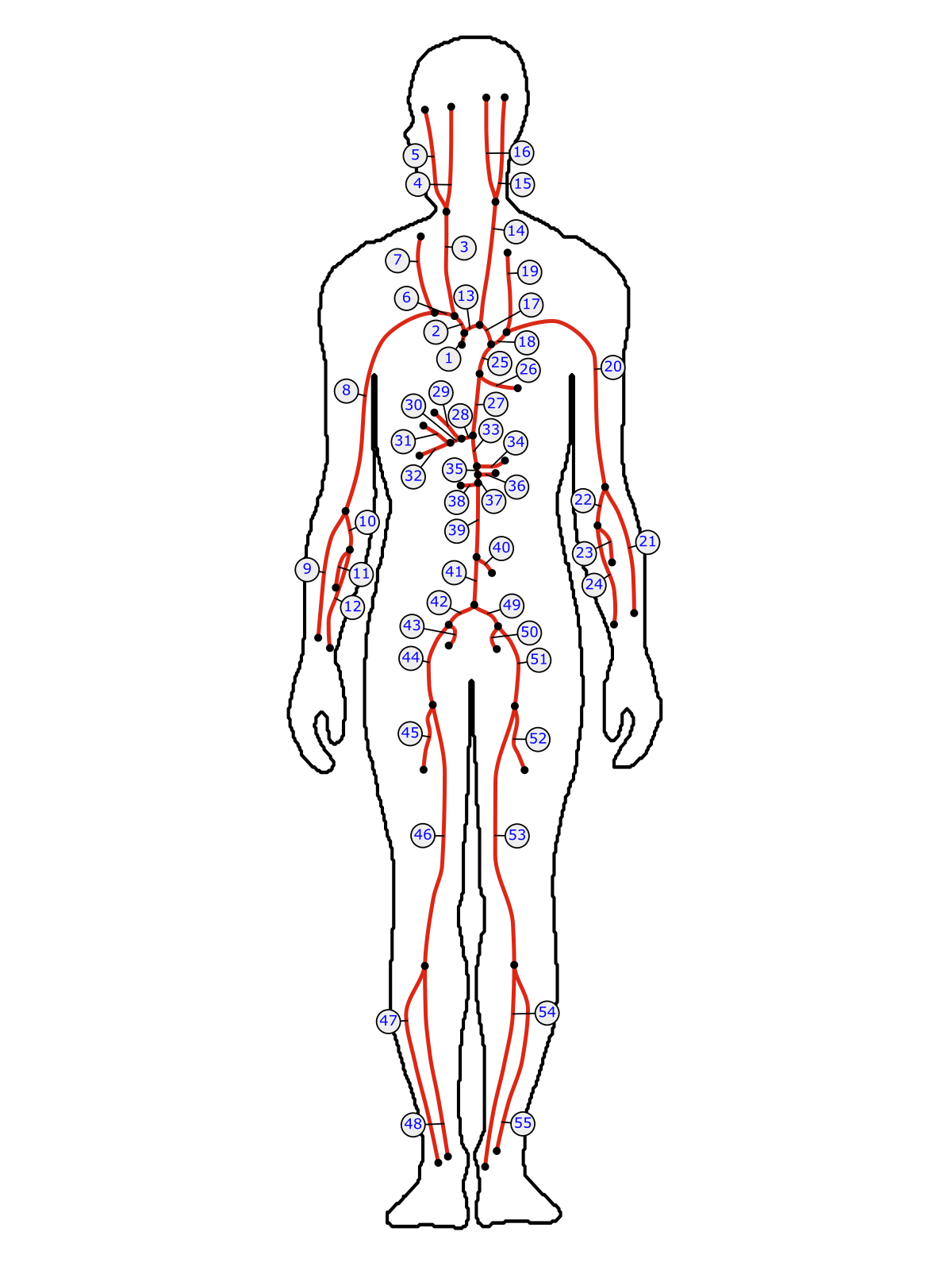}
\caption{\ra{Blood flow: arterial tree consisting of the 55 main human arteries}}
\label{fig:arterial-graph-exp3}
\end{figure}

The second numerical experiment, inspired by \cite{lombardi2014inverse}, aims to estimate the aortic stiffness of a one-dimensional reduced blood flow model using peripheral measurements. This model represents the blood flow in the 55 main human arteries (see Figure \ref{fig:arterial-graph-exp3}). It is derived from the three-dimensional incompressible Navier--Stokes equations by assuming that (1) vessel curvature is negligible, (2) vessels are axisymmetric, and (3) vessel properties are constant (see \cite{formaggia2010cardiovascular} for details). Under these assumptions, the mass and momentum equations reduce to
\begin{equation}
\label{eq:hemo_conservative}
\frac{\partial \mathbf{U}}{\partial t} + \frac{\partial \mathbf{F}(\mathbf{U})}{\partial x} = \mathbf{S}(\mathbf{U}),
\end{equation}
where
\begin{equation*}
\mathbf{U} = 
\begin{bmatrix}
A \\
u
\end{bmatrix}, \quad
\mathbf{F} = 
\begin{bmatrix}
Q \\
\frac{u^2}{2}+ \frac{p}{\rho}
\end{bmatrix}, \quad 
\mathbf{S} = 
\begin{bmatrix}
0 \\
-\nu \frac{u}{A}
\end{bmatrix},
\end{equation*}
with $x$ the axial coordinate along the vessel, $A$ the cross-sectional area, $u$ the mean fluid velocity, $Q=Au$ the flow rate, $p$ the fluid pressure, $\rho = \qty{1050}{\kg \per \cubic \m}$ the fluid density, and $\nu = \qty{9e-6}{\square \m \per \s}$ the viscous resistance. This system \eqref{eq:hemo_conservative} is closed by assuming that the vessel wall follows an elastic law, leading to the pressure--area relation:
\begin{equation}
p = p_{\textrm{ext}} + \beta (\sqrt{A}-\sqrt{A_0}),
\end{equation}
where $p_{\textrm{ext}} = \qty{11465.692}{\Pa}$ is the external pressure, $\beta$ is the elastic coefficient, and $A_0$ is the cross-sectional area at rest. The conservative system \eqref{eq:hemo_conservative} is hyperbolic because the Jacobian matrix
\begin{equation*}
\frac{\partial \mathbf{F}}{\partial \mathbf{U}} =
\begin{bmatrix}
u & A \\
\frac{\beta}{2\sqrt{A}} & u
\end{bmatrix}
\end{equation*}
has two real eigenvalues $\lambda_{1,2} = u \pm c$, and the characteristic variables are $W_{1,2} = u \pm 4c$, with $c = \sqrt{\frac{\beta}{2\rho}} A^{1/4}$ the wave speed propagation. In physiological conditions, the flow is subsonic with $c \gg u$. System \eqref{eq:hemo_conservative} is completed with suitable boundary conditions. At the inlet of the ascending aorta, we prescribe a periodic inflow of the form
\begin{equation}
\label{eq:proximal}
A(t) \;=\; a_0 + \sum_{k=1}^{4} \Bigl(a_k \cos(k \omega t) + b_k \sin(k \omega t)\Bigr),
\end{equation}
where $a_0 = \num{7.441e-4}$, $a_1 = \num{-2.809e-5}$, $b_1 = \num{4.699e-5}$, $a_2 = \num{-3.094e-5}$, $b_2 = \num{3.497e-6}$, $a_3 = \num{-5.753e-6}$, $b_3 = \num{-1.405e-5}$, $a_4 = \num{1.557e-6}$, $b_4 = \num{-2.932e-6}$, and $\omega = 2\pi/T$ with $T= \qty{0.8}{\s}$ the period of a cardiac cycle.
At the terminal nodes, non-reflecting boundary conditions are imposed:
\begin{equation}
\label{eq:non_reflection}
\frac{\rd W_2}{\rd t} = -\nu \frac{u}{A}.
\end{equation}
At the bifurcation nodes, the conservation of mass flux and the continuity of total pressure reduce to the conditions
\begin{equation}
\label{eq:bifurcation}
\begin{gathered}
Q_\textrm{in} = Q_{\textrm{out}_1} + Q_{\textrm{out}_2}, \\
\left(\rho\frac{u^2}{2}+p\right)_\textrm{in} = \left(\rho\frac{u^2}{2}+p\right)_{\textrm{out}_1} = \left(\rho\frac{u^2}{2}+p\right)_{\textrm{out}_2}.
\end{gathered}
\end{equation}
To ensure that the number of equations matches the number of unknowns at vessel boundaries and bifurcations, systems \eqref{eq:proximal}--\eqref{eq:bifurcation} are completed with the standard compatibility relations (see \cite{formaggia2010cardiovascular} for details). These relations enforce the correct coupling of characteristic variables and lead to the solution of local nonlinear systems of size $2 \times 2$ (at vessel boundaries) and $6 \times 6$ (at bifurcations). Physiological properties of the arterial tree are taken from \cite{liang2009biomechanical} and summarised in Table \ref{tab:arteries}. The resulting blood flow model is integrated in time using the forward Euler method and fixed time-step size $\delt = \num{5e-5}$. In space, a second-order finite volume method based on the MUSCL reconstruction \cite{TowardsTheUltVanLe1979} and local Lax-Friedrichs numerical flux \cite{NumericalMethoLevequ1992} is used with the uniform mesh size $\Delta x = \num{2e-3}$, leading to $d = \num{3542}$ DOFs.

The objective is now to estimate the elastic coefficient of the ascending aorta (n$^{\circ}$1) and aortic arch I (n$^{\circ}$13), i.e. $\theta^{\true} = \begin{bmatrix} \beta_1& \beta_{13} \end{bmatrix}$, while assuming all other parameters are known. For this purpose, we collect measurements of fluid velocity and cross-sectional area at the midpoint of three arteries\,--\,the left external carotid (n$^{\circ}$15), left radial (n$^{\circ}$21), and right femoral (n$^{\circ}$46)\,--\,yielding $k=6$ observations at each time-step. The observation noise covariance is
\begin{equation*}
\Gamma = 
\begin{bmatrix}
\sigma_A \bI_{3} & \\
& \sigma_u \bI_{3}
\end{bmatrix},
\end{equation*}
where $\sigma_A = 10^{-2}$ and $\sigma_u = 10^{-7}$, corresponding to about 5\% noise. The ensemble is initialised with ${\rc{N}} = 100$ particles sampled as $\htheta_0^{(i)} \sim \Ncal(\theta^{\mathrm{perturbed}}, \sigma_\theta^2 \bI_2)$, where $\theta^{\mathrm{perturbed}} \sim \Ncal(\theta^{\true}, \sigma_\theta^2 \bI_2)$ and $\sigma_\theta = 10^{5}$. In order to initialise the state associated to each particle, we start from the system at rest (i.e., $A = A_0$ and $u = 0$), and we perform two cardiac cycles (corresponding to $\qty{1.6}{\s}$ or \num{32000} time-steps) to reach a quasi-periodic flow. Then, we collect observations, generated using the true parameters, and assimilate them during five cardiac cycles (i.e., $\qty{4}{\s}$ or \num{80000} time-steps).

\begin{figure}[htbp]
\centering
\subfloat[$\beta_1$ (\unit{\kg \per \square \cm \per \square \s})]{\includegraphics[width=0.4\linewidth]{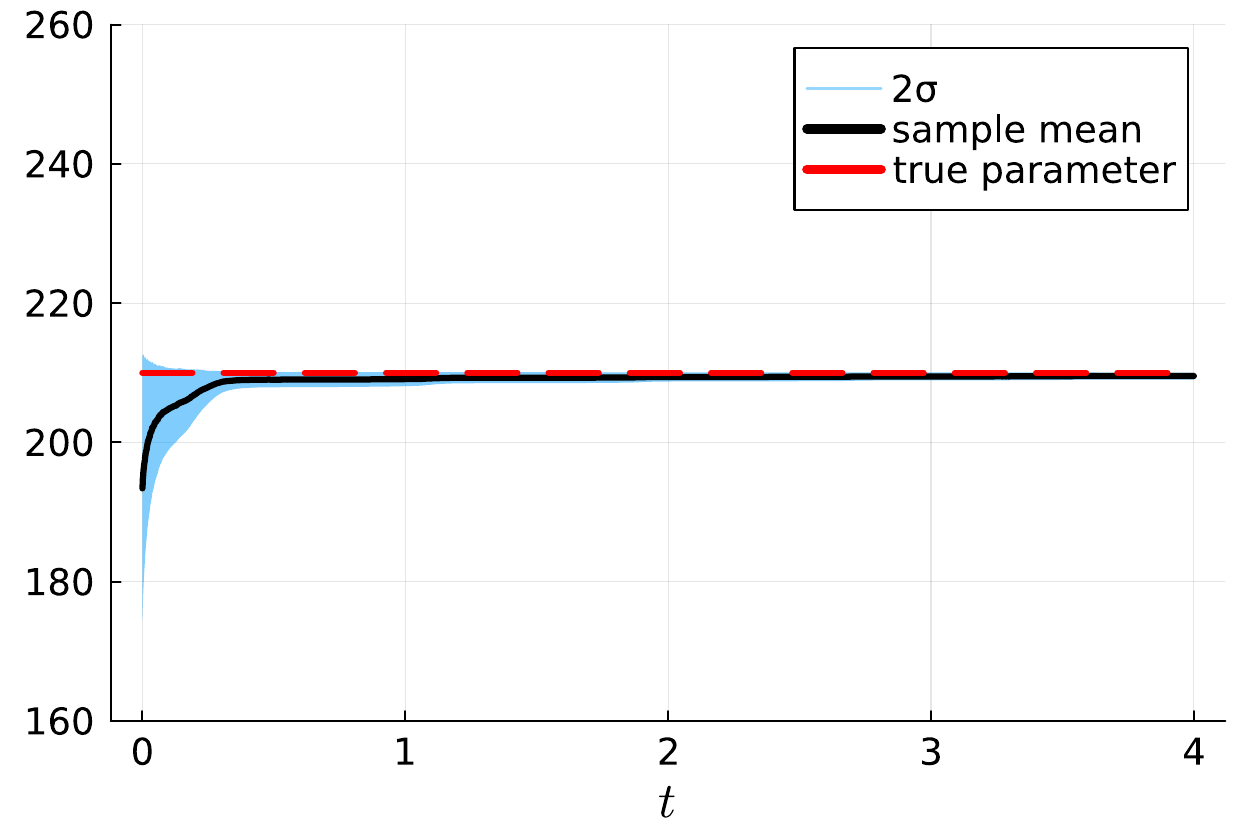}} \quad
\subfloat[$\beta_{13}$ (\unit{\kg \per \square \cm \per \square \s})]{\includegraphics[width=0.4\linewidth]{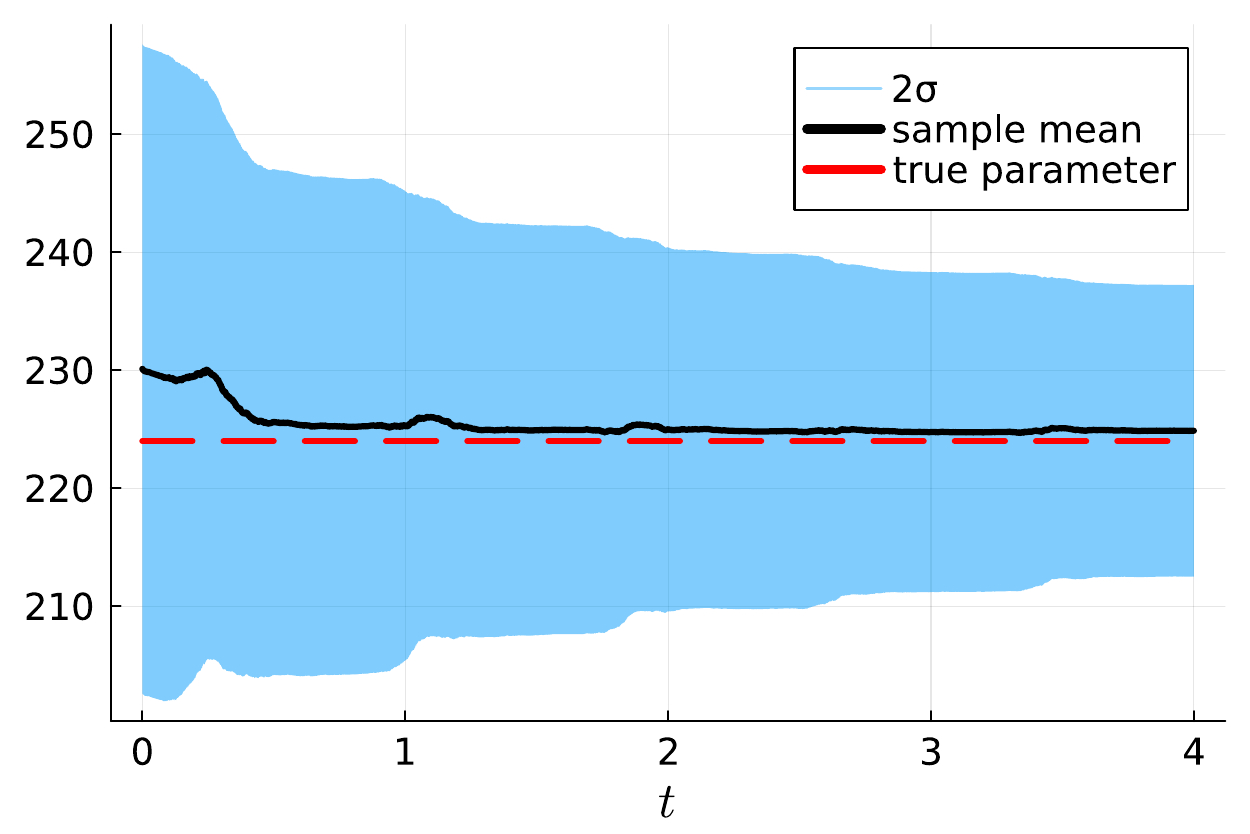}} 
\caption{\ra{Blood flow: Full Order Model \senkf~parameter identification of the elastic coefficient $\beta$ in the ascending aorta (n$^{\circ}$1) and aortic arch I (n$^{\circ}$13) with ${\rc{N}}=100$ particles}}
\label{fig:fom-parameter-identification-exp3}
\end{figure}

\begin{figure}[htbp]
\centering
\subfloat[$\beta_1$ (\unit{\kg \per \square \cm \per \square \s})]{\includegraphics[width=0.4\linewidth]{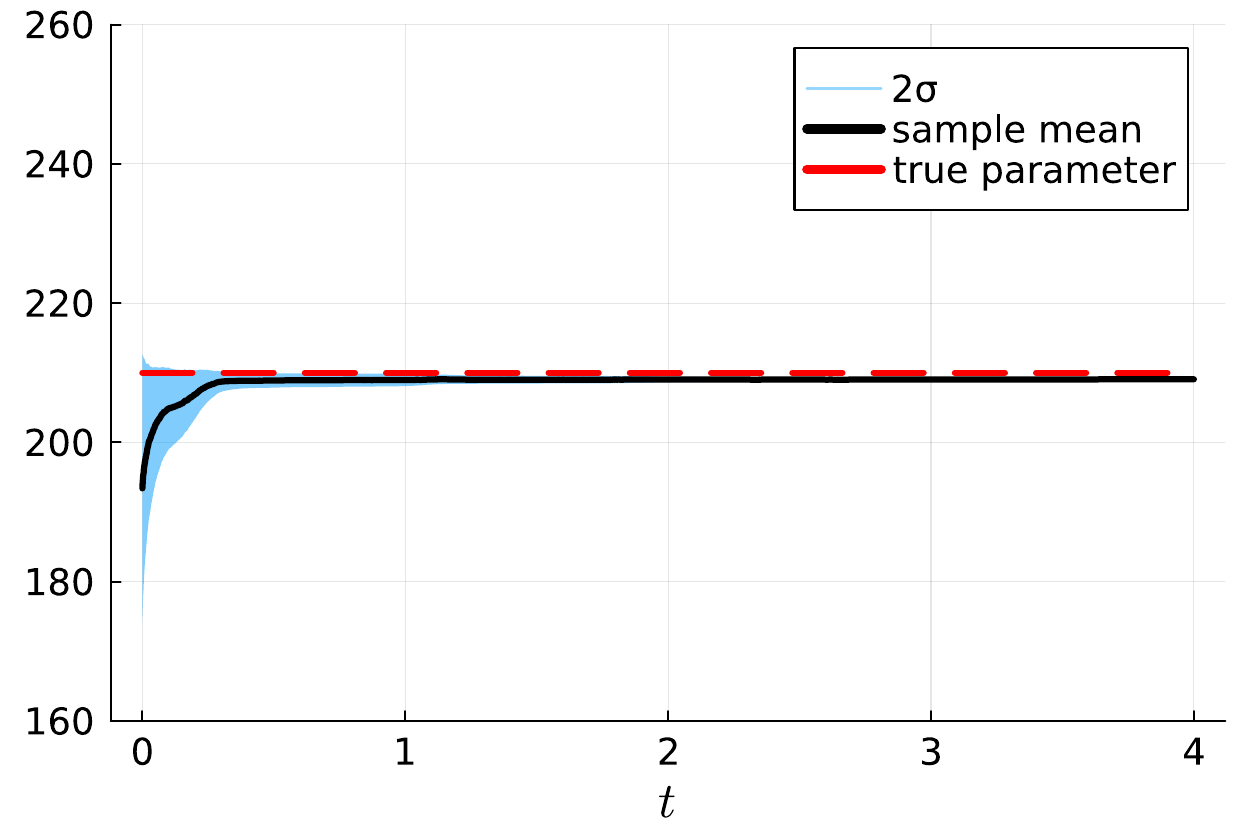}} \quad
\subfloat[$\beta_{13}$ (\unit{\kg \per \square \cm \per \square \s})]{\includegraphics[width=0.4\linewidth]{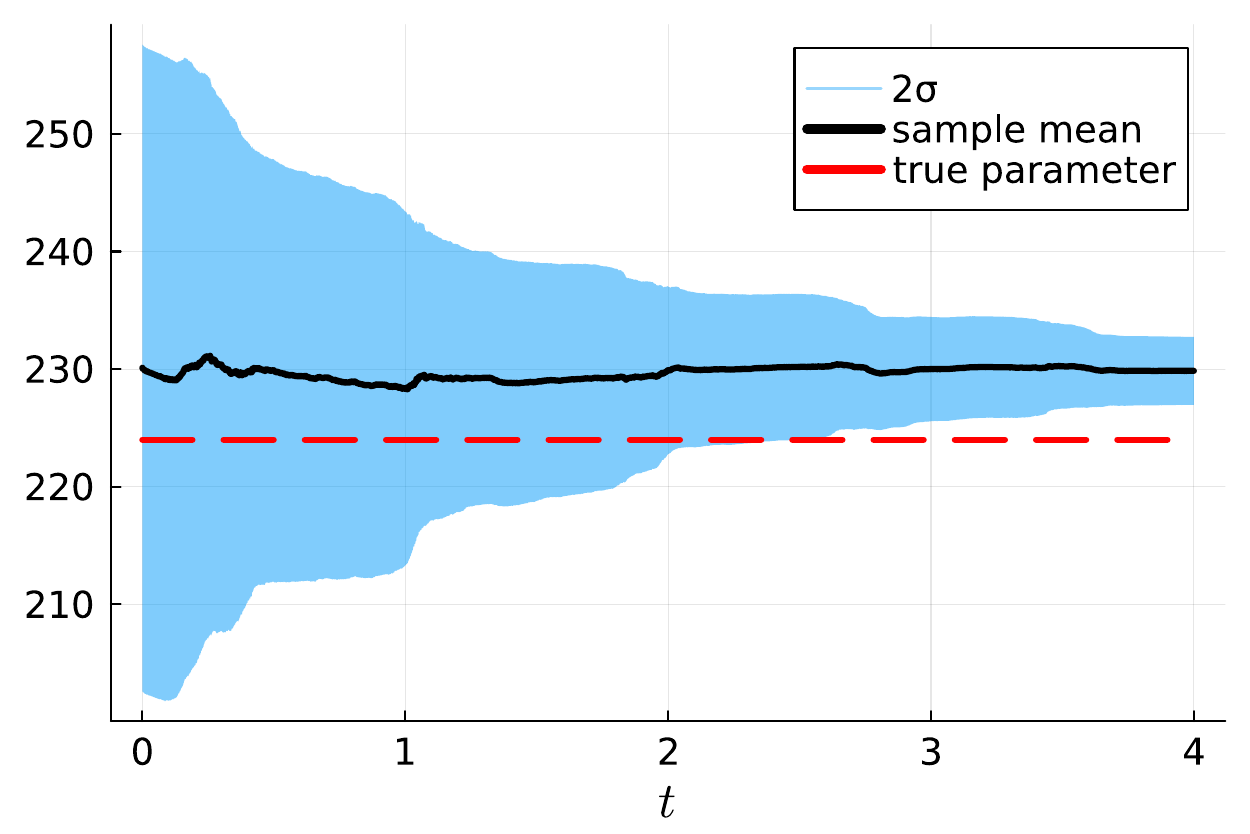}} 
\caption{\ra{Blood flow: hyper-reduced \dlr-\senkf~($R = 10$) parameter identification of the elastic coefficient $\beta$ in the ascending aorta (n$^{\circ}$1) and aortic arch I (n$^{\circ}$13) with ${\rc{N}}=100$ particles}}
\label{fig:dlr-10-parameter-identification-exp3}
\end{figure}

\begin{figure}[htbp]
\centering
\subfloat[$\beta_1$ (\unit{\kg \per \square \cm \per \square \s})]{\includegraphics[width=0.4\linewidth]{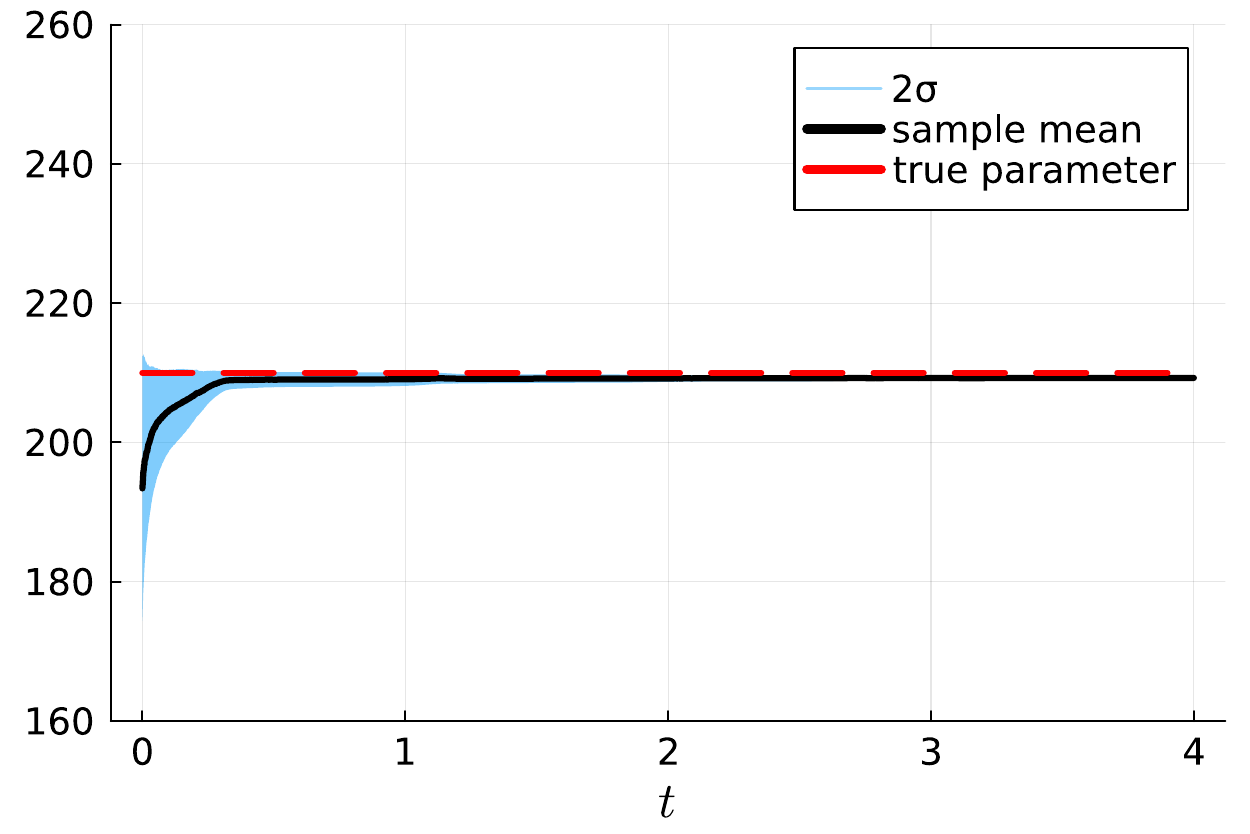}} \quad
\subfloat[$\beta_{13}$ (\unit{\kg \per \square \cm \per \square \s})]{\includegraphics[width=0.4\linewidth]{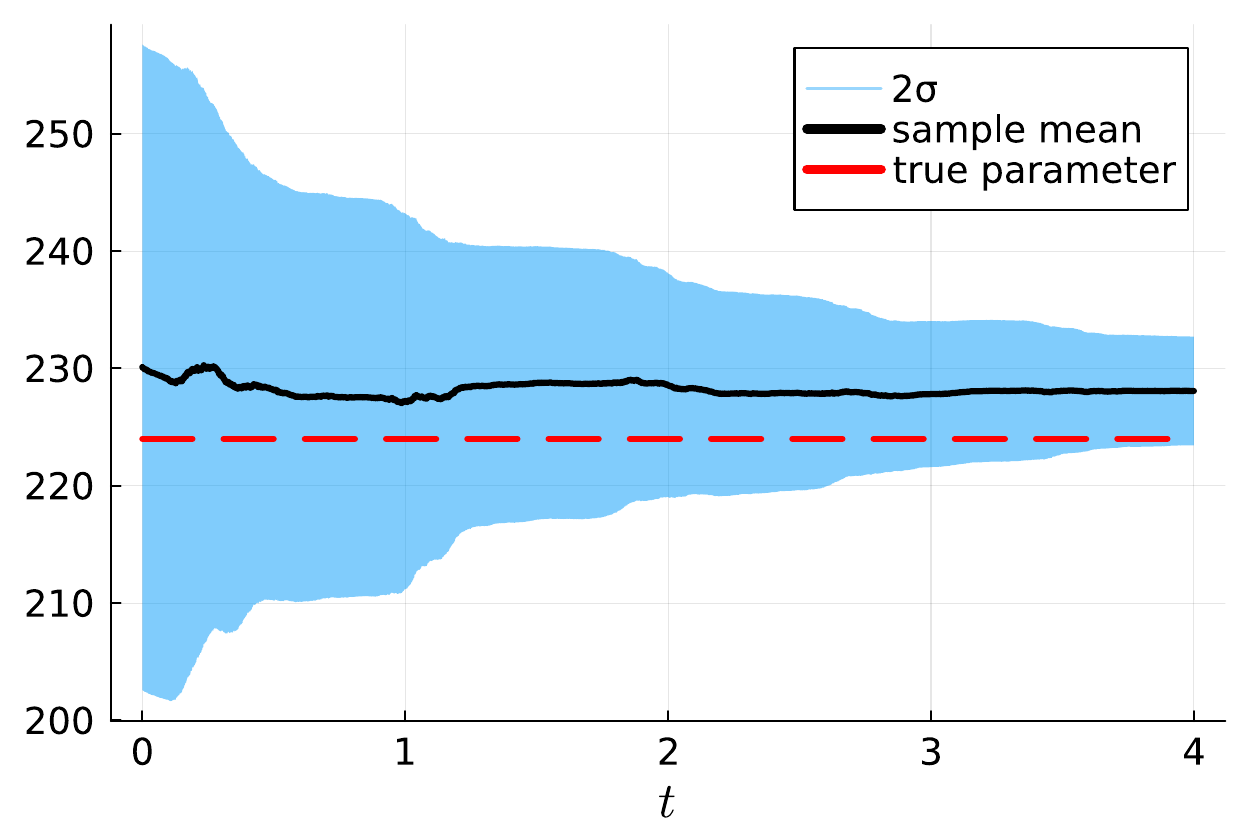}} 
\caption{\ra{Blood flow: hyper-reduced \dlr-\senkf~($R = 15$) parameter identification of the elastic coefficient $\beta$ in the ascending aorta (n$^{\circ}$1) and aortic arch I (n$^{\circ}$13) with ${\rc{N}}=100$ particles}}
\label{fig:dlr-15-parameter-identification-exp3}
\end{figure}

\begin{figure}[htbp]
\centering
\subfloat[$\beta_1$ (\unit{\kg \per \square \cm \per \square \s})]{\includegraphics[width=0.4\linewidth]{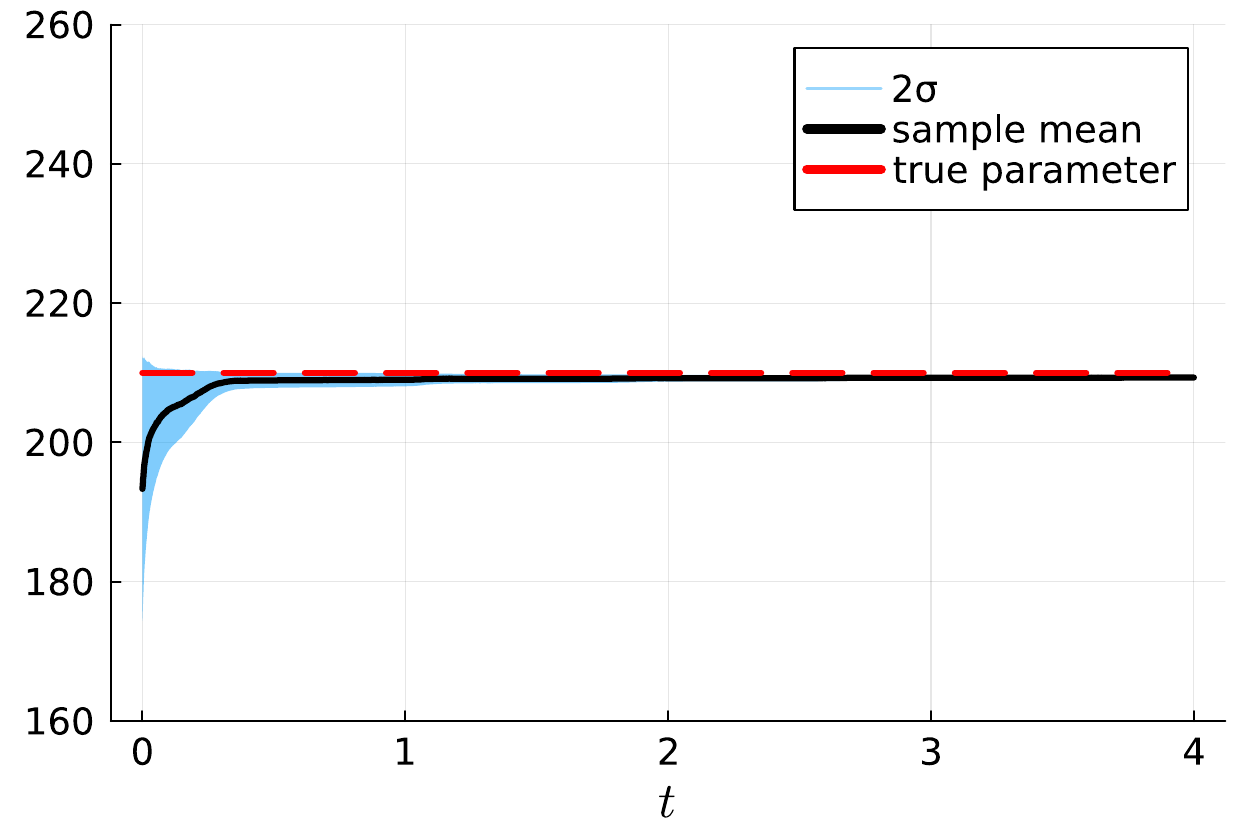}} \quad
\subfloat[$\beta_{13}$ (\unit{\kg \per \square \cm \per \square \s})]{\includegraphics[width=0.4\linewidth]{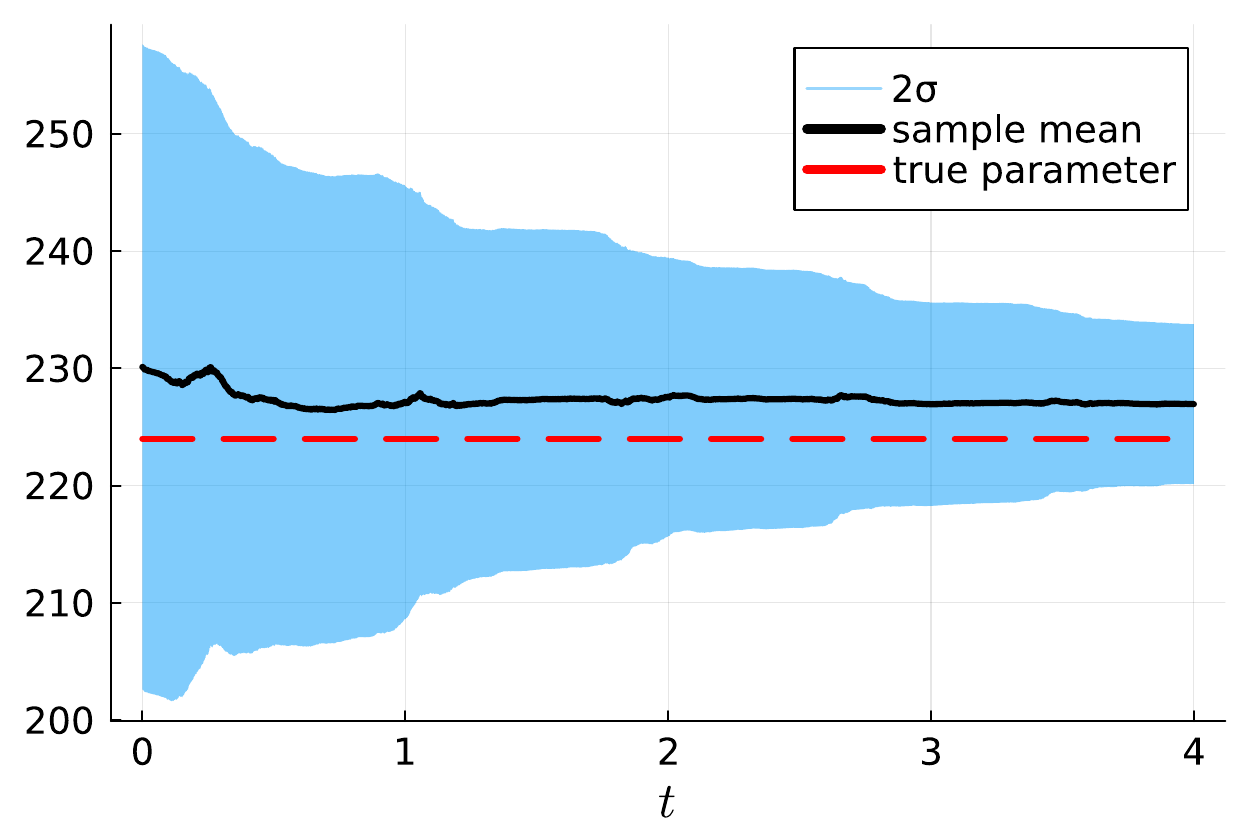}} 
\caption{\ra{Blood flow: hyper-reduced \dlr-\senkf~($R = 20$) parameter identification of the elastic coefficient $\beta$ in the ascending aorta (n$^{\circ}$1) and aortic arch I (n$^{\circ}$13) with ${\rc{N}}=100$ particles}}
\label{fig:dlr-20-parameter-identification-exp3}
\end{figure}

Figure~\ref{fig:fom-parameter-identification-exp3} shows the parameter identification achieved by the \senkf~(the results for the \venkf~and \denkf~are similar). The filter correctly identifies $\beta_1$, while there is more uncertainty for $\beta_{13}$ due to the limited number of measurements. Similarly, Figures~\ref{fig:dlr-10-parameter-identification-exp3}, \ref{fig:dlr-15-parameter-identification-exp3}, and~\ref{fig:dlr-20-parameter-identification-exp3} display the parameter identification achieved by the \dlr-\senkf~\rb{using the DEIM/CUR hyper-reduction} for ranks $R \in \{10,15,20\}$. We observe that $R=20$ correctly captures the filtering density over time, whereas for $R=10$, the estimate of $\beta_{13}$ is biased and the variance is underestimated.

Table~\ref{tab:performances-exp3-error} reports the relative error of the different methods. For $\beta_1$, the \dlr-\enkf~versions achieve the same level of accuracy as their full-order \enkf~counterparts even with small ranks. For $\beta_{13}$, the error tends to decrease as the rank increases. \rb{Moreover, the \dlr-\enkf~versions using the DEIM/CUR hyper-reduction achieve essentially the same accuracy as their non-hyper-reduced counterparts.}

Table~\ref{tab:performances-exp3-timings} reports the runtime of the different methods. \rb{The comparison with the \dlr-\enkf~without hyper-reduction shows that the low-rank formulation alone is not sufficient to reduce the computational cost in this example. Indeed, without hyper-reduction, the nonlinear terms are still evaluated at the full-order level, while the \dlr-\enkf~also requires additional algebraic operations, such as QR factorizations and SVD truncations in the BUG integrator. As a result, the non-hyper-reduced \dlr-\enkf~is slower than the full-order \enkf. In contrast, the \dlr-\enkf~using the DEIM/CUR hyper-reduction effectively reduces the runtime without significantly affecting the accuracy. Since the spatial discretisation used here is of moderate size, the relative gain of the hyper-reduction is expected to become more pronounced for finer discretisations, provided that the effective rank remains moderate. For a fixed spatial discretisation, the computational cost of the hyper-reduced \dlr-\enkf~increases with the rank, since larger ranks lead to larger reduced systems and more DEIM/CUR sampling indices.}

\begin{figure}[htbp]
\centering
\subfloat[$Q_1$ (\unit{\l \per \s})]{\includegraphics[width=0.4\linewidth]{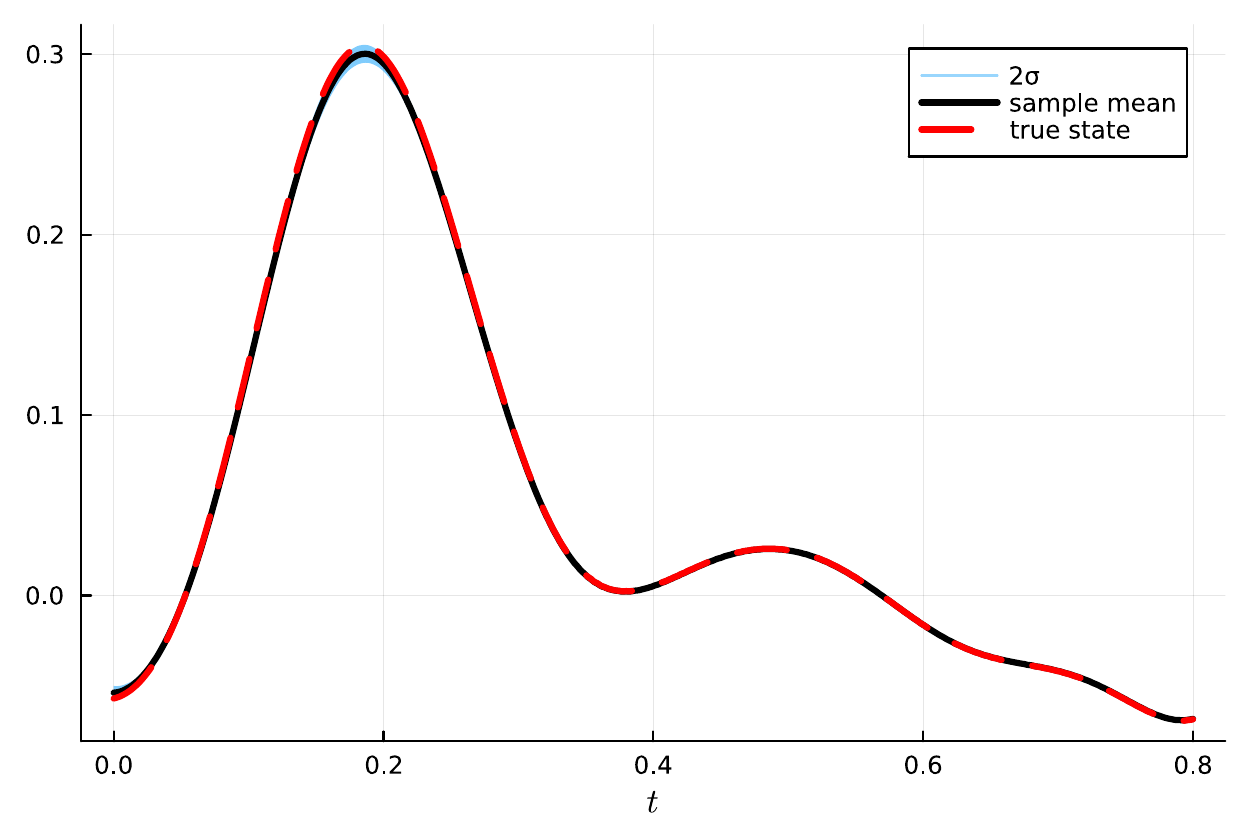}} \quad
\subfloat[$p_1$ (\unit{\mmHg})]{\includegraphics[width=0.4\linewidth]{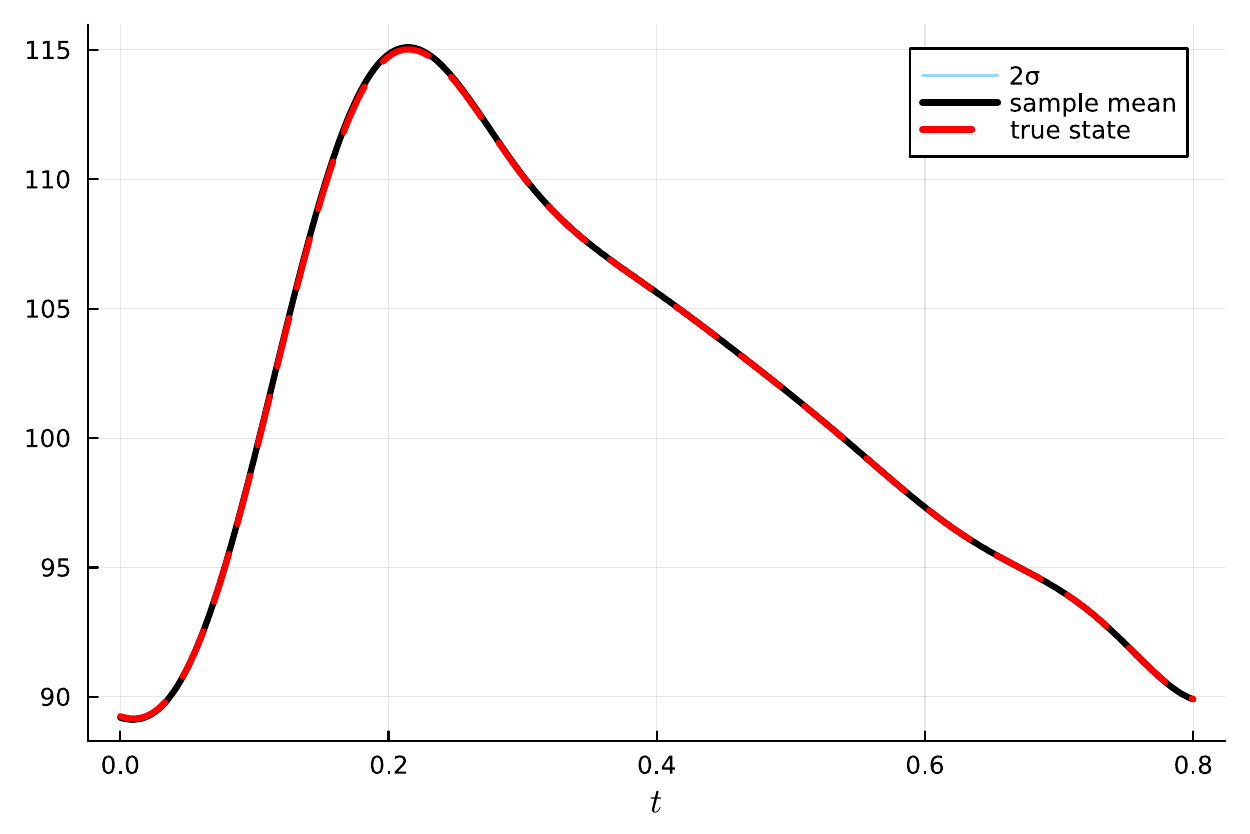}} \\
\subfloat[$Q_{13}$ (\unit{\l \per \s})]{\includegraphics[width=0.4\linewidth]{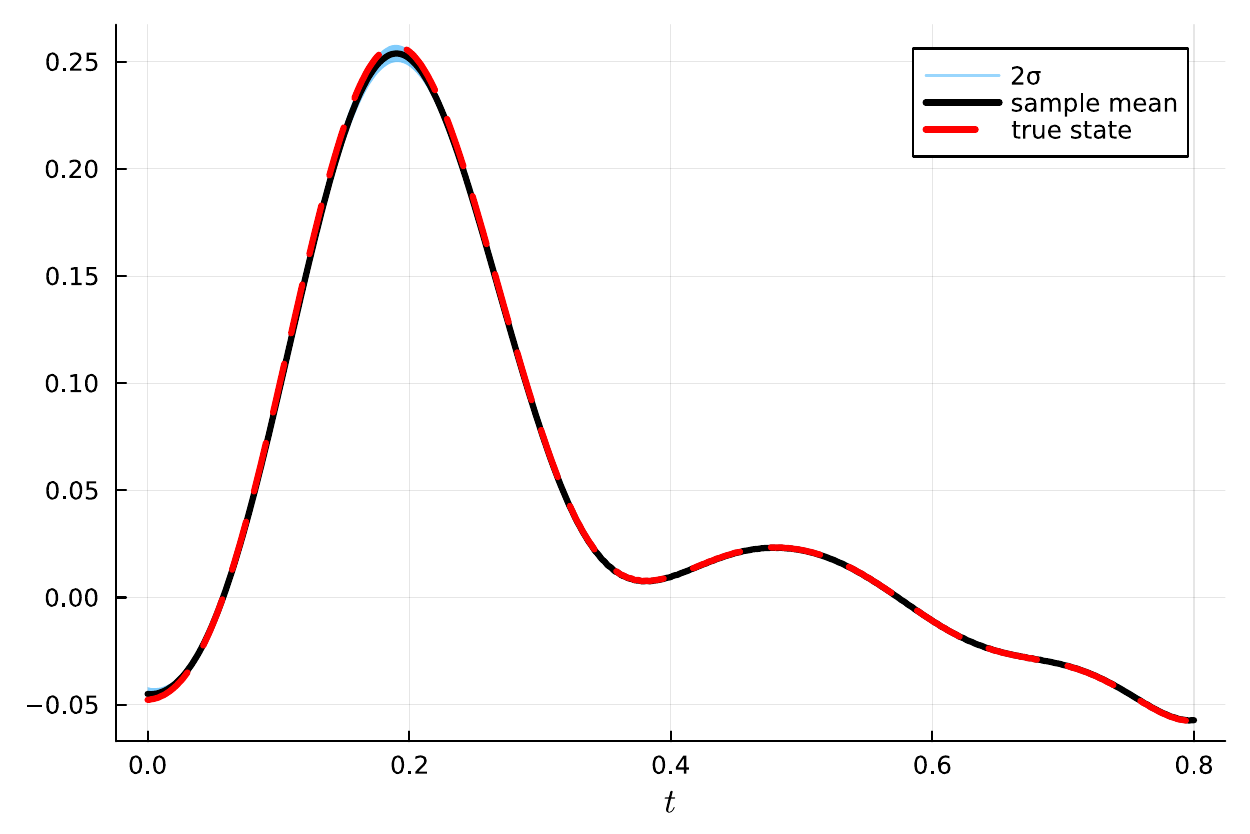}} \quad
\subfloat[$p_{13}$ (\unit{\mmHg})]{\includegraphics[width=0.4\linewidth]{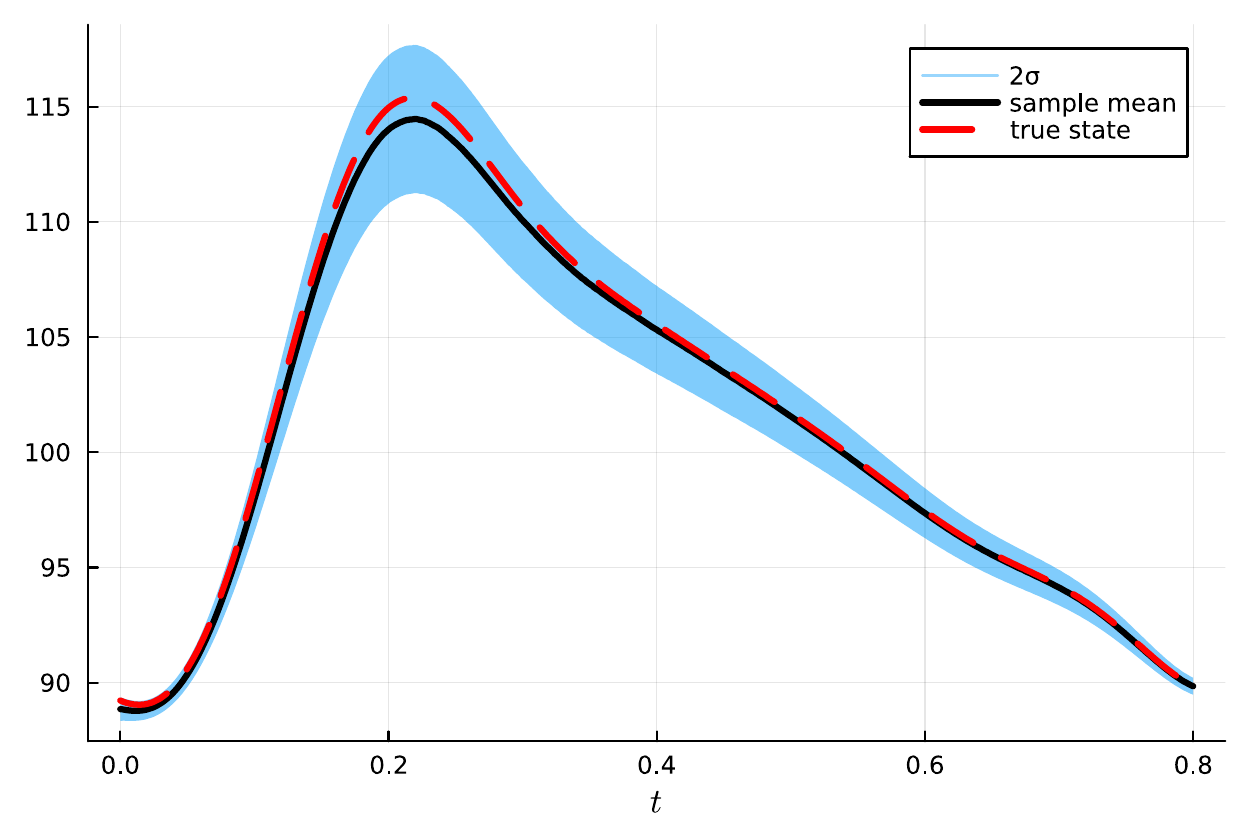}}
\caption{\ra{Blood flow: Full Order Model \senkf~state estimation of the flow rate $Q$ and pressure $p$ in the ascending aorta (n$^{\circ}$1) and aortic arch I (n$^{\circ}$13) with ${\rc{N}}=100$ particles}}
\label{fig:fom-state-estimation-exp3}
\end{figure}

\begin{figure}[htbp]
\centering
\subfloat[$Q_1$ (\unit{\l \per \s})]{\includegraphics[width=0.4\linewidth]{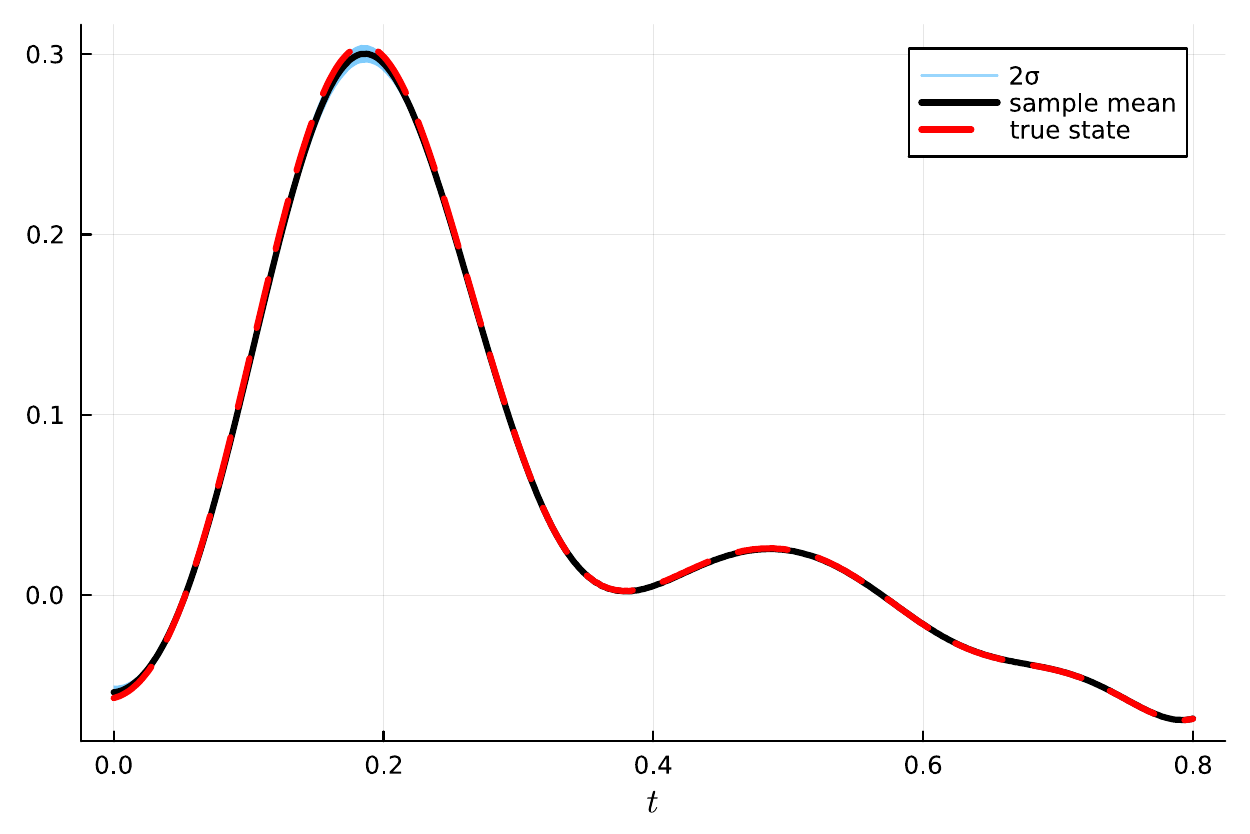}} \quad
\subfloat[$p_1$ (\unit{\mmHg})]{\includegraphics[width=0.4\linewidth]{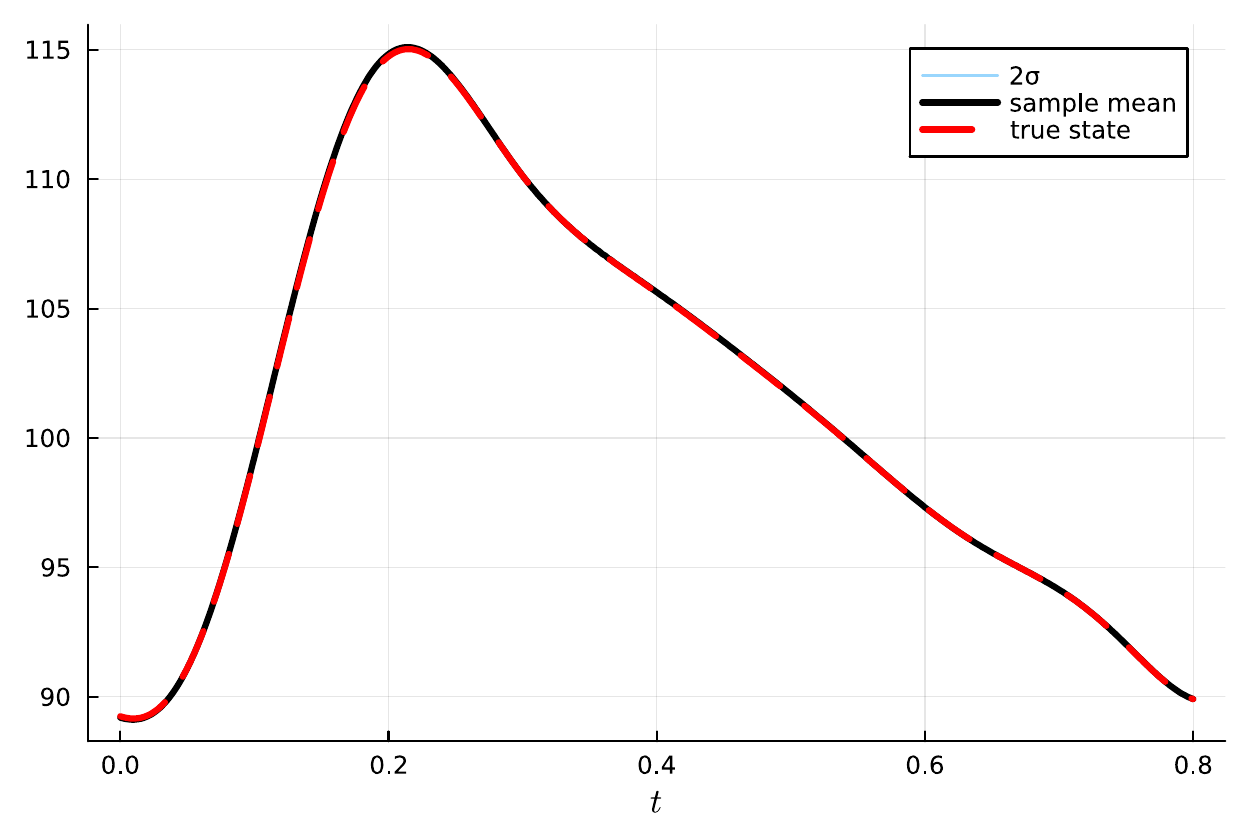}} \\
\subfloat[$Q_{13}$ (\unit{\l \per \s})]{\includegraphics[width=0.4\linewidth]{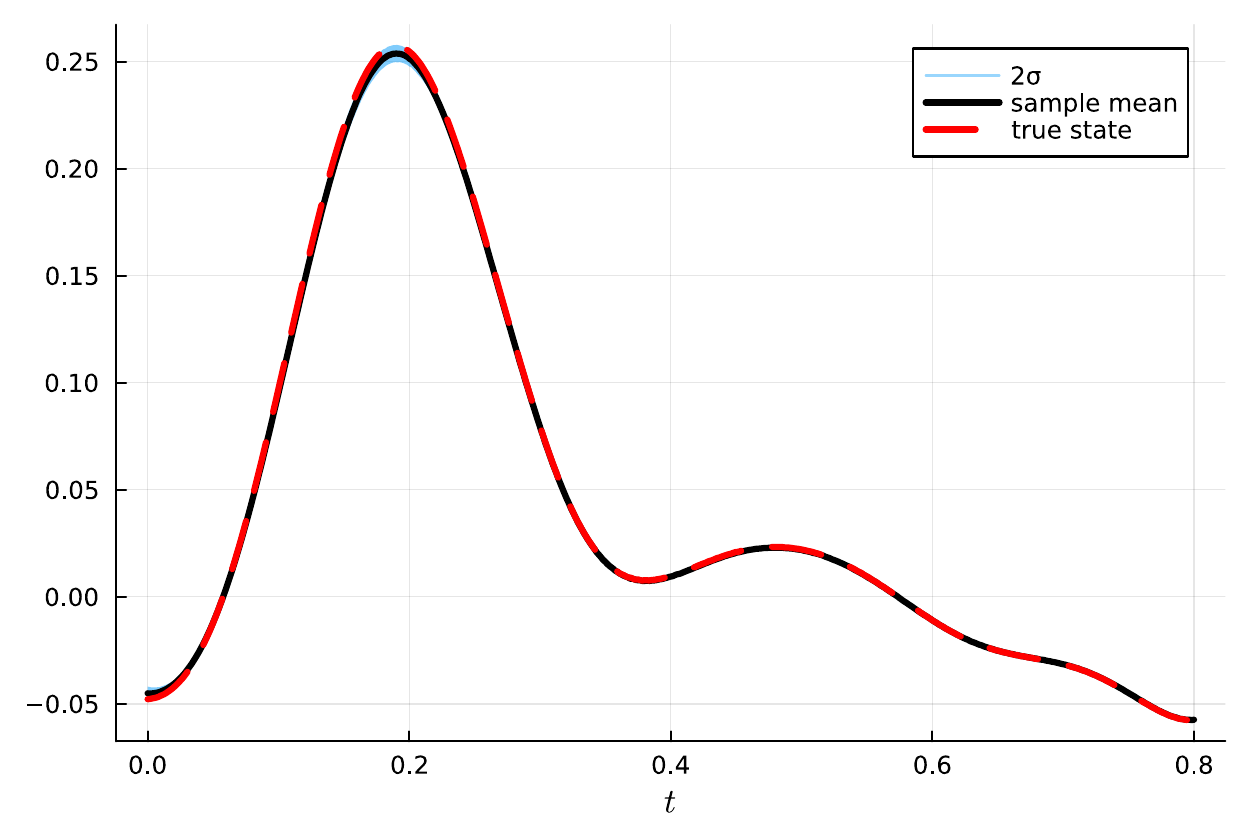}} \quad
\subfloat[$p_{13}$ (\unit{\mmHg})]{\includegraphics[width=0.4\linewidth]{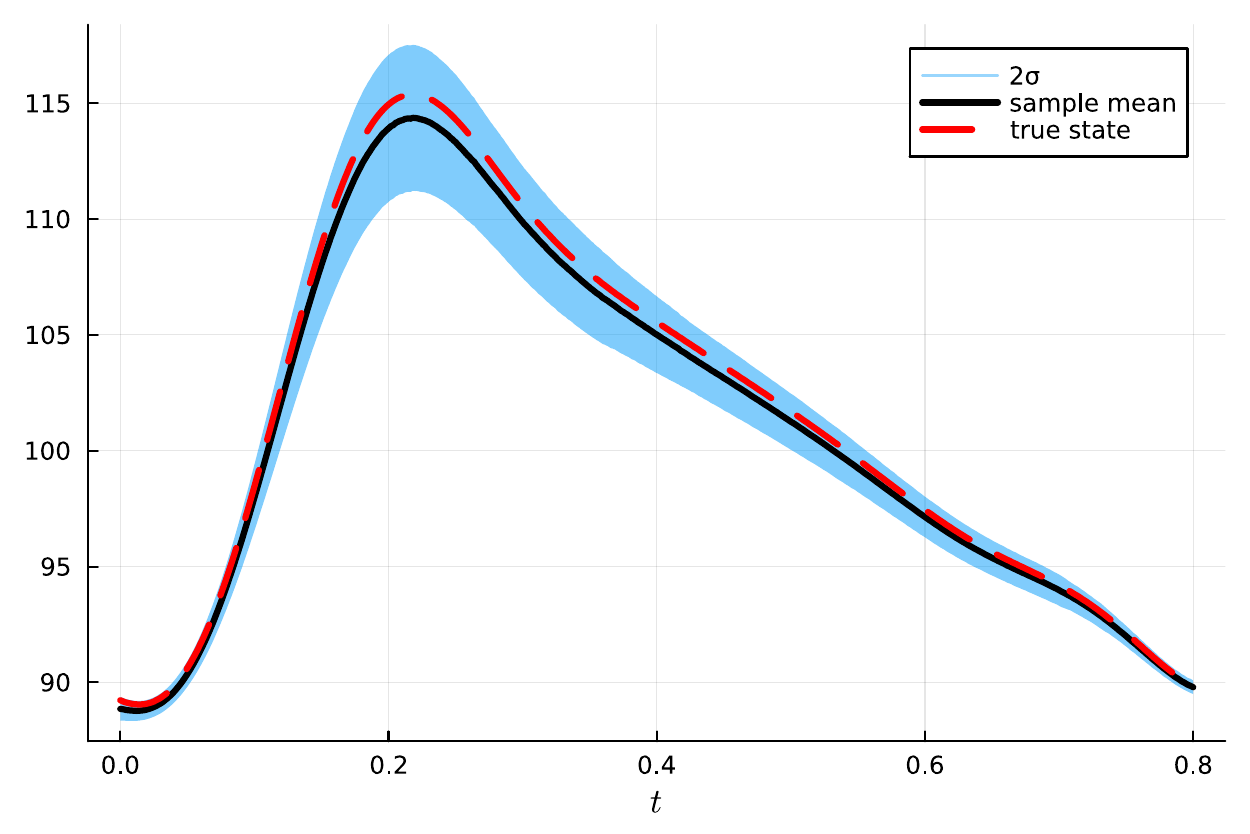}}
\caption{\ra{Blood flow: hyper-reduced \dlr-\senkf~($R = 10$) state estimation of the flow rate $Q$ and pressure $p$ in the ascending aorta (n$^{\circ}$1) and aortic arch I (n$^{\circ}$13) with ${\rc{N}}=100$ particles}}
\label{fig:dlr-state-estimation-exp3}
\end{figure}

Finally, Figures~\ref{fig:fom-state-estimation-exp3} and~\ref{fig:dlr-state-estimation-exp3} compare the reference flow rate and pressure in the ascending aorta (n$^{\circ}$1) and aortic arch I (n$^{\circ}$13) with the estimates obtained from the \senkf~and the \rb{hyper-reduced} \dlr-\senkf~with $R=10$. Both methods accurately recover the flow rate and pressure in each artery after one cardiac cycle.

\begin{table}
\setlength{\tabcolsep}{4pt}
\centering
\subfloat[$\beta_1$]{
\begin{tabular}{c|ccc|ccc|c}
\multicolumn{1}{l|}{} & \multicolumn{3}{c|}{\rb{DLR without hyper-reduction}} & \multicolumn{3}{c|}{\rb{DLR with DEIM/CUR}} & FOM \\ \midrule
\multicolumn{1}{l|}{} & \rb{$R = 10$} & \rb{$R =15$} & \rb{$R =20$} & \rb{$R = 10$} & \rb{$R =15$} & \rb{$R =20$} & \\ \midrule
\senkf & \rb{0.003} & \rb{0.003} & \rb{0.002} & \rb{0.004} & \rb{0.004} & \rb{0.003} & 0.002 \\ \midrule
\venkf & \rb{0.001} & \rb{0.001} & \rb{0.001} & \rb{0.001} & \rb{0.001} & \rb{0.001} & 0.001 \\ \midrule
\denkf & \rb{0.001} & \rb{0.001} & \rb{0.001} & \rb{0.001} & \rb{0.001} & \rb{0.001} & 0.001 \\ \bottomrule
\end{tabular}
} \quad\qquad
\subfloat[$\beta_{13}$]{
\begin{tabular}{c|ccc|ccc|c}
\multicolumn{1}{l|}{} & \multicolumn{3}{c|}{\rb{DLR without hyper-reduction}} & \multicolumn{3}{c|}{\rb{DLR with DEIM/CUR}} & FOM \\ \midrule
\multicolumn{1}{l|}{} & \rb{$R = 10$} & \rb{$R =15$} & \rb{$R =20$} & \rb{$R = 10$} & \rb{$R =15$} & $\rb{R =20}$ & \\ \midrule
\senkf & \rb{0.014} & \rb{0.010} & \rb{0.012} & \rb{0.034} & \rb{0.023} & \rb{0.017} & 0.008 \\ \midrule
\venkf & \rb{0.018} & \rb{0.013} & \rb{0.013} & \rb{0.022} & \rb{0.030} & \rb{0.027} & 0.010 \\ \midrule
\denkf & \rb{0.016} & \rb{0.016} & \rb{0.015} & \rb{0.029} & \rb{0.033} & \rb{0.022} & 0.011 \\ \bottomrule
\end{tabular}
}
\caption{\ra{Blood flow: relative error $\norm{\E_{\rc{N}}[{\htheta}_{T}] - \theta^{\true}} / \norm{\theta^{\true}}$ (averaged over 3 runs)}}
\label{tab:performances-exp3-error}
\end{table}

\begin{table}
\setlength{\tabcolsep}{4pt}
\centering
\subfloat{
\begin{tabular}{c|ccc|ccc|c}
\multicolumn{1}{l|}{} & \multicolumn{3}{c|}{\rb{DLR without hyper-reduction}} & \multicolumn{3}{c|}{\rb{DLR with DEIM/CUR}} & FOM \\ \midrule
\multicolumn{1}{l|}{} & \rb{$R = 10$} & \rb{$R =15$} & \rb{$R =20$} & $\rb{R = 10}$ & $\rb{R =15}$ & $\rb{R =20}$ & \\ \midrule
\senkf & \rb{81.7} & \rb{86.9} & \rb{90.9} & \rb{30.1} & \rb{47.0} & \rb{62.6} & \rb{73.9} \\ \midrule
\venkf & \rb{82.0} & \rb{86.6} & \rb{90.2} & \rb{30.1} & \rb{47.2} & \rb{62.6} & \rb{73.8} \\ \midrule
\denkf & \rb{81.9} & \rb{86.6} & \rb{90.6} & \rb{30.0} & \rb{47.1} & \rb{62.7} & \rb{74.3} \\ \bottomrule
\end{tabular}
}
\caption{\ra{Blood flow: timings [\unit{\minute}]}}
\label{tab:performances-exp3-timings}
\end{table}

\begin{table}[ht]
\centering
\setlength{\tabcolsep}{20pt}
\resizebox{\textwidth}{!}{
\begin{tabular}{r l c c c}
\hline
\multicolumn{1}{c}{No.} & \multicolumn{1}{c}{Artery} & Length (\unit{\m}) & $A_0$ (\unit{\square \m}) & $\beta$ (\unit{\kg \per \square \m \per \square \s}) \\
\hline
1  &  Ascending Aorta  &  \num{2.0e-2}  &  \num{6.812e-4}  &  \num{210e4} \\
2  &  Brachiocephalic  &  \num{3.4e-2}  &  \num{1.267e-4}  &  \num{652e4} \\
3  &  R. Carotid  &  \num{17.6e-2}  &  \num{0.466e-4}  &  \num{1079e4} \\
4  &  R. Internal Carotid  &  \num{17.6e-2}  &  \num{0.26e-4}  &  \num{2324e4}\\
5  &  R. External Carotid  &  \num{17.6e-2}  &  \num{0.126e-4}  &  \num{4042e4} \\
6  &  R. Subclavian I  &  \num{3.4e-2}  &  \num{0.544e-4}  &  \num{797e4} \\
7  &  R. Vertebral  &  \num{13.4e-2}  &  \num{0.126e-4}  &  \num{5505e4} \\
8  &  R. Subclavian II  &  \num{39.8e-2}  &  \num{0.319e-4}  &  \num{1077e4} \\
9  &  R. Radial  &  \num{22.0e-2}  &  \num{0.078e-4}  &  \num{7704e4} \\
10  &  R. Ulnar I  &  \num{6.6e-2}  &  \num{0.145e-4}  &  \num{4248e4} \\
11  &  R. Interosseous  &  \num{7.0e-2}  &  \num{0.031e-4}  &  \num{28722e4} \\
12  &  R. Ulnar II  &  \num{17.0e-2}  &  \num{0.115e-4}  &  \num{9714e4} \\
13  &  Aortic Arch I  &  \num{3.0e-2}  &  \num{5.992e-4}  &  \num{224e4} \\
14  &  L. Carotid  &  \num{20.8e-2}  &  \num{0.466e-4}  &  \num{1079e4} \\
15  &  L. External Carotid  &  \num{17.6e-2}  &  \num{0.126e-4}  &  \num{4042e4} \\
16  &  L. Internal Carotid  &  \num{17.6e-2}  &  \num{0.26e-4}  &  \num{2324e4} \\
17  &  Aortic Arch II  &  \num{4.0e-2}  &  \num{5.26e-4}  &  \num{239e4} \\
18  &  L. Subclavian I  &  \num{3.4e-2}  &  \num{0.544e-4}  &  \num{797e4} \\
19  &  L. Vertebral  &  \num{13.4e-2}  &  \num{0.126e-4}  &  \num{5505e4} \\
20  &  L. Subclavian II  &  \num{39.8e-2}  &  \num{0.319e-4}  &  \num{1077e4} \\
21  &  L. Radial  &  \num{22.0e-2}  &  \num{0.078e-4}  &  \num{7704e4} \\
22  &  L. Ulnar I  &  \num{6.6e-2}  &  \num{0.145e-4}  &  \num{4248e4} \\
23  &  L. Interosseous  &  \num{7.0e-2}  &  \num{0.031e-4}  &  \num{28722e4} \\
24  &  L. Ulnar II  &  \num{17.0e-2}  &  \num{0.115e-4}  &  \num{9714e4} \\
25  &  Thoracic Aorta I  &  \num{5.4e-2}  &  \num{4.412e-4}  &  \num{261e4} \\
26  &  Intercostals  &  \num{7.2e-2}  &  \num{0.283e-4}  &  \num{2008e4} \\
27  &  Thoracic Aorta II  &  \num{10.4e-2}  &  \num{3.294e-4}  &  \num{302e4} \\
28  &  Celiac I  &  \num{2.0e-2}  &  \num{0.332e-4}  &  \num{1252e4} \\
29  &  Hepatic  &  \num{6.4e-2}  &  \num{0.216e-4}  &  \num{2124e4} \\
30  &  Celiac II  &  \num{2.0e-2}  &  \num{0.238e-4}  &  \num{1843e4} \\
31  &  Gastric  &  \num{5.4e-2}  &  \num{0.126e-4}  &  \num{2426e4} \\
32  &  Splenic  &  \num{5.8e-2}  &  \num{0.083e-4}  &  \num{3801e4} \\
33  &  Abdominal Aorta I  &  \num{5.2e-2}  &  \num{2.438e-4}  &  \num{351e4} \\
34  &  Superior Mesenteric  &  \num{5.0e-2}  &  \num{0.442e-4}  &  \num{1052e4} \\
35  &  Abdominal Aorta II  &  \num{1.4e-2}  &  \num{2.143e-4}  &  \num{375e4} \\
36  &  L. Renal  &  \num{3.0e-2}  &  \num{0.238e-4}  &  \num{1577e4} \\
37  &  Abdominal Aorta III  &  \num{1.5e-2}  &  \num{2.026e-4}  &  \num{385e4} \\
38  &  R. Renal  &  \num{3.0e-2}  &  \num{0.238e-4}  &  \num{1577e4} \\
39  &  Abdominal Aorta IV  &  \num{12.4e-2}  &  \num{1.581e-4}  &  \num{436e4} \\
40  &  Inferior Mesenteric  &  \num{3.8e-2}  &  \num{0.11e-4}  &  \num{2468e4} \\
41  &  Abdominal Aorta V  &  \num{8.0e-2}  &  \num{1.088e-4}  &  \num{526e4} \\
42  &  R. Common Iliac  &  \num{5.8e-2}  &  \num{0.466e-4}  &  \num{931e4} \\
43  &  R. Internal Iliac  &  \num{4.4e-2}  &  \num{0.126e-4}  &  \num{6043e4} \\
44  &  R. External Iliac  &  \num{14.4e-2}  &  \num{0.367e-4}  &  \num{1722e4} \\
45  &  R. Deep Femoral  &  \num{11.2e-2}  &  \num{0.126e-4}  &  \num{3678e4} \\
46  &  R. Femoral  &  \num{44.2e-2}  &  \num{0.272e-4}  &  \num{2640e4} \\
47  &  R. Posterior Tibial  &  \num{34.4e-2}  &  \num{0.096e-4}  &  \num{9717e4} \\
48  &  R. Anterior Tibial  &  \num{32.2e-2}  &  \num{0.196e-4}  &  \num{4533e4} \\
49  &  L. Common Iliac  &  \num{5.8e-2}  &  \num{0.466e-4}  &  \num{931e4} \\
50  &  L. Internal Iliac  &  \num{4.4e-2}  &  \num{0.126e-4}  &  \num{6043e4} \\
51  &  L. External Iliac  &  \num{14.4e-2}  &  \num{0.367e-4}  &  \num{1722e4} \\
52  &  L. Deep Femoral  &  \num{11.2e-2}  &  \num{0.126e-4}  &  \num{3678e4} \\
53  &  L. Femoral  &  \num{44.2e-2}  &  \num{0.272e-4}  &  \num{2640e4} \\
54  &  L. Anterior Tibial  &  \num{34.4e-2}  &  \num{0.196e-4}  &  \num{4533e4} \\
55  &  L. Posterior Tibial  &  \num{32.2e-2}  &  \num{0.096e-4}  &  \num{9717e4} \\
\hline
\end{tabular}}
\caption{\ra{Blood flow: physiological properties of the arterial tree. The data are taken and adapted from \cite{liang2009biomechanical}}}
\label{tab:arteries}
\end{table}

\section{Conclusion}

\modb{In this work, we have proposed a Dynamical Low-Rank Ensemble Kalman Filter (\dlr-\enkf) for joint state--parameter estimation with time-continuous observations. The idea is to constrain the state ensemble to evolve in a low-dimensional subspace that adapts over time, while using an \enkf-type update to correct the augmented state. To this end, we propagate a low-rank approximation of the ensemble through the forecast and analysis steps: the forecast relies on a robust dynamical low-rank integrator inspired by the BUG scheme, while the analysis performs an \enkf-type update that exploits the resulting low-rank structure. To efficiently handle nonlinearities, we further incorporate a DEIM/CUR hyper-reduction method, which evaluates the nonlinear terms only on a small set of selected rows and columns.}

\modb{The proposed method has been evaluated on the Fisher--KPP equation and on a 1D hyperbolic blood flow model. The numerical experiments show that, for sufficiently large ranks $R$, the \dlr-\enkf~achieves an accuracy level comparable to that of its full-order \enkf~counterpart, while reducing the computational cost. The results also highlight the impact of rank truncation on parameter assimilation and identification: if the rank $R$ is too small, the truncation can lead to biased parameter estimates and an underestimation of the ensemble covariance, especially in partially observed regimes; whereas increasing $R$ tends to improve accuracy towards the full-order level.}


\rb{
In conclusion, the present work demonstrates the potential of DLR approximations for joint state-parameter estimation. 
Several directions remain open for future research. 
A first research direction concerns the use of more refined rank-adaptive strategies, such as mechanisms that are less constrained by the rank-doubling limitation of the BUG update, or that provide sharper information on the error components orthogonal to the augmented BUG subspace. 
In the filtering context, such strategies could also be adapted to the assimilation procedure, for example by enriching the low-rank subspace using information from the current observations. 
A further research direction is the integration of filter stabilisation techniques, such as additive variance inflation and covariance localisation, directly into the \dlr-\enkf~framework. Additive variance inflation is expected to be relatively straightforward to include, whereas covariance localisation is more challenging, since it may destroy the low-rank structure of the covariance and therefore reduce the computational gains.
}


\section*{Declaration of generative AI and AI-assisted technologies in the manuscript preparation process}
During the preparation of this work, the authors used large language models for text editing and polishing, and coding agents to refactor the code and execute numerical experiments.
The authors reviewed and edited the output as needed and take full responsibility for the content of the published article.

\section*{Declaration of competing interests}
The authors declare that they have no known competing financial interests or personal relationships that could have appeared to influence the work reported in this paper.

\section*{Acknowledgements} 
\modb{This work was funded by the Swiss National Science Foundation project ``Dynamical low rank methods for uncertainty quantification and data assimilation'' (n. 200518).}

\section*{Data availability}
The research code to reproduce the results of this paper is available at \url{https://zenodo.org/records/18655790}.

\bibliographystyle{elsarticle-num}
\bibliography{bibliography}

\end{document}